\documentclass[10pt,a4paper,reqno,oneside,authoryear]{elsarticle.modified}

\usepackage{amsmath}
\usepackage{amssymb}
\usepackage{mathscinet}
\usepackage{accents}
\usepackage{mathrsfs}
\usepackage{enumitem}
\usepackage{ifthen}
\usepackage{pifont}
\usepackage{bbding}
\usepackage{textcomp}
\usepackage{xspace}
\usepackage[pagebackref]{hyperref}
\usepackage{comment}
\usepackage[T1]{fontenc}
\usepackage{fixme} 
\usepackage{calc}

\newtheorem{lem}{Lemma}

\newtheorem{lemclaim}{Claim}[lem]
\newtheorem{lemsubclaim}{Subclaim}[lemclaim]

\newtheorem{prop}[lem]{Proposition} 
\newtheorem{claim}{Claim}
\newtheorem{thm}{Theorem}
\newtheorem{question}{Question}
\newtheorem{cor}[lem]{Corollary}

\newdefinition{defn}[lem]{Definition }
\newdefinition{example}[lem]{Example }

\newproof{remark}[lem]{Remark} 
\newproof{proof}{Proof}
\newenvironment{pf}{\begin{proof}}{\end{proof}}
\newtheorem{termgy}[lem]{Terminology}
\newtheorem{notn}[lem]{Notation}

\makeatletter

\newenvironment{enumeq}{
 \begin{enumerate}[label=(\arabic*), ref=\arabic*, widest=10]
     \setcounter{enumi}{\value{equation}}}{\setcounter{equation}{\value{enumi}}\end{enumerate}}

\newlength{\label@width}
\newlength{\label@sep}
\newlength{\left@margin}
\newcounter{@temp}

\newenvironment{principle}[1]{
\setlength{\label@sep}{\labelsep}
\begin{list}{#1}{\setlength{\rightmargin}{0pt}
\settowidth{\labelwidth}{#1}\settowidth{\labelsep}{\quad}\setlength{\leftmargin}{\labelwidth+\labelsep}}
\item
\begin{minipage}[t]{\linewidth}
}{
\end{minipage}
\end{list}
}

\newcommand{\ch}{\mathrm{CH}}
\newcommand{\gch}{\mathrm{GCH}}

\newcommand{\zfc}{\mathrm{ZFC}}

\newcommand{\pfa}{\mathrm{PFA}}
\newcommand{\mm}{\mathrm{MM}}

\newcommand{\pstar}{(*)}
\newcommand{\pstarc}{(\gdot_c)}

\newcommand{\dom}{\operatorname{dom}}

\newcommand{\supp}{\operatorname{supp}}

\newcommand{\game}{\Game}

\newcommand{\rank}{\operatorname{rank}}
\newcommand{\length}{\operatorname{len}}

\newcommand{\downcl}[1]{{\downarrow\hspace{-2.5pt}#1}}
\newcommand{\upcl}[1]{{\uparrow\hspace{-2.5pt}#1}}

\newcommand{\gdot}{\star}
\newcommand{\power}{\P}

\newcommand{\inv}{^{-1}}

\newcommand{\powfin}[1]{[#1]^{<\aleph_0}}
\newcommand{\powcnt}[1]{[#1]^{\leq\aleph_0}}
\newcommand{\powcif}[1]{[#1]^{\aleph_0}}

\newcommand{\otp}{\operatorname{otp}}

\newcommand{\hit}{\operatorname{ht}}

\newcommand\bigext{^\frown}

\newcommand{\lands}{\mathrel{\land}}

\renewcommand{\div}{\mathbin{/}}

\newcommand{\Gen}{\operatorname{Gen}}
\newcommand{\Genc}{\operatorname{Gen}^+}

\newcommand{\id}{\operatornamewithlimits{id}}

\newcommand{\ifm}{\bigwedge}
\newcommand{\spr}{\bigvee}

\renewcommand{\restriction}{\mathbin{\upharpoonright}}

\newcommand{\@@eq}[2]{=_{#1}^{#2}}
\newcommand{\@@neq}[2]{\ne_{#1}^{#2}}
\newcommand{\@@l}[2]{<_{#1}^{#2}}
\newcommand{\@@nl}[2]{\nless_{#1}^{#2}}
\newcommand{\@@le}[2]{\le_{#1}^{#2}}
\newcommand{\@@nle}[2]{\nleq_{#1}^{#2}}
\newcommand{\@@equiv}[2]{\equiv_{#1}^{#2}}

\newcommand{\@in}[1]{\in_#1}
\newcommand{\@eq}[1]{=_#1}
\newcommand{\@le}[1]{\le_{#1}}
\newcommand{\@l}[1]{<_#1}
\newcommand{\@nle}[1]{\nleq_#1}
\newcommand{\@ge}[1]{\ge_#1}

\newcommand{\u@le}[1]{\le^#1}
\newcommand{\u@l}[1]{<^#1}
\newcommand{\u@eq}[1]{=^#1}

\newcommand{\@lefnt}[1]{\@@le#1*}
\newcommand{\@nlefnt}[1]{\@@nle#1*}
\newcommand{\@eqfnt}[1]{\@@eq#1*}

\newcommand{\subseteqfnt}{\subseteq^*}

\newcommand{\supseteqfnt}{\mathrel{\supseteq^*}}

\newcommand{\len}[1]{\@le{#1}}
\newcommand{\lln}[1]{\@l#1}
\newcommand{\gen}[1]{\@ge{{#1}}}
\newcommand{\altle}{\lesssim}

\newcommand{\str}{\mathrm{str}}
\newcommand{\sep}{\mathrm{sep}}
\newcommand{\asym}{\mathrm{asym}}
\newcommand{\trn}{\mathrm{tran}}
\newcommand{\negstr}{-\str}

\newcommand{\ggen}[1]{\gg_{#1}}

\newcommand{\strlen}[1]{\@@le{#1}\str}
\newcommand{\nstrlen}[1]{\@@nle{#1}\str}
\newcommand{\strln}[1]{\@@l{#1}\str}
\newcommand{\nstrln}[1]{\@@nl{#1}\str}
\newcommand{\strequiv}[1]{\@@equiv{#1}\str}

\newcommand{\lenn}[2]{\@@le{#1}{#2}}
\newcommand{\nlenn}[2]{\@@nle{#1}{#2}}
\newcommand{\lnn}[2]{\@@l{#1}{#2}}
\newcommand{\nlnn}[2]{\@@nl{#1}{#2}}

\newcommand{\eqnn}[2]{\@@eq{#1}{#2}}
\newcommand{\neqnn}[2]{\@@neq{#1}{#2}}
\newcommand{\asymle}{\u@le\asym}
\newcommand{\asyml}{\u@l\asym}

\newcommand{\strtrnlen}[1]{\@@le{#1}{\str,\trn}}
\newcommand{\strtrnln}[1]{\@@l{#1}{\str,\trn}}
\newcommand{\trnle}{\u@le{\trn}}
\newcommand{\trnl}{\u@l{\trn}}
\newcommand{\trneq}{\u@eq{\trn}}
\newcommand{\negstrlen}[1]{\@@le{#1}\negstr}
\newcommand{\nnegstrlen}[1]{\@@nle{#1}\negstr}
\newcommand{\nnegstrln}[1]{\@@nl{#1}\negstr}
\newcommand{\negstreq}[1]{\@@eq{#1}\negstr}

\renewcommand{\aa}{\mathrm{aa}}
\newcommand{\eqaa}{\@eq\aa}
\newcommand{\inaa}{\@in\aa}

\newcommand{\leaa}{\@le\aa}
\newcommand{\nleaa}{\@nle\aa}
\newcommand{\geaa}{\@ge\aa}
\newcommand{\lefntaa}{\@lefnt\aa}
\newcommand{\nlefntaa}{\@nlefnt\aa}
\newcommand{\eqfntaa}{\@eqfnt\aa}

\newcommand{\naa}{\mathrm{naa}}
\newcommand{\eqnaa}{\@eq\naa}
\newcommand{\eqfntnaa}{\@eqfnt\naa}
\newcommand{\lenaa}{\@le\naa}
\newcommand{\lefntnaa}{\@lefnt\naa}

\newcommand{\eqae}{\@eq\ae}
\newcommand{\leae}{\@le\ae}

\newcommand{\forces}{\mskip 5mu plus 5mu\|\hspace{-2.5pt}{\textstyle \frac{\hspace{4.5pt}}{\hspace{4.5pt}}}\mskip 5mu plus 5mu}

\newcommand{\Iff}{\espc\mathrm{iff}\espc}

\newcommand{\impls}{\espc\mathrm{implies}\espc}

\renewcommand{\And}{\espc\mathrm{and}\espc}

\renewcommand{\and}{\hespc\mathrm{and}\hespc}

\renewcommand{\iff}{\leftrightarrow}

\newcommand{\p}{\mathfrak p}

\newcommand{\cof}{\operatorname{cof}}

\newcommand{\one}{\mathbf{1}}

\newcommand{\ord}{\mathrm{On}}

\newcommand{\pmax}{\mathbb P_{\mathrm{max}}}

\DeclareFontFamily{U}{cmsy}{}
\DeclareFontShape{U}{cmsy}{m}{n}{<12> sfixed * [10] cmsy10 
<10> <9> <8> <7> <6> <5> sfixed * [10] cmsy10}{}
\DeclareSymbolFont{customtwo}{U}{cmsy}{m}{n} 
\DeclareMathSymbol{\sctn}{\mathord}{customtwo}{"78}

\DeclareFontFamily{U}{cmmi}{}
\DeclareFontShape{U}{cmmi}{m}{n}{<20> sfixed * [11] cmmib10 <12> sfixed * [10]
cmmi10 <10> <9> <8> sfixed * [6] cmmi6 <5> <6> <7> sfixed * [5] cmmi5}{}
\DeclareSymbolFont{custom}{U}{cmmi}{m}{n}
\DeclareMathSymbol{\rharpoon}{\mathord}{custom}{"2A}
\newlength{\widt}
\newlength{\widttwo}
\newlength{\hgt}

\newcommand{\oone}{{\omega_1}}

\newcommand{\<}{\langle}
\renewcommand{\>}{\rangle}

\newcommand{\espc}{\quad}
\newlength{\@@hespc}
\newcommand{\hespc}{\hspace{\@@hespc}}
\newcommand{\overbar}[1]{\,\overline{\!{#1}}}

\newcommand{\overbarg}[1]{\hsp\overline{\hspb{#1}}}
\newcommand{\Th}{{^{\mathrm{th}}}}

\renewcommand{\ae}{\mathrm{ae}}

\newcommand{\ulc}{\ulcorner}
\newcommand{\urc}{\urcorner}

\newcommand{\tu}{\textup}

\newcommand{\A}{\mathcal A}
\newcommand{\B}{\mathcal B}
\newcommand{\C}{\mathcal C}
\newcommand{\D}{\mathcal D}
\newcommand{\E}{\mathcal E}
\newcommand{\F}{\mathcal F}

\newcommand{\ideal}{\I}
\newcommand{\I}{\mathcal I}
\newcommand{\J}{\mathcal J}
\newcommand{\K}{\mathcal K}
\renewcommand{\L}{\mathcal L} 
\newcommand{\G}{\mathcal G}
\renewcommand{\H}{\mathcal H}
\newcommand{\M}{\mathcal M}
\renewcommand{\O}{\mathcal O}
\renewcommand{\P}{\Pcal}
\newcommand{\Pcal}{\mathcal P}
\newcommand{\Q}{\mathcal Q}
\newcommand{\R}{\mathcal R}
\newcommand{\Section}{{\char'237}}
\renewcommand{\S}{\mathcal S}
\newcommand{\T}{\mathcal T}

\newcommand{\X}{\mathcal X}
\newcommand{\Y}{\mathcal Y}
\newcommand{\Z}{\mathcal Z}

\newcommand{\Cb}{\mathbb C}

\newcommand{\Eb}{\mathbb E}
\newcommand{\Fb}{\mathbb F}
\newcommand{\Gb}{\mathbb G}
\newcommand{\Pb}{\mathbb P}
\newcommand{\Qb}{\mathbb Q}

\newcommand{\Dfrak}{\mathfrak D}
\newcommand{\Efrak}{\mathfrak E}

\newcommand{\Ifrak}{\mathfrak I}

\newcommand{\hsp}{\mspace{1.5mu}}
\newcommand{\hspb}{\mspace{-1.5mu}}
\newcommand{\spc}{\,\,\,}

\newcounter{saveenumi}
\newcommand{\save}{\setcounter{saveenumi}{\value{enumi}}}
\newcommand{\restore}{\setcounter{enumi}{\value{saveenumi}}}

\makeatletter

\DeclareFontFamily{T1}{greek}{}
\DeclareFontShape{T1}{greek}{m}{n}{ <-> psyr }{}
\DeclareFontShape{T1}{greek}{o}{n}{ <-> [0.9] psyr }{}
\DeclareFontShape{T1}{greek}{m}{sl}{ <-> rpsyro }{}

\DeclareFontFamily{T1}{title}{}
\DeclareFontShape{T1}{title}{m}{n}{ <-> ptmr }{}
\DeclareFontShape{T1}{title}{m}{sc}{ <-> ptmrc }{}
\DeclareFontShape{T1}{title}{m}{it}{ <-> ptmri }{}
\DeclareFontShape{T1}{title}{b}{n}{ <-> ptmb }{}

\DeclareFontFamily{T1}{computer}{}
\DeclareFontShape{T1}{computer}{m}{n}{ <-> pcrr }{}

\newcommand{\aupalpha}{\textsf{a}}

\newcommand{\aalpha}{\textsf{\textsl{a}}}

\newcommand{\pstarsplus}{(\raisebox{-1pt}{\text{\scriptsize\EightStarTaper}}_{\hspace{-1pt}s})}
\newcommand{\pstarsm}{(\raisebox{-1pt}{\text{\scriptsize\EightStarTaper}}^-_{\hspace{-1pt}s})}

\newcommand{\deecmp}{$\mathbb D$-complete\xspace}
\newcommand{\alphaproper}{$\aupalpha$-proper\xspace}

\newcommand{\lrg}{\operatornamewithlimits{Lrg}}

\newcommand{\shade}{\partial}
\newcommand{\derivative}{\partial}

\newcommand{\poset}{\Q}

\newcommand{\spo}{\R}
\renewcommand{\ggen}{\Game_{\mathrm{gen}}}

\renewcommand{\gen}{\operatorname{gen}}
\newcommand{\genc}{\gen^+}

\newcommand{\genmod}{\M}
\newcommand{\idealmod}{\Dfrak}

\newcommand{\trsup}{\operatorname{trsup}}

\newcommand{\ns}{\mathrm{NS}}

\renewcommand{\diamond}{\diamondsuit}

\newcommand{\nnr}{\mathrm{NNR}}
\newcommand{\pproper}{($\p$)\xspace}
\newcommand{\liminv}[2]{\varprojlim{}_{{#1}<{#2}}}
\renewcommand{\injlim}{\varinjlim}

\newcommand{\ppcif}[1]{\power(\powcif{#1})}

\newcommand{\tar}{\,\textrightarrow\,}
\newcommand{\subseteqtl}[1]{\subseteq^{#1}_{\mathrm{tl}}}
\newcommand{\eqtl}[1]{=^{#1}_{\mathrm{tl}}}

\newlength{\cmd}
\newcommand{\idealmodp}[2][]{\settowidth{\cmd}{$#1$}
\ifthenelse{\lengthtest{\cmd = 0pt}}{\Dfrak_{#2}}{\Dfrak_{#2}^{#1}}}
\newcommand{\embeds}[1][]{\preccurlyeq_{#1}}
\newcommand{\imbeds}[1][]{\preccurlyeq_{#1}^i}

\newcommand{\bigexta}{\hspace{1pt}\bigext}

\newcommand{\cofle}{\lesssim}
\newcommand{\cofeq}{\cong}
\newcommand{\code}[1]{{\Cb}^{#1}}
\newcommand{\codep}[1]{{\Cb}_{P}^{#1}}
\newcommand{\gcode}[1]{{\Gb}^{#1}}
\newcommand{\gcodep}[1]{{\Gb}^{#1}_{P}}
\newcommand{\pcode}[1]{{\Pb}^{#1}}
\newcommand{\pcodep}[1]{{\Pb}^{#1}_{P}}
\newcommand{\qcode}[1]{{\Qb}^{#1}}
\newcommand{\qcodep}[1]{{\Qb}^{#1}_{P}}
\newcommand{\codeg}[1]{C(#1)}
\newcommand{\codegfin}[2]{C(#1,#2)}
\newcommand{\genfn}{\Eb}
\newcommand{\genfnf}{\Fb}
\newcommand{\qsep}{\sim_{\tu{sep}}}
\newcommand{\qasym}{\sim_{\tu{asym}}}
\newcommand{\trind}{\operatorname{trind}}

\setlist{leftmargin=*,noitemsep}

\reversemarginpar

\showhyphens{chooser}

\begin{document}

\title{A strong antidiamond principle compatible with   $\ch$ \tnoteref{t1}}
\author{James Hirschorn\fnref{f1}}
\tnotetext[t1]{First draft, July 4, 2007. Final version, March 23, 2008.}
\begin{keyword}
Diamond principle \sep continuum hypothesis \sep properness
  parameter. 

\MSC Primary 03E35; Secondary 03E05, 03E15, 03E50.
\end{keyword}
\ead{j\_hirschorn@yahoo.com}
\ead[url]{http://logic.univie.ac.at/~hirschor}
\address{Thornhill, ON, CANADA 
\\ \vspace*{20pt}{\usefont{T1}{title}{m}{it}\normalsize
In memory of Sylvester Hirschorn}
\vspace*{-29pt}}
\fntext[f1]{The author began a study of the parameterized properness
   theory  while being
   supported by Consorcio Centro de Investigaci\'on Matem\'atica, Spanish
  Government grant No.~SB2002-0099. The initial stages of this research
  were supported by Fonds zur F\"orderung der wissenschaftlichen
  Forschung, Lise Meitner Fellowship,  Project No.~M749-N05.}

\begin{abstract}
A strong antidiamond principle $\pstarc$ is shown to be
consistent with~$\ch$. This principle can be stated as a ``$P$-ideal dichotomy'':
\emph{every $P$-ideal on $\oone$ 
\tu(i.e.~an ideal that is $\sigma$-directed
under inclusion modulo finite\tu) 
either has a closed unbounded subset of $\oone$ locally inside of it, 
or else has a stationary subset of $\oone$ orthogonal to it}. 
We rely on Shelah's theory of \emph{parameterized properness} 
for $\nnr$ iterations,
and make a contribution to the theory
with a method of constructing the properness parameter
simultaneously with the iteration. 
Our handling of the application of the $\nnr$ iteration theory involves definability of
forcing notions in third order arithmetic, analogous to Souslin
forcing in second order arithmetic. 
\end{abstract}

\maketitle
\section{Introduction}
\label{sec:introduction}

It is a remarkable fact (i.e.~theorem of $\zfc$ and the existence of
some large cardinals) that if $\phi$ and $\psi$ are two $\Pi_2(\ns)$
sentences in the language of set theory, both of which can
individually be forced to hold in the structure
$(H_{\aleph_2},\in,\ns)$ ($\ns$ denotes the ideal on nonstationary
subsets of $\oone$), then their conjunction can also be forced to hold
in this structure. Indeed Woodin has constructed a canonical model
$\pmax$ where the $\Pi_2$ theory over $(H_{\aleph_2},\in,\ns)$ is
maximal~(cf.~\cite{MR1713438}). In this model Cantor's Continuum Hypothesis
($\ch$) is false. The question of whether the $\Pi_2$ theory can be
maximized over structures $(H_{\aleph_2},\in,\ns)$ satisfying $\ch$,
is a major obstacle to further progress on the Continuum Hypothesis. This is closely
related to the question of whether there are forcing axioms analogous
to the Proper Forcing Axiom ($\pfa$) or Martin's Maximum ($\mm$) that
are consistent with $\ch$. 

Specifically, there is the test question below of Shelah and
Woodin asking whether the above mentioned remarkable fact still holds
if we take the conjunctions of $\phi$ and $\psi$ with $\ch$: Let
$\Pi_2(\ns)$ denote the collection of all $\Pi_2$ sentences in the language of
set theory (i.e.~of the form $\forall x\,\exists y\spc \varphi(x,y)$
where $\varphi$ has no unbounded quantifiers) with the added unary
predicate $\ns$. 

\begin{question}
\label{q-1}
Are there $\Pi_2(\ns)$ sentences $\phi$ and $\psi$ such that both
\begin{enumerate}[label=\tu{(\arabic*)}, ref=\arabic*, widest=2]
\item $(H_{\aleph_2},\in,\ns)\models\ulc \phi \lands \ch\urc$ and
\item $(H_{\aleph_2},\in,\ns)\models\ulc \psi\lands \ch\urc$
\end{enumerate}
can be forced, yet provably $(H_{\aleph_2},\in,\ns)\models\ulc\phi\lands\psi\to\lnot\ch\urc$\tu?
\end{question}

\noindent Woodin has conjectured a positive answer, which would
indicate that the $\Pi_2(\ns)$ theory over $H_{\aleph_2}$ cannot be
maximized for models of $\ch$, and thus there are ``disjoint'' $\Pi_2$-rich models of $\ch$. 

There has been much 
work done on maximizing the $\Pi_2$ theory in the presence
of $\ch$, where the idea to show that some `strong' $\Pi_2$ statement
is consistent with $\ch$. A major breakthrough in this line of
research was the Abraham--Todor\v cevi\'c $P$-ideal dichotomy $\pstar$ that
implies many well known $\Pi_2$ consequences of $\pfa$, and yet was
shown to be consistent with $\ch$ (\cite{MR1441232}). 
In the present paper, we push the boundary of maximizing the $\Pi_2$
theory over $\ch$, by proving that the strengthening $\pstarc_\oone$ (see
below) of $\pstar_\oone$ is consistent with $\ch$. 
(Technically speaking, $\pstarc$ is a variant of $\pstar$, 
but once we can obtain a model of $\pstarc$ with CH 
we can easily obtain $\pstar$ simultaneously, whereas the converse is
false.)

We believe there are good indications that $\pstarc_\oone$ is strong
enough to serve as one of the two $\Pi_2(\ns)$ statements in giving a
positive answer to question~\ref{q-1}. In particular, there is
example~\ref{x-6} discussed below, which tells us that in an iterated
forcing construction, given a $P$-ideal $\ideal$ 
where the second alternative of $\pstarc$ fails, depending on the
initial stages of the iteration, 
we may or may not be able to force the first alternative while at the
same time making sure we do not add reals. And this
suggests serious difficulties in obtaining a forcing axiom consistent
with $\ch$ that would imply $\pstarc$. 

Consider now the following dichotomy of Eisworth, based on the
$P$-ideal dichotomy~$\pstar$.

\begin{principle}{$\pstarc$}
For every ordinal $\theta$ of uncountable cofinality,
every $\sigma$-directed downward closed (i.e.~under subsets) subfamily $\I$ of
$(\powcif\theta,\subseteqfnt)$ has either
\begin{enumerate}[labelsep=\label@sep, label=(\arabic*), ref=\arabic*]
\item\label{item:12} a closed uncountable subset of $\theta$ locally in $\I$,
\item\label{item:164} a stationary subset of $\theta$ orthogonal to $\I$,
\end{enumerate}
\end{principle}
where $C\subseteq\theta$ is \emph{locally in $\I$} means
$\powcif{C}\subseteq\I$, and $S\subseteq\theta$ is \emph{orthogonal to
  $\I$} means $S\cap x$ is finite for all $x\in\I$. For some fixed
ordinal $\theta$, $\pstarc_\theta$ denotes the restriction of
$\pstarc$ to $\theta$. The original principle $\pstar$ is also a
dichotomy, where in the first alternative~\eqref{item:12}, ``closed
uncountable'' is weakened to ``uncountable''; and the second
alternative~\eqref{item:164} is strengthened to the existence
of a countable decomposition of~$\theta$ into pieces orthogonal to $\I$.
Other similar variations are possible such as the principle
$\pstarsplus$ studied in~\cite{Hir-comb} (actually this is a weakening
of $\pstarc$ optimal with respect to permitting the existence of a nonspecial Aronszajn tree). 

The main result of this research is that $\pstarc_\oone$ is consistent with $\ch$.

\begin{thm}
\label{u-3}
$\pstarc_\oone$ is consistent with $\ch$ relative to $\zfc$. 
\end{thm}

\noindent This answers Shelah's question~\cite[Question~2.17]{MR1804704}.
The methods here can also be modified in the straightforward
manner to obtain the
consistency of the unrestricted principle $\pstarc$ with $\ch$
relative to a supercompact cardinal.

It was already known that $\pstarc$ is consistent with the failure of~$\ch$. 
The following theorem is due to Eisworth, at least in the case $\theta=\oone$, 
and is proved in~\cite{H2}. 

\begin{thm}
\label{u-2}
$\pfa$ implies $\pstarc$. 
\end{thm}

The principle $\pstar$ is already very powerful with applications to
uncountable objects appearing in other areas of mathematics such as
measure theory. 
The principle $\pstarc$ moreover brings into play the most significant structural
property of $H_{\aleph_2}$, as compared to $H_{\aleph_1}$; thus,
unlike $\pstar$, it is not a $\Pi_2$ statement (i.e.~without the
predicate $\ns$). 
Let us briefly consider a couple of examples of how such
principles are applied to combinatorial objects. 

\begin{example}
\label{x-12}
As demonstrated in~\cite{MR1441232}, 
to any tree $\T=(T,\len T)$ we can associate the ideal $\T^\perp$ of all
countable subsets $x$ of $T$ perpendicular to the tree, i.e.~every node
has at most finitely many predecessors in $x$. Then, for example, if
$\T$ has all levels countable then $\T^\perp$ is a $P$-ideal. And
$\T^\perp$ has an uncountable set orthogonal to it iff $\T$ has an
uncountable branch. 
\end{example}

\begin{example}
\label{x-13}
If $\vec x=(x_\delta:\delta<\theta$ with $\cof(\delta)=\omega)$ is a
sequence where each $x_\delta\subseteq\delta$ is a cofinal subset of order
type $\omega$, then we can associate an ideal $\vec x^\perp$ of all
countable $y\subseteq\theta$ orthogonal to $\vec x$, i.e.~$x_\delta\cap y$ is finite for
all $\delta$. Then $\vec x^\perp$ is a $P$-ideal, with no orthogonal
subset of $\theta$ of order type $\omega^2$. And $\vec x$ is a
club-guessing sequence, in the weak sense, iff it has no closed unbounded subset of
$\theta$ locally in $\vec x^\perp$. See
e.g.~\cite{H2},~\cite[Ch.~XVIII,~Problem~1.9]{MR1623206}.
\end{example}

Let us mention some of the challenges that need to be overcome to
prove theorem~\ref{u-3}. First of all, it is known that
$\pstarc_\oone$ negates the relatively weak consequence of $\diamond$ 
that there is a club guessing sequence on $\oone$ (see example~\ref{x-13},~\cite{H2}). 
Therefore, we cannot use $\aupalpha$-proper forcing to obtain
theorem~\ref{u-3}. Moreover, $\pstarc_\oone$ implies that all
Aronszajn trees are special (\cite{H2}), and thus there are
significant difficulties in using 
Shelah's theory in~\cite[Ch.~XVIII, \Section2]{MR1623206} that he
developed for negating
club guessing with $\ch$. We use his newest $\nnr$ (no new reals)
iteration theory from~\cite{math.LO/0003115}, called parameterized
properness, which was developed in
order to obtain the negation of club guessing sequences together with
all Aronszajn trees being special simultaneously with $\ch$. This
involved devising new techniques for constructing the properness
parameters. We also discuss the possibility of using the $\nnr$ iteration theory
in~\cite[Ch.~XVIII, \Section2]{MR1623206} (cf.~\Section\ref{sec:trind-properness}). 

We say that two families
$\H,\I\subseteq\powcif\theta$ are \emph{orthogonal}, 
written $\H\perp\I$ if $x\cap y$ is
finite for all $x\in\H$ and $y\in\I$. The following example can be
obtained by a straightforward construction of an $(\oone,\oone)$ gap
in $(\powcif\oone,\subseteqfnt)$.

\begin{example}[$\diamond$]
\label{x-6}
There exist two $\sigma$-directed subfamilies $\H$ and $\I$ of
$(\powcif\oone,\subseteqfnt\nobreak)$ such that $\H\perp\I$ and neither has a
stationary set orthogonal to it. 
\end{example}

It follows from this (see~\cite{math.LO/0003115}) that we cannot
obtain a model of $\pstarc_\oone+\ch$ by a straightforward
iteration, where at each stage a club is forced locally inside a 
$P$-ideal with no stationary subset of $\oone$ orthogonal to it. The
above example shows that we will run into the so called ``disjoint clubs''. 

Many forcing notions, such as Cohen forcing and random forcing, 
can be represented as sets of reals that have simple definitions.
This fact has been well used to obtain results in the
descriptive set theory of the reals. For example,
in~\cite[\Section5]{rst:S} the simplicity of the 
representations of random forcing and amoeba forcing as 
sets of reals, respectively, is used (rather spectacularly) to construct a
nonmeasurable $\Sigma^1_3$-set of reals from some real computing
$\oone$ over $L$. The property of being simply definable as a set of
reals is particular useful in the iteration of such forcing
notions. Judah--Shelah gave a systematic treatment of these forcing
notions in~\cite{MR973109}, where they are named \emph{Souslin
  forcing}, with the emphasis on the iteration of Souslin forcing
notions. 

In overcoming the ``disjoint clubs'' obstacle by constructing a
properness parameter suitable for our iteration, we entered an
analogous situation but in the realm of third order arithmetic, instead
of the second order realm of set theory of the reals. We used the
fact that our forcing notions can be represented as simply definable
subsets of $\power(\oone)$ to establish nice properties of their
iterations. For example, we were able to show that our forcing
notions, which have cardinality $2^{\aleph_1}$, satisfy properties
such as commutativity, of both their iterations and their generic
objects (analogous to the fact that if $r$ is a random
real over $V$ and $s$ is a random real over $V[r]$ then $r$ is a
random real over $V[s]$).

\subsection{Terminology}
\label{sec:terminology}

We use standard order theoretic notation and terminology. Thus for a
family $\F$ of subsets of some fixed set $S$, we let $\downcl\F$
denote the downwards closure in the inclusion order,
i.e.~$\downcl\F=\bigcup_{x\in\F}\power(x)$. The definition of the
upwards closure $\upcl\F$ is symmetric. When want to take the
downwards closure with respect to some other quasi ordering $\altle$ 
of $\power(S)$, we write $\downcl{(\F,\altle)}$. E.g.~we will consider
the \emph{almost inclusion} quasi ordering $\subseteqfnt$, where
$x\subseteqfnt y$ means $x\setminus y$ is finite. 
A \emph{$P$-ideal} is an ideal that is also
$\sigma$-directed in the $\subseteqfnt$-ordering; furthermore, a
$P$-ideal on some specified set $S$ is always assumed to contain every
finite subset of $S$. A subset $A\subseteq
Q$ of a quasi order $(Q,\altle)$ 
is \emph{cofinal} if every $q\in Q$
has an $a\in A$ with $q\le a$. While a subset $A\subseteq P$ of a
strict partial order $(P,<)$ is \emph{cofinal} if every $p\in P$ has
an $a\in A$ with $p<a$, e.g.~we will consider cofinal subsets of some
structure $(M,\in)$. 

A subfamily of $\H\subseteq\powcif S$ is called \emph{cofinal} if it
is cofinal in the inclusion ordering, i.e.~for
all $a\in\powcif S$ there exists $b\in\H$ with $a\subseteq b$. It is
\emph{closed} if whenever $a_0\subseteq a_1\subseteq\cdots $ is a sequence of elements
of $\H$ then so is $\bigcup_{n<\omega}a_n\in\H$, and \emph{stationary}
if it intersects every closed set. 

We write $q\ge p$ for $q$ extends $p$, i.e.~carries more information
than $p$. This is clearly more natural than, the perhaps more common,
$q\le p$, especially in the context of $\aupalpha$-properness and
more generally parameterized properness (cf.~definition~\ref{d-3}).
As usual, $\Gen(M,P)$ denotes the family of ideals $G\subseteq P$
that are generic over $M$, while $\gen(M,P)$ is the set of all
$(M,P)$-generic elements of $P$. And $\Genc(M,P)$ is the subfamily of
all $G\in\Gen(M,P)$ that have a common extension in $P$, and $\genc(M,P)$
is set of all \emph{completely $(M,P)$-generic} elements $q$ of $P$,
meaning $q$ extends some member of $D$ for every dense $D\subseteq P$
in $M$. Every $q\in\genc(M,P)$ uniquely determines a member of
$\Genc(M,P)$, namely $\{p\in P\cap M:p\le q\}$, which
we denote as $\dot G_P[M,q]$. \emph{Complete properness} 
has the same formulation as
properness using countable elementary submodels, but with
$(M,P)$-generic replaced by completely $(M,P)$-generic. Thus a forcing
notion is completely proper iff it is proper and adds no new reals.
For more about proper forcing see e.g.~\cite{hbst}. 

Unless otherwise stated, for some function $f$ and some subset
$X\subseteq\dom(f)$, we write $f[X]$ for the image $\{f(x):x\in X\}$
of $X$ under $f$. Hopefully this will not cause confusion, because we
also use square brackets for generic interpretations. 

\section{Parameters for Properness}
\label{sec:params}

Although the following definition looks slightly different, it is almost
the same as the definition of a ``reasonable parameter'' 
from~\cite[\Section1]{math.LO/0003115}. The main difference is that we
only require cofinality in $H_\lambda$ (cf.~\eqref{item:61}); indeed,
it is noted in~Remark~1.10(2) of that article that this is sufficient.

\begin{defn}
\label{d-2}
For a regular cardinal $\lambda$, 
a \emph{$\lambda$-parameter for properness} is 
a pair $(\vec\A,\idealmod)$ for which there exists a sequence of
regular cardinals $\vec\mu=(\mu_\alpha:\alpha<\oone)$
with $\mu_0\ge\lambda$ and 
$H_{\mu_\alpha}\in H_{\mu_\beta}$ for all $\alpha<\beta$, such
that
\begin{enumerate}[label=(\roman*), ref=\roman*, widest=ii]
\item\label{item:59} $\vec\A$ is an $\oone$ sequence where
  $\A_\alpha\subseteq\powcif{H_{\mu_\alpha}}$ is stationary for
  all $\alpha<\oone$, and for every $M\in\A_\alpha$, 
  \begin{enumerate}[widest=b]
  \item\label{item:52} $M\prec H_{\mu_\alpha}$,
  \item\label{item:70}    $(\vec\mu\restriction\alpha,\vec\A\restriction\alpha)\in M$,
  \end{enumerate}
\save
\end{enumerate}
$\vec\A$ is called the \emph{skeleton} of the parameter, and the
$\rank$ function on $\injlim\A=\bigcup_{\alpha<\oone}\A_\alpha$ is defined by
\begin{equation}
  \label{eq:29}
  M\in\A_{\rank(M)}.
\end{equation}
We use the notation $\A_{<\alpha}=\bigcup_{\xi<\alpha}\A_\xi$. 
\begin{enumerate}[label=(\roman*), ref=\roman*, widest=ii]
\restore
  \item\label{item:61} For all $M\in\A$,
    if $\rank(M)=0$ then
    $\idealmod(M)=\{0\}$; and if $\rank(M)>0$ then $\idealmod(M)$
    is a \emph{nonempty} collection of subsets of $\A_{<\rank(M)}\cap M$
    so that for each element  $X\in\idealmod(M)$, every $\xi<\rank(M)$ and every
    $a\in M\cap H_\lambda$ has a $K\in X$ such that
   \begin{enumerate}[widest=b]
   \item\label{item:88} $\rank(K)\ge\xi$, 
    \item\label{item:79} $X\cap K\in\idealmod(K)$,
    \item\label{item:80} $a\in K$. 
    \end{enumerate}
\save
\end{enumerate}
\end{defn}

\fxnote{Comment}

\begin{prop}
\label{p-15}
$\rank(M)<\oone\cap M$. 
\end{prop}
\begin{pf}
By~\eqref{item:70}.
\end{pf}

\begin{prop}
\label{p-31}
$\rank(K)\le\rank(M)$ for all $K\in\injlim\A\cap M$.
\end{prop}
\begin{pf}
By~\eqref{item:59}.
\end{pf}

\begin{prop}
\label{p-3}
For all $M\in\injlim\A$ with $\rank(M)>0$, 
$\idealmod(M)$ is closed under supersets in $\power(\A_{<\rank(M)}\cap
M)$. 
\end{prop}
\begin{pf}
A simple induction on $\rank(M)$. 
\end{pf}

The properness parameter is often utilized through the following game.
It is a simplification of the game in \cite[Definition~1.5]{math.LO/0003115},
that appears to serve the same purpose. In any case, the two games are
equivalent for the properness parameters we will be using (cf.~definition~\ref{ex:1}). 

\begin{defn}
\label{d-7}
Let $(\vec\A,\idealmod)$ be a properness parameter.
For each $M\in\injlim\A$ of positive rank, the \emph{Chooser Game}
$\game(M)=\game(\vec\A,\idealmod)(M)$ is defined as follows.
It is a two player game of length~$\omega$, 
where the \emph{challenger} moves first against the \emph{chooser}. 
On the $k\Th$ move (setting $X_{-1}=\A_{<\rank(M)}\cap M$):
\begin{itemize}
\item The challenger plays $X_k\subseteq X_{k-1}$ in $\idealmod(M)$.
\item The chooser plays $K_k\in X_k$, and $Y_k\subseteq X_k$ in $\idealmod(K_k)$.
\end{itemize}
The chooser wins the game 
if $\bigcup_{k<\omega} Y_k\cup\{K_k\}\in\idealmod(M)$. 
Otherwise the challenger wins. 

We say that the chooser has a
\emph{global winning \tu(nonlosing\tu) strategy in the game
  $\game(\vec\A,\idealmod)$} if the chooser has a winning
(nonlosing) strategy in the game\linebreak $\game(\vec\A,\idealmod)(M)$ for
all $M\in\injlim\A$ with $\rank(M)>0$.
\end{defn}

Note that the chooser always has a valid move:

\begin{lem}
\label{l-12}
Every $X\in\idealmod(M)$ has a $Y\subseteq X$ in $\idealmod(M)$ such
that $Y\cap K\in\idealmod(K)$ for all $K\in Y$.
\end{lem}
\begin{pf}
The proof is by induction on $\rank(M)$. 
For each $a\in M\cap H_\lambda$ and each $\xi<\nobreak\rank(M)$, 
there exists $K_{a\xi}\in X$ with $a\in K_{a\xi}$,
$X\cap K_{a\xi}\in\idealmod(K_{a\xi})$ and $\rank(K_{a\xi})\ge\xi$. 
Then by the induction hypothesis,
there exist $Y_{a\xi}\subseteq X\cap K_{a\xi}$ in $\idealmod(K_{a\xi})$ 
such that $Y_{a\xi}\cap J\in\idealmod(J)$ for all $J\in Y_{a\xi}$. 
Now $\bigcup_{a\in M\cap H_\lambda}\bigcup_{\xi<\rank(M)}Y_{a\xi}\allowbreak\cup\{K_{a\xi}\}$ 
is in $\idealmod(M)$ and satisfies the desired conclusion, by proposition~\ref{p-3}.
\end{pf}

\begin{example}
\label{x-4}
The present example corresponds to a \emph{standard parameter}
in\linebreak \cite{math.LO/0003115}. 
Given a skeleton $\vec\A$ for a properness parameter,
define $\idealmod(M)$ by recursion on $\rank(M)$ as follows. 
Let $\idealmod(M)=\{0\}$ if $\rank(M)=0$, and $\idealmod(M)$ consist
of all subsets of $\A_{<\rank(M)}\cap M$ satisfying~\eqref{item:61} otherwise.
Condition~\eqref{item:59} on the
skeleton ensures that $\idealmod(M)\ne\emptyset$ for all $M\in\injlim\A$,
and thus this does indeed define a parameter for properness. 
The chooser has a global winning strategy in the
corresponding Chooser Game, because given $M\in\injlim\A$ of positive
rank, if $(a_n:n<\omega)$ enumerates $M\cap H_\lambda$ and
$\lim_{n\to\omega}\xi_n+1=\rank(M)$ (cf.~remark~\ref{r-1}),
then playing $K_n$ with $a_0,\dots,a_n\in K_n$ and
$\rank(K_n)\ge\xi_n$ (and $Y_k$ arbitrary) 
defines a winning strategy for the chooser in the
game $\game(\vec\A,\idealmod)(M)$.
\end{example}

The following definition is the case $\alpha=\beta$ 
of~\cite[Definition~2.8]{math.LO/0003115}. 
 
\begin{defn}
\label{d-3}
Let $(\vec\A,\idealmod)$ be a $\lambda$-parameter for properness.
A poset $P$ 
is \emph{proper for the parameter $(\vec\A,\idealmod)$} or
\emph{$(\vec\A,\idealmod)$-proper} if $P\in
H_\lambda$\footnote{In~\cite{math.LO/0003115} it is required that in
  fact $\power(P)\in H_\lambda$. Although it is harmless to ask for this,
  we have left it out.} and
there exists $a\in H_\lambda$ such that  
for all $M\in\injlim\A$ with $a,P\in M$:\label{pproper}
\begin{principle}{($\p$)}
for all  $p\in P\cap M$ and all $X\in\idealmod(M)$, 
there is an $(M,P)$-generic extension $q$ of $p$ 
such that
\begin{equation}
  \label{eq:27}
  \genmod(q)\cap X\in\idealmod(M),
\end{equation}
\end{principle}
where $\genmod(q)=\bigl\{M\in \powcif{H_\lambda}:q\in\gen(M,P)\bigr\}$.
\end{defn}

Note that equation~\eqref{eq:27} is trivial when $\rank(M)=0$. 

\begin{example}
\label{x-2}
Suppose $P$ is an $\aupalpha$-proper forcing notion with $\power(P)\in
H_\lambda$. 
Then $P$ is $(\vec\A,\idealmod)$-proper for every
$\lambda$-parameter for properness. 
\end{example}

\subsubsection{Tails}
\label{sec:tails}

\begin{defn}
\label{d-12}
For any $X\in\idealmod(M)$, a \emph{tail} of $X$ is a subset of the
form $\{K\in X:a\in K\}$ where $a\in M\cap H_\lambda$.
\end{defn}

\begin{prop}
\label{p-2}
The intersection of two tails of $X$ is itself a tail of $X$. 
\end{prop}

\begin{prop}
\label{p-13}
For all $X\in\idealmod(M)$ and all $J\in\injlim\A\cap M$, there exists
a tail $Y$ of $X$ such that $K\notin J$ for all $K\in Y$. 
\end{prop}

Any `reasonable' parameter has each $\idealmod(M)$ closed under
taking tails, but this is not a requirement.

\fxnote{Do not seem to need the tail qo.}

\fxnote{Comment}

\subsection{Properness parameters for shooting clubs}
\label{sec:prop-param-shoot}

When forcing a club subset of $\theta$, 
if $p$ is generic over $M$ then $p$ forces that
$\sup(\theta\cap M)$ is in the club. This explains why
$\aupalpha$-properness is unsuitable for purposes such as destroying a
club guessing sequence, because it can be used to guess the generic
club in the ground model. The following class of properness parameters is
designed to handle this difficulty by putting restrictions on these
suprema. For any family $\M$, the \emph{trace of the suprema of $\M$
  on $\theta$} is
\begin{equation}
  \label{eq:20}
  \trsup_\theta(\M)=\{\sup(\theta\cap M):M\in\M\}.
\end{equation}

\fxnote{Did not end up using this approach}

\begin{defn}
\label{ex:1}
Suppose $\theta$ is an ordinal of uncountable cofinality and $\vec\A$ is
a skeleton of a $\lambda$-parameter for properness for some $\lambda$.
For each $M\in\injlim\A$,
let us be given a countable family $\Omega(M)\subseteq\power(\theta)$
(normally, $\Omega(M)\subseteq\power(\theta\cap M)$).
For each $M\in\injlim\A$,  $\idealmodp{\Omega}(\vec\A)(M)=\idealmodp{\Omega}(M)$ 
is defined by recursion on $\rank(M)$. 
If $\rank(M)=0$ then define $\idealmodp{\Omega}(M)=\{\emptyset\}$; otherwise,
it is defined as the family of subsets of
$\A_{<\rank(M)}\cap M$ containing a subset
of the form
\begin{equation}
  \label{eq:75}
  X=\bigcup_{n<\omega}X_n\cup\{K_n\}
\end{equation}
where
\begin{enumerate}[label=(\roman*), ref=\roman*, widest=iii]
  \item\label{item:39} $K_0\in K_1\in\cdots$ is cofinal in $(M\cap H_\lambda,\in)$
    with $\lim_{n\to\omega}\rank(K_n)+1=\rank(M)$,
  \item\label{item:53} $X_n\in\idealmodp{\Omega}(\vec\A)(K_n)$ for all $n$,
  \item\label{item:34} every $x\in\Omega(M)$ has a tail $Y$ of $X$ with
    $\trsup_\theta(Y)\subseteq x$.
\save
\end{enumerate}
\fxnote{Comment}
\end{defn}

Condition~\eqref{item:34} is a geometrical restriction on the trace of
the suprema.

\begin{remark}
\label{r-1}
The limit in condition~\eqref{item:39} has its usual topological
meaning. Thus for any $f:\omega\to\ord$,
$\lim_{n\to\omega}f(n)=\alpha$ iff every $\xi<\alpha$ has a $k<\omega$
such that $f(n)\in (\xi,\alpha]$ for all $n\ge k$. Also,
$\limsup_{n\to\omega}f(n)$ should be interpreted as
$\lim_{n\to\omega}\sup_{i\ge n}f(i)$. 
\end{remark}

\begin{lem}
\label{p-11}
When $\rank(M)>0$, $\idealmodp{\Omega}(M)$ is closed under taking tails.
\end{lem}
\begin{pf}
The proof is by induction on $\rank(M)$.
Suppose $X\in\idealmodp\Omega(M)$
satisfies\linebreak \eqref{item:39}--\eqref{item:34},
and let $Y=\{K\in X:a\in K\}$ be a tail of $X$ for
some $a\in M\cap H_\lambda$.
Condition~\eqref{item:39} holds because $K_n\in Y$ for all but
finitely many $n$; condition~\eqref{item:53} holds for $Y$ by the induction hypothesis;
and condition~\eqref{item:34} is by proposition~\ref{p-2}.
\end{pf}

\begin{prop}
\label{p-16}
If $X\in\idealmodp\Omega(M)$, and $Y\subseteq X$ and 
$Y^K\subseteq X$  \tu($K\in Y$\tu) satisfy
\begin{enumerate}[label=\tu{(\alph*)}, ref=\alph*, widest=b]
\item\label{item:28} $Y^K\in\idealmodp{\Omega}(K)$, 
\item\label{item:26} for all $a\in M\cap H_\lambda$ and all\/~$\xi<\rank(M)$, there
  exists $K\in Y$ with $a\in K$ and $\rank(K)\ge\xi$,
\end{enumerate}
then $Y\cup\bigcup_{K\in Y}Y^K\in\idealmodp\Omega(M)$. 
\end{prop}

\begin{prop}
\label{p-38}
$(\vec\A,\idealmodp\Omega)$ is a properness parameter 
whenever $\idealmodp\Omega(M)\ne\emptyset$ for all $M\in\injlim\A$.
\end{prop}

An arbitrary mapping $\Omega$ may fail to define a properness parameter, 
i.e.~the $\idealmodp{\Omega}(M)$ may be empty for some $M$. 
We provide a general construct to avoid this.

\begin{defn}
\label{d-9}
A map $\Lambda:\injlim\A\to\power(\powcif\theta)$ is said to
\emph{instantiate} $\Omega$ if every $M\in\injlim\A$ with
$\rank(M)>0$, every finite $A\subseteq\Omega(M)$, every
$y\in\Lambda(M)$, every $a\in M\cap H_\lambda$ and 
every $\xi<\rank(M)$ has a $K\in\A_{<\rank(M)}\cap M$ such that 
\begin{enumerate}[label=(\roman*), ref=\roman*, widest=iii]
\item\label{item:10} $a\in K$,
\item\label{item:48} $\rank(K)\ge\xi$,
\item\label{item:38} $\sup(\theta\cap K)\in y\cap \bigcap A$,
\item\label{item:57} $y\cap\bigcap A\in\Lambda(K)$.
\end{enumerate}
\end{defn}

\begin{remark}
\label{r-7}
In all of our applications of instantiations, we will have
$\Lambda(M)=\upcl{\Omega(M)}$ for all $M$, and thus we can omit the
`$y$' in~\eqref{item:38} and~\eqref{item:57} in the verification. 
See e.g.~example~\ref{d-5}.
\end{remark}

\begin{lem}
\label{l-5}
If there exists a map instantiating $\Omega$, 
then $(\vec\A,\idealmodp{\Omega})$ is a $\lambda$-proper\-ness parameter,
i.e.~$\idealmodp{\Omega}(M)\ne\emptyset$ for all $M\in\injlim\A$. 
\end{lem}
\begin{pf}
Assume $\Lambda$ instantiates $\Omega$.
We proceed by induction on $\rank(M)$, with the induction hypothesis
that for all $y\in\Lambda(M)$, 
there exists $X\in\idealmodp{\Omega}(M)$ with $\trsup(X)\subseteq y$. 
Suppose $M\in\A_\alpha$ where $\alpha>0$, and $y\in\Lambda(M)$.
Let $(x_n:n<\omega)$ enumerate
$\Omega(M)$, letting $x_n=\theta\cap M$ in case $\Omega(M)=\emptyset$.
Let $(a_n:n<\omega)$ enumerate $M\cap H_\lambda$,
and fix a sequence $\xi_n<\nobreak\alpha$
($n<\omega$) such that $\limsup_{n\to\omega}\xi_n+1=\alpha$. We
recursively choose $K_{n+1}\ni K_n$ in $\A_{<\alpha}\cap M$ with $a_{n}\in K_{n}$,
$\rank(K_{n})\ge\xi_n$, $\sup(\theta\cap K_{n})\in y\cap \bigcap_{i=0}^n x_i$
and $y\cap\bigcap_{i=0}^n x_i\in\Lambda(K_n)$.
This is possible by~\eqref{item:10}--\eqref{item:57}. And for each~$n$, 
there exists $X_n\in\idealmodp{\Omega}(K_n)$ with $\trsup(X_n)\subseteq
y\cap \bigcap_{i=0}^n x_i$ by the induction
hypothesis. 
Then putting
 $X=\bigcup_{n<\omega}X_n\cup\{K_n\}$, 
 conditions~\eqref{item:39}--\eqref{item:34} 
 of definition~\ref{ex:1} are clearly satisfied,
 i.e.~$X\in\idealmodp{\Omega}(M)$, and also $\trsup(X)\subseteq y$, completing the induction. 
\end{pf}

\fxnote{Not in the APAL version.}

\fxnote{Cannot recall what this was for.}

\begin{example}
\label{d-5}
Suppose $\E\subseteq\{M\in\powcif{H_\lambda}:M\prec H_\lambda\}$ is
stationary, and $\vec\A$ is
a skeleton of a $\lambda$-properness parameter with $M\cap H_\lambda\in\E$
for all $M\in\injlim\A$.  Suppose $\H\subseteq\powcif\theta$ has no
stationary orthogonal set.
If $y_M\in\downcl\H$ ($M\in\E$) satisfies $x\subseteqfnt y_M$ for all
$x\in\H\cap M$, then defining $\Omega:\injlim\A\to\powcnt{\powcif\theta}$ by 
\begin{equation}
  \label{eq:35}
  \Omega(M)=\{y_M\setminus s:s\subseteq y_M\text{ is finite}\},
\end{equation}
we have that $(\vec\A,\idealmodp\Omega)$ is a
$\lambda$-properness parameter. This is instantiated by $\Lambda$, 
where $\Lambda(M)=\upcl{\Omega(M)}$ for all $M$. 
\end{example}
\begin{pf}
We apply lemma~\ref{l-5}. 
Thus given $M\in\A_\alpha$ with $\alpha>0$, a nonempty finite
$A\subseteq\Omega(M)$,  $a\in
M\cap H_\lambda$ and $\xi<\alpha$, 
we need to show there exists $K\in\A_{<\alpha}\cap M$
satisfying~\eqref{item:10}--\eqref{item:57}. 

Put $\B=\{K\in\A_\xi: a\in K\}$.
Since $\B\in M$ and $\B$ stationary,
 $\trsup(\B)\in\nobreak M$ is a stationary subset of $\theta$.
Thus there exists $\delta\in\trsup(\B)\cap\bigcap A\cap M$ by
elementarity, since $\H$ has no stationary orthogonal set.  
And by elementarity, we can find $K\in\nobreak\B\cap\nobreak M$ 
with $\sup(\theta\cap K)=\delta$, and hence $K$
satisfies~\eqref{item:10}--\eqref{item:38}. Condition~\eqref{item:57}
is satisfied because $y_K\subseteqfnt y_M$.
\end{pf}

\begin{prop}
\label{l-19}
Assume that $\Omega$ 
does define a properness parameter in definition~\tu{\ref{ex:1}}.
Then the chooser has a global winning strategy in the game
$\game(\vec\A,\idealmodp{\Omega})$. 
\end{prop}
\begin{pf}
It is immediate from the definition of $\idealmodp{\Omega}(M)$ and the
payout of the game, that the chooser wins so long as
$\lim_{n\to\omega}\rank(K_n)+1=\rank(M)$ 
and $(K_n:n<\omega)$ is cofinal in $M\cap H_\lambda$,
where $(K_n,Y_n)$ denotes
the chooser's $n\Th$ move. The chooser can guarantee this in the
obvious manner.
\end{pf}

\subsubsection{Direction constraints}
\label{sec:imag-subf-power}

In addition to the limitations imposed on the trace of the suprema, 
we shall want additional control over the `direction' in which the 
members of each $\idealmod(M)$ can grow. This simply means that we
want the set in equation~\eqref{eq:75} to be contained in some
subfamily of $\injlim\A$, but what is more, this family lives in some
forcing extension (and is thus `imaginary'). 

\begin{defn}
\label{d-28}
Let $\Z$ be a collection of pairs of the form $(P,\dot\B)$,
where $P$ is a forcing notion and $\dot\B$ is a $P$-name for a subset
of $\injlim\A$. 
Then we define subfamilies 
$\idealmodp\Omega(\Z)(M)=\idealmodp\Omega(\vec\A;\Z)(M)
\subseteq\idealmodp\Omega(\vec\A)(M)$ by recursion on $\rank(M)$ as follows.
If $\rank(M)=0$ then $\idealmodp\Omega(\vec\A;\Z)(M)=\{\emptyset\}$;
otherwise, it is the family of all members 
of $\idealmodp\Omega(\vec A)(M)$ such that the set $X$ in
equation~\eqref{eq:75} additionally satisfies
\begin{enumerate}[label=(\roman*), ref=\roman*, widest=iii]
\restore
\item\label{item:131} every $p\in P\cap M$ has a $q\in\gen(M,P,p)$
 and a tail $Y\subseteq X$ such that $q\forces Y\subseteq\dot\B$,
\end{enumerate}
for all $(P,\dot\B)\in\Z\cap M$.
\end{defn}

\begin{lem}
\label{l-2}
When $\rank(M)>0$, $\idealmodp\Omega(\vec\A;\Z)(M)$ is closed under taking tails.
\end{lem}
\begin{pf}
By lemma~\ref{p-11}. 
\end{pf}

Typically, we have $\dot\B$ a subset of the models over which $\dot
G_P$ is generic, 
in which case we automatically have that $P$ is $\Omega(\vec\A;\Z)$-proper. 

\begin{cor}
\label{p-29}
Suppose that 
$\idealmodp\Omega(\vec\A;\Z)$ is a properness parameter, 
i.e.~$\idealmodp\Omega(\vec\A;\allowbreak\Z)(M)\allowbreak\ne\emptyset$
for all $M\in\injlim\A$.
If $(P,\dot\B)\in\Z$ and $P\forces\dot\B\subseteq\{M\in\injlim\A:\allowbreak
\dot G_P\in\Gen(M,P)\}$, 
then $P$ is $\idealmodp\Omega(\vec\A;\Z)$-proper. 
\end{cor}

We generalize lemma~\ref{l-5} to the present setting. 

\begin{lem}
\label{l-61}
Assume that for every $M\in\injlim\A$ with $\rank(M)>0$, every
$(P,\dot\B)\in\Z\cap M$ has a\/ $q^M_{P,\dot\B}(p)\in\gen(M,P,p)$ for each $p\in P\cap M$,
such that every finite $A\subseteq \Omega(M)$, every finite
sequence
$(P_0,\dot\B_0),\dots,(P_{m-1},\dot\B_{m-1})$ in $\Z\cap M$, 
all finite $O_i\subseteq P_i\cap M$  \tu($i=0,\dots,m-1$\tu),
every $a\in M\cap H_\lambda$ and every $\xi<\rank(M)$ has a
$K\in\A_{<\rank(M)}\cap M$ satisfying
\begin{enumerate}[label=\tu{(\roman*)}, ref=\roman*, widest=iii]
\item\label{item:68} $a,(P_0,\dot\B_0),\dots,(P_{m-1},\dot\B_{m-1})\in K$,
\item\label{item:76} $\rank(K)\ge\xi$,
\item\label{item:83} $\sup(\theta\cap K)\in\bigcap A$,
\item\label{item:75} $\bigcap A\in\upcl{\Omega(K)}$,
\item\label{item:132} $q^M_{P_i,\dot\B_i}(p)\forces K\in\dot\B_i$ 
for all $p\in O_i$ for all $i=0,\dots,m-1$,
\item\label{item:1} $q^M_{P_i,\dot\B_i}(p)\ge q^K_{P_i,\dot\B_i}(p)$
  for all $p\in O_i$ for all $i=0,\dots,m-1$. 
\end{enumerate}
Then $\idealmodp\Omega(\vec\A;\Z)$ is a properness parameter.
\end{lem}
\begin{pf}
We establish the lemma 
by induction on $\rank(M)$, with the induction hypothesis
that there exists $X\in\idealmodp\Omega(\Z)(M)$ such that
$q^M_{P,\dot\B}(p)$ ($p\in P_i\cap M$) from the hypothesis of the lemma
witnesses  definition~\ref{d-28}\eqref{item:131} for $X$, 
for all $(P,\dot\B)\in\Z\cap M$. 

Suppose then that $M\in\injlim\A$ with $\alpha=\rank(M)>0$.
Let $(x_n:n<\nobreak\omega)$ enumerate $\Omega(M)$, 
letting $x_n=\theta\cap M$ in case $\Omega(M)=\emptyset$,
let $(a_n:n<\nobreak\omega)$ enumerate $M\cap H_\lambda$,
let $(P_n,\dot\B_n:n<\omega)$ enumerate $\Z\cap M$,
and let $(p_{ni}:i<\nobreak\omega)$ enumerate $P_n\cap M$ for each $n$.
Fix a sequence $\xi_n<\nobreak\rank(M)$ ($n<\omega$) 
such that $\limsup_{n\to\omega}\xi_n+1=\alpha$. 
We recursively choose $K_{n+1}\ni K_n$ in $\A_{<\alpha}\cap M$ 
with $a_{n},(P_0,\dot\B_0),\dots,(P_n,\dot\B_n)\in K_{n}$, $\xi_n\le\rank(K_{n})$, 
$\sup(\theta\cap K_{n})\in\bigcap_{i=0}^n x_i$,   
$\bigcap_{i=0}^n x_i\in\upcl{\Omega(K_n)}$
and $q^M_{P_i,\dot\B_i}(p_{ij})\forces K_n\in\dot\B_i$ and
$q^M_{P_i,\dot\B_i}(p_{ij})\ge q^{K_n}_{P_i,\dot\B_i}(p_{ij})$ 
for all $i,j=0,\dots,n$. 
This is possible by~\eqref{item:68}--\eqref{item:1}. 
And for each~$n$, there exists $X_n$ in $\idealmodp{\Omega}(\Z)(K_n)$ 
as in the induction hypothesis. Furthermore, by going to a tail of $X_n$, we may
assume that  $\trsup(X_n)\subseteq\bigcap_{i=0}^n x_i$ and that
$q^{K_n}_{P_i,\dot\B_i}(p_{ij})\forces X_n\subseteq\dot\B_i$ 
for all $i,j=0,\dots,n$.
Then putting  $X=\bigcup_{n<\omega}X_n\cup\{K_n\}$, 
 conditions~\eqref{item:39}--\eqref{item:34} of definition~\ref{ex:1} 
are clearly satisfied. And for condition~\eqref{item:131}, given
$(P,\dot\B)\in\Z\cap M$ and $p\in P\cap M$, say
$(P,\dot\B)=(P_i,\dot\B_i)$ and $p=p_{ij}$, 
$q^M_{P_i,\dot\B_i}(p)\forces\bigcup_{n\ge\max(i,j)}X_n\cup\{K_n\}\subseteq\dot\B_i$.
This proves that $X\in\idealmodp{\Omega}(\Z)(M)$ is as needed, completing the induction. 
\end{pf}

\subsection{The iteration}
\label{sec:iteration}

\begin{notn}
\label{o-3}
Let $\E\subseteq\powcif{H_\lambda}$. For a poset $P$ with $P\in H_\lambda$, 
if $G\subseteq P$ is a generic ideal over $V$, in $V[G]$ we define
\begin{equation}
  \label{eq:1}
  \E[G]=\{M[G]:M\in\E\text{, }P\in M\and G\in\Gen(V,P)\}. 
\end{equation}
\end{notn}

\begin{prop}
\label{p-5}
Let $P$ be a forcing notion that adds no new\/ $\omega$-sequences of ground model
elements \tu(e.g.~$P$ completely proper\tu). If
$\E\subseteq\powcif{H_\kappa}$ is closed and cofinal then so is
$\E[G]$. 
\end{prop}

A $\lambda$-properness parameter $(\vec\A,\idealmod)$ can be
interpreted in a forcing extension $V[G]$ by some forcing notion $P\in
H_\lambda$ as $(\vec\B,\Efrak)$ where $\B_\alpha=\A_\alpha[G]$ for all
$\alpha<\oone$ and $\Efrak(M[G])=\idealmod(M)$ for all $M\in\injlim\A$
with $M[G]\in\injlim\B$. Thus when we say that $P\forces\ulc\dot Q$ is
$(\vec\A,\idealmod)$-proper$\urc$, we mean that $V[G]\models\ulc\dot
Q[G]$ is $(\vec\B,\Efrak)$-proper$\urc$. 

Let us  now address the issue of ``long properness''. The following
is essentially \cite[Definition~1.8(c)]{math.LO/0003115}, but without
the requirement of complete properness. 

\begin{defn}
\label{d-31}
Let $(\vec\A,\idealmod)$ be a $\lambda$-properness parameter. 
A countable support iteration $(P_\xi:\xi\le\delta)$ is called
\emph{long $\idealmod$-proper} if $P_\delta\in H_\lambda$ and there
exists $a\in H_\lambda$ such that for all $M\in\injlim\A$ with
$a,P_\delta\in M$, for all $\xi<\delta$ in $M$, if
\begin{enumerate}[label=(\roman*), ref=\roman*, widest=iii]
\item\label{item:95} $q\in \gen(M,P_\xi)$,
\item\label{item:96} $X=\genmod(q)\cap M\in\idealmod(M)$, 
\item\label{item:89} $\dot p$ is a $P_\xi$-name 
such that 
\begin{enumerate}[widest=b]
\item\label{item:97} $q\forces\dot p\in P_\delta\cap M$,
\item\label{item:99} $q\forces\dot p\restriction\xi\in\dot G_{P_\xi}$,
\end{enumerate}
\save
\end{enumerate}
then there exists $r\in\gen(M,P_\delta)$ such that
\begin{enumerate}[label=(\roman*), ref=\roman*, widest=iii]
\restore
\item\label{item:134} $r\restriction\xi=q$, 
\item\label{item:101} $\genmod(r)\cap X\in\idealmod(M)$,
\item\label{item:133} $r\forces\dot p\in\dot G_{P_\delta}$. 
\end{enumerate}
\end{defn}

The following lemma says that the iteration is long $\idealmod$-proper
when each iterand is $\idealmod$-proper. It is proved
in~\cite[page~17, ``Proof of clause (c)'']{math.LO/0003115}. 

\begin{lem}
\label{l-68}
Assume that the chooser has a global winning strategy in the game
$\game(\idealmod)$. Suppose $(P_\xi,\dot Q_\xi:\xi<\delta)$ is a
countable support iteration such that $P_\xi\forces\ulc\dot Q_\xi$\tu{ is
$\idealmod$-proper}$\urc$ for all $\xi<\delta$. Then $(P_\xi:\xi\le\delta)$
is long $\idealmod$-proper. 
\end{lem}

Although~\cite{math.LO/0003115} is the first place we read the phrase
``long properness'', it is a familiar concept used for example in the
proof of the preservation of properness under countable support
iterations. Indeed lemma~\ref{l-68} corresponds to the ``Properness
Extension Lemma'' of~\cite[Lemma~2.8]{hbst} and the
``$\aupalpha$-Extension Property'' of~\cite[Lemma~5.6]{hbst}.

The following is a simplified, and somewhat weakened, 
version of~\cite[Main Claim~1.9]{math.LO/0003115},
which is the basic $\nnr$ iteration theorem for parameterized properness. 

\begin{thm}[Shelah]
\label{u-1}
Let $(\vec\A,\idealmod)$ be a properness parameter, 
where the choo\-ser has a global winning
strategy in  the game $\game(\idealmod)$.
Suppose $(P_\xi,\dot Q_\xi:\xi<\delta)$ is a countable
support iteration such that 
\begin{enumerate}[label=\tu{(\alph*)}, widest=b]
\item\label{item:135} $P_\xi\forces\dot Q_\xi$ is $\idealmod$-proper for all $\xi<\delta$,
\item\label{item:136} $P_\xi\forces\dot Q_\xi$ is \deecmp for all $\xi<\delta$.
\end{enumerate}
Then the limit $P_\delta$ of the iteration does not add new reals.
\end{thm}

We do not refer to the following strengthening of theorem~\ref{u-1} as
a ``theorem'' because, unlike theorem~\ref{u-1}, it does not stand alone in
the sense that the needed hypothesis is preserved at limits. 
I.e.~there is no conclusion from the hypothesis that 
the limit $P_\delta$ is $\idealmod_\delta$-proper for some properness parameter $(\vec\A,\idealmod_\delta)$; 
although in our application of lemma~\ref{u-7} this will be the case. 
In fact lemma~\ref{u-7} below is in the same spirit as the fact that an
iteration of proper \deecmp forcings of length less than $\omega^2$
adds no new reals. 

\fxnote{Perhaps define ``countable support'' here?}

\begin{lem}
\label{u-7}
Let $(\vec\A,\idealmod_\xi:\xi<\delta)$ be a sequence of properness
parameters such that the chooser has a global winning strategy in the
game $\game(\idealmod_\xi)$ for all $\xi<\delta$. 
Suppose that $(P_\xi,\dot Q_\xi:\xi<\delta)$ is a countable support
iteration satisfying
\begin{enumerate}[label=\tu{(\alph*)}, ref=\alph*, widest=b]
\item\label{item:27} $P_\xi$ is long $\idealmod_\xi$-proper for all $\xi<\delta$, 
\item\label{item:20} $P_\xi\forces\dot Q_\xi$ is \deecmp for all $\xi<\delta$.
\end{enumerate}
Then the limit $P_\delta$ does not add new reals. 
\end{lem}
\begin{pf}[Sketch of proof]
The proof is based on Abraham's proof in\linebreak \cite[\Section5]{hbst} of Shelah's
fundamental $\nnr$ iteration theorem that countable support iterations
of forcing notions that are both \alphaproper and  \deecmp do not add new reals.

There is a function $\Eb$ that is implicitly assumed to exist in
e.g.~the proof in \cite[Ch.~V,~\Section7]{MR1623206}, and is
thankfully made explicit in~\cite{hbst}. It takes arguments of the
form $(M,\vec M,(P_\xi,\dot Q_\xi:\xi<\gamma),G,p)$ as input, where
$M\prec H_\lambda$ is countable containing $(P_\xi,\dot Q_\xi:\xi<\gamma)$,
$\vec M=(M_\eta:\eta<\alpha)$ is an $\in$-tower of countable
elementary submodels with $M_0=M$,
$G\in\Gen(M,P_{\gamma_0},p\restriction\nobreak\gamma_0)\cap M_1$ 
for some $\gamma_0<\gamma$ in $M$ and $p\in P_\gamma\cap M$. 
The value $\Eb(M,\vec M,(P_\xi,\dot Q_\xi:\xi<\nobreak\gamma),G,p)$ returned
is an element $H\in\Gen(M,P_\gamma,p)$ extending $G$,
i.e.~$p\restriction\gamma_0\in G$ for all $p\in H$. 
The whole point
of introducing $\Eb$ is that it is definable from some parameters, and
thus the generic output by $\Eb$ can be found inside a suitable
elementary submodel. 

It is then shown in~\cite[Lemma~5.21]{hbst} that if the tower $\vec M$
is high enough, if $(P_\xi,\dot Q_\xi:\xi<\gamma)$ is a 
countable support iteration of \alphaproper
and \deecmp forcing notions, and if $G\in\Genc(M,P_{\gamma_0})$ is also
generic over all members of the tower, 
then $\Eb(M,\vec M,(P_\xi,\dot
Q_\xi:\xi<\nobreak\gamma),G,p)$ is completely generic over $M$, proving that $P_\gamma$
does not add new reals. 

By making the corresponding changes to the definition of $\Eb$, the
exactly analogous proof works for iterations $(P_\xi,\dot Q_\xi:\xi<\gamma)$ of
\deecmp forcing notions that are long $\idealmod$-proper for some
properness parameter $(\vec\A,\idealmod)$; we obtain a completely
generic $\Eb(M,\vec M,(P_\xi,\dot Q_\xi:\xi<\gamma),G,p)$ whenever the range of
$\vec M$ is in $\idealmod(N)$ for some $N\in\injlim\A$ of big enough
rank. 

Now consider the iteration from the hypothesis of the lemma. 
By the hypotheses~\eqref{item:27} and~\eqref{item:20}, for every $\gamma<\delta$,
$\Eb(M,\vec M,(P_\xi,\dot Q_\xi:\xi<\gamma),G,p)$ is completely generic 
for all suitable $\vec M$ and $G$. To show that $P_\delta$ does not
add reals, we want to find a completely $(M,P_\delta)$-generic
ideal. Although we cannot take $\Eb(M,\vec M,(P_\xi,\dot Q_\xi:\xi<\delta),G,p)$
since we do not have a $\idealmod_\delta$, we can still go through
the proof of~\cite[Lemma~5.21]{hbst} to obtain complete genericity,
by using $\vec M=\vec M^0\bigexta\vec M^1\bigexta\cdots$, 
where $\xi_0<\xi_1<\cdots$ is cofinal in $\delta\cap M$, $N_0\in
N_1\in \cdots$ in $\injlim\A$ are of big enough rank, 
and the range of $\vec M^n$ is  in $\idealmod_{\xi_n}(N_n)$ for all $n<\omega$. 
\end{pf}

\section{The Forcing Notions}

\label{sec:forcing-notion}

In the context of  an ordinal $\theta$ of uncountable cofinality, $\kappa$
will always denote a regular cardinal $\kappa\ge(|\theta|^{\aleph_0})^+$; 
and in the context of a cardinal $\kappa$, we let $\lambda$ denote a
regular cardinal that is sufficiently large, by which we mean 
$\lambda\ge|H_\kappa|^+=(2^{<\kappa})^+$. 
Thus in the most important case $\theta=\oone$, taking $\kappa$ and
$\lambda$ to be the least sufficiently large regular cardinals and
assuming $\ch$ and $2^{\aleph_1}=\aleph_2$ (e.g.~assuming~$\gch$),
\begin{align}
  \label{eq:88}
  |\theta|=\aleph_1&&\kappa=(\aleph_1^{\aleph_0})^+=\aleph_2&&
  \lambda=(2^{\aleph_1})^+=\aleph_3.
\end{align}

\begin{defn}
\label{d-21}
Let $S$ be a set. 
For two families $\F\subseteq\H\subseteq\power(S)$, we let
\begin{equation}
  \label{eq:142}
  \shade_\H(\F)=\{x\subseteq\theta:\upcl x\cap\F\text{ is cofinal in }(\H,\subseteqfnt)\},
\end{equation}
i.e.~the set of all $x$ such that $\{y\in\F:x\subseteq y\}$ is
$\subseteqfnt$-cofinal in $\H$. We write $\shade(\H)$
for~$\shade_\H(\H)$, and we write $\alpha\in\shade_\H(\F)$ to indicate
that $\{\alpha\}\in\shade_\H(\F)$. 
\end{defn}

\begin{prop}
\label{p-69}
$\shade_\H(\F)\subseteq\downcl\F$ whenever $\H$ is nonempty. 
\end{prop}

\subsection{Forcing notion for shooting clubs}
\label{sec:forc-noti-shoot}

The following forcing notion is equivalent 
to the forcing notion $\R(\H,\C(\downcl\H))$ from~\cite{Hir-comb}. 
Thus the forcing notion in definition~\ref{def:consistency}
is a special case of a more general class of forcing notions studied
there. Many of the main results here, with the exception of the new result
in lemma~\ref{l-14}, follow from the general theory developed
in~\cite{Hir-comb}. We will provide direct proofs for most of the
results. 

\begin{defn}
\label{def:consistency}
For some ordinal $\theta$ of uncountable cofinality, 
let $\H\subseteq\powcif\theta$ be nonempty.
Then define $\poset(\H)$ to be the poset consisting of all pairs $p=(x_p,\X_p)$ where
\begin{samepage}
 \begin{enumerate}[label=(\roman*), ref=\roman*, widest=iii]
   \item\label{item:33} $x_p\in\shade(\H)$ is a closed subset of $\theta$, 
   \item\label{item:6} $\X_p$ is a countable family of subsets of $\H$ with 
      \begin{equation*}
       x_p\in\derivative_{\H}(\J)\espc\text{ for all $\J\in\X_p$},
     \end{equation*}
 \save
 \end{enumerate}
ordered by $q$ extends $p$ if
 \begin{enumerate}[label=(\roman*), ref=\roman*, widest=iii]
\restore
   \item $x_q\sqsupseteq x_p$ (i.e.~$x_q$ end-extends $x_p$ with respect to the
     ordinal ordering),
   \item $\X_q\supseteq\X_p$.
 \end{enumerate}
For an ideal $G\subseteq\poset(\H)$, we write $C_G=\bigcup_{p\in
  G}x_p$. 
And we write $0_{\poset(\H)}$ for the condition $(\emptyset,\emptyset)$. 
\end{samepage}
\end{defn}

Our poset forces the following desired result.

\begin{prop}
\label{p-10}
$\poset(\H)\forces\ulc C_{\dot G_{\poset(\H)}}\tu{ is locally in }\H\urc$.
\end{prop}
\begin{pf}
By proposition~\ref{p-69}.
\end{pf}

\begin{lem}
\label{l-16}
Suppose $\J\subseteq\H$ is cofinal in $(\H,\subseteqfnt)$.
Then $\poset(\H)\forces\ulc \exists y\in\H\,\,C_{\dot G_{\poset(\H)}}\setminus y$\tu{ is
locally in }$\J\urc$.
\end{lem}
\begin{pf}
Observe that the set of all $p\in\poset(\H)$ containing $\K\in\X_p$ of
the form $\K_y=\{x\cup y:x\in\J\}$ for some $y\in\downcl\H$ is
dense. That is, given $p\in\poset(\H)$,
$(x_p,\X_p\cup\nobreak\{\K_{x_p}\})\in\poset(\H)$ since $x_p\in\shade_\H(\K_{x_p})$.\end{pf}

\begin{prop}
\label{p-28}
The class $\poset$ is provably equivalent to a $\Delta_0$-formula.
\end{prop}

\begin{prop}
\label{p-60}
If $p$ and $q$ are two conditions in $\poset(\H)$ such that
$x_q\sqsubseteq x_p$ and every $\J\in\X_q$ has a $\K\subseteq\J$ in
$\X_p$, then $\poset(\H)\div{\qsep}\models q\le p$. 
\end{prop}

\begin{prop}
\label{p-64}
$p$ and $q$ are compatible in $\poset(\H)$ iff $x_p$ and $x_q$ are
comparable under end-extension and $x_p\cup
x_q\in\shade_\H(\J)$ for all $\J\in\X_p\cup\X_q$. 
\end{prop}

The following game is equivalent to the game $\ggen(M,y,\H,\C(\downcl\H),p)$
from\linebreak \cite[Definition~3.11]{Hir-comb}, except for the requirement that
$p_0=p$. It is used to establish the various properties of our forcing notion. 

\begin{notn}
\label{notn:generate}
For a centered subset $C$ of some forcing notion $P$, 
we let $\<C\>$ denote the ideal on $P$
generated by $C$.
\end{notn}

\begin{defn} 
\label{d-11}
For any $M\prec H_\kappa$, with $\H\subseteq\powcif\theta$ nonempty and in $M$, 
for $y\in\powcif\theta$ and  $p\in\poset(\H)\cap M$,
we  define the game $\ggen(M,y,\H,p)$ 
with players \emph{Extender} and \emph{Complete} of length~$\omega$.
Extender plays first and on move~$0$ must play $p_0$ so that
\begin{itemize}[leftmargin=*]
\item $p_0=p$.
\end{itemize}
On Extender's $k+1\Th$ move:
\begin{itemize}[leftmargin=*]
\item Extender plays $p_{k+1}\in\poset(\H)\cap M$ satisfying
  \begin{enumerate}[leftmargin=*, label=(\arabic*)]
  \item\label{item:14}  $p_{k+1}$ extends $p_k$,
  \item\label{item:15}  $x_{p_{k+1}}\setminus x_{p_k}
    \subseteq y\setminus\bigcup_{i=0}^k s_i$.
  \end{enumerate}
\end{itemize}
On Complete's $k\Th$ move:
\begin{itemize}[leftmargin=*]
\item Complete plays a finite $s_k\subseteq y$.
\end{itemize}
This game has three possible outcomes, determined as follows:
\begin{enumerate}[leftmargin=*, label=(\roman*), ref=\roman*, widest=iii]
\item\label{item:35} 
  Extender loses (i.e.~Complete wins) if
  $\<p_k:k<\omega\>\notin\Gen(M,\poset(\H))$,
\item\label{item:36}
  the game is drawn (i.e.~a tie) if $\<p_k:k<\omega\>\in\Genc(M,\poset(\H))$,
\item\label{item:7} Extender wins the game
if $\<p_k:k<\omega\>\in\Gen(M,\poset(\H))$ 
  but~\eqref{item:36} fails.
\end{enumerate}
\end{defn}

The game $\ggen(M,y,\H,p)$ is especially interesting for us 
because a draw in this game 
corresponds precisely with complete genericity.

\begin{prop}
\label{p-4}
Let\/ $p_k$ denote Extender's\/ $k\Th$ move in the game\/ 
$\ggen(M,\allowbreak y,\allowbreak \H,\allowbreak p)$.
Then the game results in a draw iff\/ $\<p_k:k<\omega\>
\in\Genc(M,\allowbreak\poset(\H),\allowbreak p)$.
\end{prop}

\begin{prop}
\label{p-21}
At the end of the game $\ggen(M,y,\H,p)$, for every $k<\omega$,
\begin{equation}
  \label{eq:25}
  \bigcup_{n<\omega}x_{p_n}\setminus x_{p_k}\subseteq y
  \setminus\bigcup_{i=0}^k s_i.
\end{equation}
\end{prop}

\begin{prop}
\label{p-24}
Suppose $\Phi$ is a nonlosing strategy for Complete in the game $\ggen(M,y,\H,p)$. 
Then Complete does not lose if it plays $s_k\supseteq\Phi(P_k)$, 
where $P_k$ is the position after Extender's $k\Th$ move. 
\end{prop}
\begin{pf}
This is by general principles. If Complete plays as in the hypothesis,
then it is restricting Extender's moves. Since the outcome of the game
is determined solely on Extender's sequence of moves, this is
beneficial to Complete. 
\end{pf}

Similarly:

\begin{prop}
\label{p-17}
For all $y\subseteq z$, if $\Phi$ is a nonlosing strategy for Complete in the game
$\ggen(M,z,\H,p)$ then $P_k\mapsto\Phi(P_k)\cap y$ 
is a nonlosing strategy for Complete in the game
$\ggen(M,\allowbreak y,\allowbreak \H,\allowbreak p)$. 
\end{prop}

Before proceeding, recall that whenever $(\H,\subseteqfnt)$ is
$\sigma$-directed, the ideal of noncofinal subsets of $\H$ forms a $\sigma$-ideal,
i.e.~is closed under countable unions (see
e.g.~\cite[Lemma~2.2]{Hir-comb}).

\begin{lem}\label{lem:stationary2}
Suppose $\H$ is a $\sigma$-directed subfamily of $\powcif\theta$ with
no stationary orthogonal set. 
For every countable $M\prec H_\kappa$ with $\H\in M$, 
if $\J\subseteq\powcif\theta$ in $M$ is cofinal in $(\H,\subseteqfnt)$, 
then $\sup(\theta\cap M)\in\shade_{\H}(\J)$.
\end{lem}
\begin{pf}
Let $Z$ be the set of all $\alpha<\theta$ such that $\alpha\notin\shade_\H(\J)$. 
Supposing towards a contradiction that the lemma fails,
 $\sup(\theta\cap M)\in Z$, and thus $Z$ is stationary because $Z\in M$. 
By the assumption on $\H$, $Z$ is not orthogonal to $\H$, 
say $x\in\powcif Z$ with $x\subseteq y$ for some $y\in\H$.
Since $\{z\in\J:y\subseteqfnt z\}$ is cofinal in $(\H,\subseteqfnt)$
as $\H$ is directed, there must exist a finite $s\subseteq x$ such
that $\{z\in\J:x\setminus s\subseteq z\}$ is cofinal because $\H$ is
$\sigma$-directed. We have now arrived at the contradiction that
$\alpha\in\shade_\H(\J)$ for all $\alpha\in z\setminus s$.
\end{pf}

The following corollary implies that 
$\poset(\H)$ forces a closed cofinal subset of $\theta$, although it
remains to show that $\poset(\H)$ does not collapse $\aleph_1$.

\begin{cor}
\label{c-4}
Assume $\H$ is as in lemma~\tu{\ref{lem:stationary2}}.
For every $\xi<\theta$,
\begin{equation}
  \label{eq:46}
  \D_\xi=\{p\in\poset(\H):\max(x_p)\ge\xi\}
\end{equation}
is a dense subset of $\poset(\H)$.
\end{cor}
\begin{pf}
Given $p\in\poset(\H)$, find a countable elementary 
$M\prec H_\kappa$ with $\H,p,\xi\in M$, set  $\delta=\sup(\theta\cap M)$. 
For each $\J\in\X_p$, let $\K_\J=\{y\in \J:x_p\subseteq y\}$. 
Each $\K_\J$ is cofinal as $x_p\in\shade_\H(\J)$, 
and thus by lemma~\ref{lem:stationary2},
$\delta\in\shade_\H(\K_\J)$ for all $\J\in\X_p$.
This implies that $x_p\cup\{\delta\}\in\shade_\H(\J)$ for all $\J\in\X_p$, 
and therefore $q=(x_p\cup\{\delta\},\X_p)\in\poset(\H)$.
Since $\delta>\xi$, the proof is complete.  
\end{pf}

\begin{cor}
\label{c-1}
If $\H\subseteq\powcif\theta$ is $\sigma$-directed under
$\subseteqfnt$ and has no stationary orthogonal set, then
$\poset(\H)\forces C_{\dot G_{\poset(\H)}}$ is a closed cofinal subset of $\theta$.
\end{cor}

\begin{termgy}
Henceforth, we let $\varphi(\theta,\H)$ abbreviate the statement: \emph{$\H$
  is a $\sigma$-directed subfamily of $(\powcif\theta,\subseteqfnt)$
  with no stationary orthogonal set}. 
\end{termgy}

\begin{cor}
\label{c-3}
Assume $\varphi(\theta,\H)$. 
Let $\H\in M\prec H_\kappa$, $y\in\powcif\theta$ and $p\in\poset(\H)\cap M$,
and let $p_k$ denote Extender's $k\Th$ move in the game $\ggen(M,y,\H,p)$. 
Suppose that Extender does not lose the game.
Then the following are equivalent\tu:
\begin{enumerate}[label=\tu{(\alph*)}, ref=\alph*, widest=b]
\item\label{item:58} The game $\ggen(M,y,\H,p)$ is drawn.
\item\label{item:62}
  $\bigl(\bigcup_{k<\omega}x_{p_k}\cup\{\sup(\theta\cap M)\},
  \bigcup_{k<\omega}\X_{p_k}\bigr)\in\poset(\H)$.
\item\label{item:63} $\bigcup_{k<\omega}x_{p_k}\cup\{\sup(\theta\cap M)\}
  \in\shade_\H(\J)$ for all $\J\in\X_q$.
\end{enumerate}
\end{cor}
\begin{pf}
Put $\delta=\sup(\theta\cap M)$. 
By corollary~\ref{c-4}, 
for every $\xi<\theta$ in $M$, $\D_\xi\in M$ (cf.~equation~\eqref{eq:46}) is dense, and
thus $x_{p_k}\in\D_\xi$ for some $k$ since Extender did not lose.
Hence $\bigcup_{k<\omega}x_{p_k}$ is unbounded in~$\delta$. 

\eqref{item:58}\tar\eqref{item:62}: 
$\{p_k:k<\omega\}$ has a common extension, say $q$, by proposition~\ref{p-4}.
Since we must have $\delta\in x_{q}$, 
and obviously $\bigcup_{k<\omega}\X_{p_k}\subseteq\X_q$,
It clearly follows that the pair defined in~\eqref{item:62} is a
condition of $\poset(\H)$. 

The implication~\eqref{item:62}\tar\eqref{item:63} is trivial by
definition; the implication~\eqref{item:63}\tar\eqref{item:62} is
because the set is closed; and the implication~\eqref{item:62}\tar\eqref{item:58}
is an immediate consequence of proposition~\ref{p-4}.
\end{pf}

\begin{lem}
\label{l-29}
Let $\H\subseteq\powcif\theta$ be $\sigma$-directed. 
Suppose $M\prec H_\lambda$ is countable with $\H\in M$, and
$y\in\powcif\theta$ satisfies $x\subseteqfnt y$ for all $x\in\H\cap M$.
Then every\/ $p\in\poset(\H)\cap M$ and every dense
$D\subseteq\poset(\H)$ in $M$, has an extension $q$ of $p$ in
$D\cap M$ such that $x_q\setminus x_p\subseteq y$. 
\end{lem}
\begin{pf}
Suppose to the contrary that there is no $q\ge p$ in $D\cap M$ with $x_q\setminus x_p\subseteq y$. 
Define $\F$ as the set of all $K\in\powcif{H_\kappa}$ 
for which there exists $y_K\in\powcif\theta$ such that $x_p\subseteq y_K$,
there is no $q\ge p$ in $D$ with $x_q\subseteq y_K$, 
and $x\subseteqfnt y_K$ for all $x\in\H\cap K$. 
Then $\F\in M$. 

Take $K\in\powcif{H_\kappa}\cap M$. Since $\{z\in\H:x_p\subseteq z\}$
is cofinal by condition~\eqref{item:33} of the forcing notion, 
and since $\H$ is $\sigma$-directed,
there exists $z\in\H\cap M$ such that $x_p\subseteq z$ and $x\subseteqfnt z$ for
all $x\in\H\cap K$. Then there exists a finite $s\subseteq z\setminus
x_p$ such that $z\setminus s\subseteq y\cup x_p$.  
There can be no $q\ge p$ in $D$ with $x_q\subseteq z\setminus s$,
because otherwise by elementarity we could find such a $q\in D\cap M$,
contradicting our supposition. Hence $y_K=z\setminus s$ witnesses that $K\in\F$.

By elementarity, we have proved that $\F=\powcif{H_\kappa}$, 
and thus $\J=\{y_K:K\in\nobreak\F\}\subseteq\H$ is cofinal in
$(\H,\subseteqfnt)$ with $x_p\in\shade_\H(\J)$. Hence
\begin{equation}
  \label{eq:95}
  q=(x_p,\X_p\cup\{\J\})\in\poset(\H).
\end{equation}
Since $D$ is dense, there exists $q'\ge q$ in $D$. But then
$x_{q'}\in\shade_\H(\J)$, and in particular $x_{q'}\subseteq y_K$ for some
$K\in\F$, contradicting the choice of $y_K$.
\end{pf}

\begin{cor}
\label{c-7}
Let $\H$ be a $\sigma$-directed subfamily of $(\powcif\theta,\subseteqfnt)$.
Suppose $M\prec H_\lambda$ with $\H\in M$, $y\in\powcif\theta$ and $p\in\poset(\H)\cap M$.
If $x\subseteqfnt y$ for all $x\in\H\cap M$, 
then Extender has a nonlosing strategy in the game $\ggen(M,y,\allowbreak \H,p)$. 
\end{cor}
\begin{pf}
Let $(D_k:k<\omega)$ enumerate all of the dense subsets of
$\poset(\H)$ in $M$. Suppose the position in the game is
$(p_0,s_0),\dots,(p_k,s_k)$ after the $k\Th$ move. 
By the assumption on $y$,
\begin{equation}
  \label{eq:118}
  x\subseteqfnt y\setminus\bigcup_{i=0}^k s_i\espc\text{for all $x\in\H\cap M$}.
\end{equation}
On move $k+1$, by lemma~\ref{l-29}, 
Extender can thus play $p_{k+1}\ge p_k$ in $D_k\cap M$ such that
$x_{p_{k+1}}\setminus x_{p_k}\subseteq y\setminus\bigcup_{i=0}^k s_i$. 
Clearly then $\<p_k:k<\omega\>\in\Gen(M,\poset(\H))$. 
\end{pf}

\begin{lem}
\label{l-15}
Assume $\varphi(\theta,\H)$. 
Let $M\prec H_\lambda$ be  a countable elementary submodel with $\H\in M$, 
let $y\in\powcif\theta$ and let $p\in\poset(\H)\cap M$. 
If $x\subseteqfnt y$ for all\/~$x\in\H\cap M$, 
then Complete has a nonlosing strategy in the game $\ggen(M,y,\H,p)$.
\end{lem}
\begin{pf}
Suppose that $x\subseteqfnt y$ for all $x\in\H\cap M$. 
Then since $\H$ is $\sigma$-directed, we can find $z\in\H$ such that
$x\subseteqfnt z$ for all $x\in\H\cap M$. Therefore, by
proposition~\ref{p-17}, replacing $y$ with $y\cap z$ we can assume
that $y\in\downcl\H$. 

Set $\delta=\sup(\theta\cap M)$.
At the end of the game $\ggen(M,y,\H,p)$, where Extender has played
$p_k$ on its $k\Th$ move, 
we will set $x_q=\bigcup_{k<\omega}x_{p_k}\cup\{\delta\}$.
The aim of the Complete's strategy is to ensure that
$x_q\in\shade_\H(\J)$ for all $\J\in\X_{p_k}$, for all $k<\omega$. 

We know that $\bigcup_{k<\omega}\X_{p_k}$ will be countable, 
and thus we can arrange a diagonalization $(\J_k:k<\omega)$ in advance,
and since the $\X_{p_k}$'s will be increasing with $k$, 
we can also insist that $\J_k\in\X_{p_k}$ for all $k$. 
After Extender plays $p_k$ on move $k$, we take care of some $\J_k\in\X_{p_k}$
according to the diagonalization. Set 
\begin{equation}
\label{eq:42}
\K_k=\{z\in \J_k:x_{p_k}\subseteq z\}.
\end{equation}
Then $\K_k$ is $\subseteqfnt$-cofinal in $\H$ because $x_{p_k}\in\shade_\H(\J_k)$ by the
definition of the poset. 
Thus as $p_k\in M$, by lemma~\ref{lem:stationary2},
\begin{equation}
  \label{eq:80}
  \K'_k=\{z\in \K_k:\delta\in z\}
\end{equation}
is cofinal too. 

\begin{lemclaim}
There exists a finite $s_k\subseteq y$ 
such that $y\setminus s_k\in\shade_\H(\K'_k)$.
\end{lemclaim}
\begin{pf}
Since $y\in\downcl\H$, and $\H$ is directed, $\{z\in\K'_k:y\subseteqfnt z\}$
is cofinal. It then follows from the
fact that $\H$ is $\sigma$-directed that there exists a finite
$s_k\subseteq y$ such that $y\setminus s_k\in\shade_\H(\K'_k)$. 
\end{pf}

Complete plays $s_k$ as in the claim on its $k\Th$ move. 
This describes the strategy for Complete. 

If Extender loses then Complete wins, and 
thus we may assume that Extender does not lose.
Put $x_q=\bigcup_{k<\omega} x_{p_k}\cup\{\delta\}$  and $\X_q=\bigcup_{k<\omega}\X_{p_k}$.
It remains to show that the game is drawn, and thus it suffices to
show that $x_q\in\shade_\H(\J)$ for all $\J\in\X_q$ by corollary~\ref{c-3}. 
But every $\J\in\X_q$ appears as $\J_k$ for some $k$, and thus as 
$\bigcup_{n<\omega}x_{p_n}\setminus x_{p_k}\subseteq y\setminus s_k$ by proposition~\ref{p-21}, 
$y\setminus s_k\in\shade_\H(\K'_k)$ implies that $\{y\in\J_k:x_q\subseteq y\}$
is cofinal by equations~\eqref{eq:42}
and~\eqref{eq:80}, proving $x_q\in\shade_\H(\J_k)$. 
\end{pf}

The following lemma implies that our poset does not collapse $\aleph_1$.

\begin{lem}\label{lem:proper:b}
Assume that $\H$ is a $\sigma$-directed subfamily of
$(\powcif\theta,\subseteqfnt)$ with no stationary orthogonal set \tu(i.e.~$\varphi(\theta,\H)$\tu).
Then $\poset(\H)$ is completely proper.  
\end{lem}
\begin{pf}
Let $M\prec H_\lambda$ be countable with $\H\in M$. 
Since $\H$ is $\sigma$-directed there exists $y\in\H$
such that $x\subseteqfnt y$ for all $x\in\H\cap M$.
Let $p\in\poset(\H)\cap M$ be given. 
Then the hypotheses of corollary~\ref{c-7} and lemma~\ref{l-15} are
satisfied.
Therefore both Extender and Complete have nonlosing strategies in the
game $\ggen(M,y,\H,p)$. The game is played with both Extender and
Complete playing according to their respective strategies. Since the
game results in a draw, there exists $q\in\genc(M,\poset(\H),p)$ by
proposition~\ref{p-4}. This proves that $\poset(\H)$ is completely proper.
\end{pf}

\begin{termgy}
We let $\varphi_*(M,\H,y)$ abbreviate the statement:
\emph{$x\subseteqfnt y$ for all $x\in\nobreak\H\cap M$}. 
\end{termgy}

\begin{lem}
\label{l-14}
Let $(\vec\A,\idealmodp{\Omega})$ be a $\lambda$-properness parameter
determined by\/ $\Omega:\linebreak\injlim\A\allowbreak\to\powcnt{\powcif\theta}$ 
\tu(cf.~definition~\tu{\ref{ex:1})}. 
Assume $\varphi(\theta,\H)$. 
If every $M\in\injlim\A$ of positive rank 
with $\H\in M$ has a\/ $y\in\Omega(M)\cap\downcl\H$
satisfying $\varphi_*(M,\H,y)$,
then $\poset(\H)$ is  $(\vec\A,\idealmodp{\Omega})$-proper.
\end{lem}
\begin{pf}
This is proved by induction on $\rank(M)$, 
where $M$ is from the set of all $M\in\injlim\A$ with $\H\in M$.
The induction hypothesis is that for all $p\in\poset(\H)\cap\nobreak
M$,  
and all $X\in\idealmodp{\Omega}(M)$, for every
$y_M\in\downcl\H$ satisfying $\varphi_*(M,\H,y_M)$ and moreover
$y_M\in\Omega(M)$ when $\rank(M)>0$,
and every finite $t\subseteq y_M$,
there exists $q\in\nobreak\gen(M,\poset(\H),p)$ such that $\genmod(q)\cap X\in\idealmodp{\Omega}(M)$ and
\begin{equation}
x_q\setminus x_p\subseteq (y_M\setminus t)\cup\{\sup(\theta\cap M)\}.\label{eq:55}
\end{equation}
This will in particular entail that
$\poset(\H)$ is proper for the desired parameter, because by the
hypothesis on $\Omega$ there always exists such a $y_M$ and hence 
 $a=\H$ will witness that \pproper holds.

For the the base case $\rank(M)=0$, it suffices to 
show that every $p\in\poset(\H)\cap M$ and every $z\in\powcif\theta$  satisfying $\varphi_*(M,\H,z)$
have an $(M,\poset(\H))$-generic extension $q\ge\nobreak p$ 
with  $x_q\setminus\nobreak x_p\subseteq z\cup\{\sup(\theta\cap M)\}$. 
And Extender and Complete both have nonlosing strategies in
the game $\ggen(M,z,\H,p)$ by corollary~\ref{c-7} and
lemma~\ref{l-15}.  After the game is played according to
these respective strategies, with $p_k$ denoting Extender's $k\Th$ move,
we obtain $q\in\genc(M,\poset(\H),p)$ with
$x_q=\bigcup_{k<\omega}x_{p_k}\cup\{\sup(\theta\cap M)\}$ by corollary~\ref{c-3}. 
And then $x_q\setminus x_p\subseteq z\cup\{\sup(\theta\cap M)\}$ 
by proposition~\ref{p-21} with $k=0$ since $p_0=p$. 

Suppose now that $\rank(M)>0$ with $\H\in M$,
and we are given $p\in\poset(\H)\cap M$, $X\in\idealmodp{\Omega}(M)$
and $y_M\in\Omega(M)\cap\downcl\H$ satisfying $\varphi_*(M,\H,y_M)$
and a finite $t\subseteq y_M$. 
By going to a subset of $X$, 
we can assume that $X\cap
K\in\idealmodp\Omega(K)$ for all $K\in X$ by lemma~\ref{l-12}.
Moreover, since $y_M\in\Omega(M)$, by going to a tail of $X$, we can
assume that
\begin{equation}
\trsup_\theta(X)\subseteq y_M\setminus t.\label{eq:33}
\end{equation}
Let $(a_k:k<\omega)$ enumerate $M\cap H_\lambda$ and let
$(\xi_k:k<\omega)$ satisfy $\lim_{k\to\omega}\xi_k+1=\rank(M)$. 

The game $\ggen(M,y_M\setminus t,\H,p)$ 
shall be played with $(p_k,s_k)$ denoting the $k\Th$ move. 
Since $\varphi_*(M,\H,y_M\setminus t)$,
Complete has a nonlosing strategy in this game, 
which it plays by. 
After the $k\Th$ move has been played, 
we can find $K_k\in X$ such that $a_k,p_k \in K_k$, $\rank(K_k)\ge\xi_k$ and moreover
\begin{equation}
  \label{eq:57}
  \sup(\theta\cap K_{k})\notin\bigcup_{i=0}^k s_i.
\end{equation}
Since $y_{K_k}\in\downcl\H$, we can find a finite $u_k\subseteq y_{K_k}$ 
such that $y_{K_k}\setminus u_k\subseteq y_M$. 
Now by the induction hypothesis, there exists
$p_{k+1}\ge p_k$ in $\gen(K_{k},\poset(\H))$ such that 
\begin{equation}
Y_k=\genmod(p_{k+1})\cap X\in\idealmodp{\Omega}(K_{k}),\label{eq:34}
\end{equation}
and $x_{p_{k+1}}\setminus x_{p_k}
\subseteq \bigl(y_{K_k}\setminus \bigl(\bigcup_{i=0}^k s_i \cup t\cup u_k\bigr)\bigr)
\cup\{\sup(\theta\cap K_{k})\}
\subseteq y_M\setminus\bigl(\bigcup_{i=0}^k s_i \cup t\bigr)
\cup\{\sup(\theta\cap K_k)\}$. 
Then in fact $x_{p_{k+1}}\setminus x_{p_k}
\subseteq y_M\setminus\bigl(\bigcup_{i=0}^k s_i\cup t\bigr)$ 
by~\eqref{eq:33} and~\eqref{eq:57}, and thus  $p_{k+1}$ is a valid move for
Extender.

At the end of the game, since $(K_{k}:k<\omega)$ is cofinal in
$M\cap H_\lambda$ and each $p_{k+1}$ is $(K_{k},\poset(\H))$-generic,
it follows that the ideal $\<p_k:k<\omega\>$ is in $\Gen(M,\poset(\H))$, 
and thus Extender does not lose. 
Since Complete also does not lose, the conditions
$\{p_k:k<\omega\}$ have a common extension $q$. 
Now $\bigcup_{k<\omega}Y_k\cup\{K_{k}\}\subseteq\genmod(q)\cap X$, and clearly
$\bigcup_{k<\omega}Y_k\cup\{K_{k}\}\in\idealmodp{\Omega}(M)$.  
Moreover we can assume that
$x_q=\bigcup_{k<\omega}x_{p_k}\cup\{\sup(\theta\cap M)\}$ by
corollary~\ref{c-3}, 
and thus  $x_q\setminus x_p\subseteq (y_M\setminus t)
\cup\{\sup(\theta\cap\nobreak M)\}$, 
completing the induction.
\end{pf}

For definitions of $\mathbb D$-completeness 
we refer the reader to~\cite{Hir-comb} and/or \cite{hbst}. 
In the present paper we say that a poset is \deecmp, if it has 
a simply definable $\aleph_1$-complete completeness system.
To avoid a complicated proof we only prove that $\poset(\H)$ has a
simply definable $\aleph_0$-complete completeness system. 
If an $\aleph_1$-complete system is desired, one can use
the notion of a \emph{forward strategy} introduced there; in
particular, lemma~\ref{l-8} can be obtained as an application
of~\cite[Lemma~3.39]{Hir-comb}. 

\begin{lem}
\label{l-8}
Let $\H\subseteq\powcif\theta$ be $\sigma$-directed with no stationary
orthogonal set \tu(i.e. $\varphi(\theta,\H)$\tu).
Then $\poset(\H)$ is \deecmp.
\end{lem}
\begin{pf}[for $\aleph_0$-completeness]
The completeness system receives as input a countable $M\prec
H_\lambda$,  
a family $\H\subseteq\powcif\theta$ in $M$ and $p\in\poset(\H)\cap M$. 
We fix a suitably definable way of coding
\begin{itemize}
\item a subset $y_Z$ of $\theta\cap M$,
\item a partial function $\Phi_Z$ on $M$ with
  $\Phi_Z(a)\in\powfin{\theta\cap M}$ for all $a\in\dom(\Phi_Z)$,
\end{itemize}
by subsets $Z\subseteq M$. 
The second order formula $\varphi$ defining the family of subsets of
$\Gen(M,\poset(\H),p)$ is given by $\ulc$if
\begin{enumerate}[label=(\alph*), ref=\alph*, widest=b]
\item\label{item:163} $x\subseteqfnt y_Z$ for all $x\in\H$,
\save 
\end{enumerate}
there exists\footnote{Note that this is a second order quantifier,
  so that e.g.~the sequence $(p_k:k<\omega)$ need not be an element of $M$.}
a sequence $(p_k:k<\omega)$ of conditions in~$\poset(\H)$ 
and a sequence $(s_k:k<\omega)$ of finite subsets of $\theta$ such that
\begin{enumerate}[label=(\alph*), ref=\alph*, widest=b]
\restore
\item\label{item:2} $(p_k,s_k)$ is valid for move $k$ of
  the game  $\ggen(M,y_Z,\H,p)$,
\item\label{item:11} 
$\vec a=\bigl((p_0,s_0),\dots,(p_{k-1},s_{k-1}),p_k\bigr)\in\dom(\Phi_Z)$
  and $s_k\supseteq\Phi_Z(\vec a)\urc$.
\end{enumerate}
Thus the family coded by some $Z\subseteq M$ is
\begin{equation}
  \label{eq:2}
  \G_Z=\{G\in\Gen(M,\poset(\H),p):M\models\varphi(G,Z;\H,p)\}.
\end{equation}

First we verify $\aleph_0$-completeness. Let $Z_0,\dots,Z_{n-1}$ be given
subsets of $M$. We can assume without loss of generality that
condition~\eqref{item:163} holds for all $j=0,\dots,n-1$. 
The game $\ggen\bigl(M,\bigcap_{j=0}^{n-1}y_{Z_j},\H,p\bigr)$ 
is played with $(p_k,s_k)$ being the $k\Th$ move. 
By condition~\eqref{item:163}, $x\subseteqfnt\bigcap_{j=0}^{n-1} y_{Z_j}$
for all $x\in\H\cap M$, and hence by corollary~\ref{c-7} 
Extender has a nonlosing strategy in this game, which it plays by.
For each $j=0,\dots,n-1$, we recursively choose  for each $k<\omega$, 
$t_{jk}$ so that $\vec a_{jk}=\bigl((p_0,t_{j0}),\dots,(p_{k-1},t_{j(k-1)})\bigr)$ 
is a valid position in the game $\ggen(M,y_{Z_j},\H,p)$; 
its definition is $t_{jk}=\Phi_{Z_j}(\vec a_{jk})$.     
On move $k$, Complete plays
$s_k=\bigl(\bigcup_{j=0}^{n-1}t_{jk}\bigr)\cap\bigcap_{j=0}^{n-1}y_{Z_j}$,
which ensures that Extender's move $p_{k+1}$ in the former game
is also a valid move in each of the games $\ggen(M,y_{Z_j},\H,p)$ ($j=0,\dots,n-1$).
Let $G=\<p_k:k<\omega\>$. Then $G\in\Gen(M,\poset(\H),p)$ because Extender
does not lose. And thus for each $j$, 
$M\models\varphi(G,Z_j;\H,p)$ as witnessed by
$(p_0,t_{j0}),(p_1,t_{j1}),\dots$; hence,
$\bigcap_{j=0}^{n-1}\G_{Z_j}\ne\emptyset$
as wanted. 

For $\mathbb D$-completeness, it remains to find a $Z\subseteq M$ such
that $\G_Z\subseteq\Genc(M,\allowbreak\poset(\H),\allowbreak p)$. 
However, choosing any $y\in\powcif\theta$ satisfying
$\varphi_*(M,\H,y)$, 
Complete has a nonlosing
strategy $\Phi$ in the game $\ggen(M,y,\H,p)$ by lemma~\ref{l-15}. 
Find $Z\subseteq
M$ such that $y_Z=y$ and $\Phi_Z=\Phi$. 
Now suppose that $G\in\G_Z$, witnessed by $(p_k:k<\omega)$ and $(s_k:k<\omega)$. 
Then by~\eqref{item:2} and~\eqref{item:11}, 
and proposition~\ref{p-24}, 
Complete does not lose the game $\ggen(M,y,\H,p)$ where $(p_k,y_k)$ is played on move $k$. 
Since Complete does not lose, 
and we already know that $G\in\Gen(M,\poset(\H))$, 
we must have $G\in\Genc(M,\poset(\H))$. 
\end{pf}

In the case $\theta=\oone$, assuming $\ch$, our poset $\poset(\H)$
clearly has the $\aleph_2$-cc and thus does not collapse cardinals. 
However, if we want to avoid using an inaccessible cardinal, 
we need that iterated forcing constructions using our poset also have the
$\aleph_2$-cc, which is not in general preserved under countable
support iterations. The usual approach in this
situation is to use the \emph{properness isomorphism condition}, and
apply the theory from~\cite[Ch.~VIII,~\Section2]{MR1623206}. By the
properness isomorphism condition, we mean the $\aleph_2$-pic there;
and there is a theorem that under $\ch$, a countable support
iteration of length at most $\omega_2$ of posets satisfying the
$\aleph_2$-pic has the $\aleph_2$-cc.\label{pic} As an alternative to
using lemma~\ref{l-13}, one can
always iterate up to a strongly inaccessible cardinal $\mu$ since our
posets will all have the~$\mu$-cc. 

We will not give the actual definition of the properness isomorphism condition here, but
instead refer the reader to either Shelah's book,~\cite{Hir-comb} or~\cite{hbst}. We
also do not provide a proof of the following lemma, because as is
usual with this property, it a straightforward modification of the
proof of properness. One can also 
obtain lemma~\ref{l-13} by applying~\cite[Corollary~3.54.1]{Hir-comb}.

\begin{lem}
\label{l-13}
Assume that $\H$ is a $\sigma$-directed subfamily of
$(\powcif\oone,\subseteqfnt)$ with no stationary orthogonal set 
\tu(i.e.~$\varphi(\oone,\H)$\tu).
Then $\poset(\H)$ satisfies the properness isomorphism condition. 
\end{lem}

\subsection{Forcing notion for shooting 
non-clubs}
\label{sec:auxill-forc-noti}

The following is perhaps the most natural forcing notion for forcing
an uncountable set locally in some $\sigma$-directed subfamily of
$(\powcif S,\subseteqfnt)$, for some set~$S$. 

\begin{defn}
\label{d-19}
For $\H\subseteq\power(S)$,
let $\spo(\H)$ be the poset consisting of all pairs $p=(x_p,\X_p)$ where
\begin{enumerate}[label=(\roman*), ref=\roman*, widest=iii]
\item\label{item:44} $x_p\in\shade(\H)$,
\item\label{item:45} $\X_p$ is a countable family of cofinal subsets
  of   $(\H,\subseteqfnt)$ with
  \begin{equation}
    \label{eq:84}
    x_p\in\shade_\H(\J)\espc\text{for all $\J\in\X_p$},
  \end{equation}
\save
\end{enumerate}
ordered by $q$ extends $p$ if
\pagebreak
\begin{samepage}
\begin{enumerate}[label=(\roman*), ref=\roman*, widest=iii]
\restore
\item \label{item:60} $x_q\supseteq x_p$,
\item \label{item:65} $\X_q\supseteq\X_p$.
\end{enumerate}
\end{samepage}
For an ideal $G\subseteq\spo(\H)$, 
we write $X_G=\bigcup_{p\in G}x_p$. We write $0_{\spo(\H)}$ for the
condition~$(\emptyset,\emptyset)$. 
\end{defn}

\begin{prop}
\label{p-51}
$\spo(\H)\forces\ulc X_{\dot G_{\spo(\H)}}$\tu{ is locally in }$\downcl\H\urc$.
\end{prop}

\begin{prop}
\label{p-56}
The class $\spo$ is provably equivalent to a $\Delta_0$-formula.
\end{prop}

\begin{prop}
\label{p-61}
If $p$ and $q$ are two conditions in $\spo(\H)$ such that
$x_q\subseteq x_p$ and every $\J\in\X_q$ has a $\K\subseteq\J$ in
$\X_p$, then $\spo(\H)\div{\qsep}\models q\le p$. 
\end{prop}

The only significant difference between definition~\ref{d-19} and the forcing
notion called ``$\spo(\H)$'' in~\cite[Definition~3.1]{Hir-comb}, which
is itself very closed based on the original forcing notion from~\cite{MR1441232},
is that condition~\eqref{item:60} is not required to be end-extension
as in $\poset(\H)$. This weakens the compatibility relation as follows.

\begin{prop}
\label{p-50}
$p$ and $q$ are compatible in $\spo(\H)$ iff $x_p\cup
x_q\in\shade_\H(\J)$ for all $\J\in\X_p\cup\X_q$. 
\end{prop}

Condition~\eqref{item:60}, however, leaves many properties
unaffected. For example, the following facts can be established with
exactly the same proofs as in~\cite{Hir-comb}. Assume for now that
$\H$ is a $\sigma$-directed subfamily of $(\powcif S,\subseteqfnt)$,
where $S$ is some uncountable set.

\begin{lem}
\label{l-24}
$\spo(\H)$ is $\aalpha$-proper.
\end{lem}

\noindent Note that, unlike with $\poset(\H)$, 
we do not need any additional requirements on $\H$ for properness as
in lemma~\ref{lem:proper:b}.

\begin{lem}
\label{l-39}
$\spo(\H)$ is \deecmp.
\end{lem}

\begin{lem}
\label{l-9}
$\spo(\H)$ satisfies the properness isomorphism condition. 
\end{lem}

\begin{lem}
\label{l-7}
If $S$ cannot be decomposed into countably many pieces that are
orthogonal to $\H$, then $\spo(\H)$ forces that $X_{\dot
  G_{\spo(\H)}}$ is uncountable.
\end{lem}

The following is established in~\cite{Hir-comb}. 

\begin{lem}
\label{l-70}
Let $\theta$ be an ordinal of uncountable cofinality, and let $\H$ be
a $\sigma$-directed subfamily of $(\powcif\theta,\subseteqfnt)$. 
Let $S\subseteq\theta$ be stationary.
If $S$ has no stationary subset orthogonal to $\H$, then
$\spo(\H)$ forces that $X_{\dot G_{\spo(\H)}}\cap S$ is stationary. 
\end{lem}

\section{Absolute antichains}
\label{sec:absolute-antichains}

Suppose $\H$ is a subfamily of $\powcif\theta$. 
Suppose that $W$ is an outer model of $V$. 
Since $\poset$ and~$\spo$ are considered as classes, we can interpret $\poset(\H)$ and $\spo(\H)$ in $W$. 
And by propositions~\ref{p-28} and~\ref{p-56}, respectively, we have
\begin{equation}
  \label{eq:70}
  \poset(\H)^V\subseteq\poset(\H)^W\And\spo(\H)^V\subseteq\spo(\H)^W,
\end{equation}
and thus for $\O=\poset,\spo$, $\O(\H)^V$ is a suborder of
$\O(\H)^W$ (recall that  $P$ is a \emph{suborder} of $Q$ when~$P\subseteq Q$ and
$p\len P q\iff p\len Q q$ for all $p,q\in P$).
Since $(x_p\cup x_q,\X_p\cup\X_q)$ is a common extension of any two
compatible conditions $p$ and $q$, it follows that `$p\perp q$' is
absolute between $V$ and $W$, for either forcing notion. Here we are
interested in having $\O(\H)^V$ generically included in $\O(\H)^W$---see
definition~\ref{d-16} where we obtain an approximation to this
property---and therefore we are interested in the upwards absoluteness
of $\ulc A$ is a maximal antichain of $\O(\H)\urc$ for both classes of
forcing notions $\O=\poset,\spo$. 

This is a familiar scenario. The concept of \emph{Souslin
  forcing} was introduced in~\cite{MR973109}, concerning a certain
class of forcing notions that can be represented as definable
subsets of the reals, or more generally, definable over $(H_{\aleph_1},\in)$. 
These Souslin forcing notions can thus be
interpreted in any outer model, and they enjoy many nice absoluteness
properties, which are particularly useful in the iteration of Souslin
forcing notions.
For example, the maximality of antichains of Souslin ccc
forcing notions is upwards absolute. In our case, say for
$\theta=\oone$ and assuming $\ch$, 
our forcing notions are $\aleph_2$-cc and representable as subsets of
$\power(\oone)$, and simply definable over $(H_{\aleph_2},\in)$ by
propositions~\ref{p-28} and~\ref{p-56}, respectively; 
and we also want to establish the upwards
absoluteness of antichains. However, in the present case we shall rely
on combinatorial arguments rather than absoluteness results of second
order arithmetic. In the process, we shall observe that $\spo(\H)$ and
$\poset(\H)$ have very nice properties, such as commutativity, 
that are typically associated with certain Souslin forcing notions.

\subsection{Embeddings}
\label{sec:embeddings}

We write $P\embeds Q$ to specify 
that a forcing notion $P$ \emph{generically embeds} into a forcing notion
$Q$, which as usual we mean that for all $G\in\Gen(V,Q)$, 
$V[G]\models\ulc\Gen(V,P)\ne\emptyset\urc$, i.e.~every generic for $Q$
induces a generic for $P$. A \emph{generic embedding} between two forcing
notions has the usual meaning, 
i.e.~they are called \emph{complete embeddings} 
in~\cite[Ch.~VII,~\Section7]{MR597342}. 
We write $P\cong Q$ to indicate that $P$ and $Q$ are isomorphic as
forcing notions, i.e.~$P\embeds Q$ and $Q\embeds P$. 
Recall that $P\embeds Q$ iff there exists a
\emph{generic embedding} $e:P\div{\sim_{\sep}}\to Q$, where the
separative quotient is indicated in the domain.
Also recall that when $P$ is a suborder of $Q$,
the inclusion map $i:P\to Q$ is a generic embedding iff $i[A]$ is a maximal antichain of
$Q$ for every maximal antichain $A$ of~$P$. We want to generalize the
notion $P\embeds Q$, where $Q$ is allowed to be outside of some universe. 

\fxnote{1. Mention somewhere that $V$ is some special implicit
  reference ground model? 2. Define $\Gen(M)$ so that does not matter
  outside $M$.}
\begin{defn}
\label{d-16}
Let $M$ be a model (typically transitive), 
and let $P$ and $Q$ be forcing notions with $P\in M$.
We say that \emph{$P$ generically embeds into $Q$ over~$M$}
if for all $G\in\Gen(V,Q)$, $V[G]\models\ulc\Gen(M,P)\ne\emptyset\urc$. 
We write $P\embeds[M] Q$.  We say that $P$ is \emph{generically
  included in $Q$ over $M$} if $P$ is a suborder of $Q$ that
generically embeds over $M$. We write $P\imbeds[M] Q$.                     
And $e:P\to Q$ is a \emph{generic embedding over $M$} if it is 
order preserving, i.e.~$q\ge p$ implies $e(q)\ge e(p)$, 
$q\perp p$ implies $e(q)\perp e(p)$
and for every maximal antichain $A$ of~$P$ in $M$, $e[A]$ is a maximal
antichain of $Q$. 
\end{defn}

Thus $V\models\ulc P\embeds[V] Q\urc$ iff $V\models\ulc P\embeds Q\urc$. 

\begin{prop}
\label{p-55}
$P\imbeds[M] Q$ iff the inclusion map $i$ is a generic embedding over~$M$. 
\end{prop}

\begin{notn}
\label{o-6}
For a separative poset $P$, let $\overbar P$ denote its completion.
\end{notn}

\begin{lem}
\label{p-32}
Let $M$ be a transitive model of enough of $\zfc$. 
Then $P\embeds[M] Q$ iff there exists $e:P \to\overline{Q\div {\qsep}}$
such that $e$ is a generic embedding over $M$. 
\end{lem}
\begin{pf}
If $e:P \to\overline{Q\div{\sim_{\tu{sep}}}}$ is a generic embedding
over $M$, and $G\in\Gen(V,\allowbreak\overline{Q\div{\qsep}})$, then
$e\inv[G]\subseteq P$ is downwards closed and upwards linked (i.e.~pairwise
compatible), and intersects every maximal antichain of $P$ in $M$.
It thus follows from well known facts that $e\inv[G]$ generates a
member of $\Gen(M,P)$. Thus
$P\embeds[M]\overline{Q\div{\qsep}}\cong Q$. 

Conversely, if $P\embeds[M] Q$ then there is a $Q$-name $\dot G$ for a
member of $\Gen(M,P)$. Then $e:P\to\overline{Q\div{\qsep}}$ defined by
\begin{equation}
  \label{eq:123}
  e(p)=\|p\in\dot G\|
\end{equation}
is a generic embedding over $M$. 
\end{pf}

\begin{prop}
\label{l-54}
Suppose that $P,Q\in M$ and $M$ is a model of enough of $\zfc$.
If $e:P\to Q$ is a generic embedding over $M$,
and $H\in\Gen(M,Q)$, then $e\inv[H]\in\Gen(M,P)$. 
\end{prop}

Recall the following basic fact. 

\begin{prop}
\label{p-54}
If $e:P\to Q$ is a generic embedding, then every $q\in Q$ has a $p_q\in
P$ 
such that $e(p')$ is compatible with $q$  for all $p'\ge p_q$. 
\end{prop}

\begin{lem}
\label{l-71}
Let $P$ and $Q$ be separative forcing notions with $P\in\nobreak M$,
where $M$ is a transitive model of enough of $\zfc$. 
If $e:P\to Q$ is a generic embedding over $M$, then so is
$\tilde e:P\div{\qsep}\to Q\div{\qsep}$ given by $\tilde e([p])=[e(p)]$. 
\end{lem}
\begin{pf}
First we observe that if $P\div{\qsep}\models p\ge q$ then
$Q\div{\qsep}\models e(p)\ge e(q)$. For suppose to the contrary that
$Q\div{\qsep}\models e(p)\ngeq e(q)$. Then there exists $r\ge e(p)$ 
in~$Q$ that is $Q$-incompatible with $e(q)$. In $M$, let $A$ be a
maximal antichain with $q\in A$. 
Then since $e[A]$ is a maximal antichain, there exists $q'\in A$ such
that $e(q')$ is $Q$-compatible with $r$. If $q'\perp q$ then $q'\perp p$
as $P\div{\qsep}\models p\ge q$, and hence $e(q')\perp e(p)$ implies
$e(q')\perp r$. But then $q'\in A$ implies $q'=q$, contradicting that
$e(q)\perp r$. 

The preceding observation obviously implies that $\tilde e$ is well
defined and order preserving. That $\tilde e$ preserves maximal
antichains follows immediately from the fact that $p\perp q$
implies $[p]\perp[q]$. 
\end{pf}

\begin{prop}
\label{p-62}
Let $P$ and $Q$ be separative forcing notions with $P\in\nobreak M$,
where $M$ is a transitive model of enough of $\zfc$. 
If $e:P\to Q$ is a generic embedding over $M$, then so is
$\bar e:\overbar P\to \overbarg Q$ given 
by $\bar e(\bar p)=\ifm\{e(p):p\in P$\tu, $p\ge\bar p\}$.
\end{prop}

The notion of a \emph{projection} is used in~\cite{hbst} as a map from
$Q$ into $P$ witnessing $P\embeds Q$. We weaken the requirements on
projections here for brevity, but only use them as inverses of generic
embeddings (noting that it would have been better to do it the former way).

\begin{defn}
\label{d-25}
We say that
$\pi:Q\to P$ is a \emph{projection} if $\pi$ is an order preserving surjection.
\end{defn}

\begin{defn}
\label{d-30}
Let $\kappa$ be a cardinal. A forcing notion $P$ is said to be
\emph{$\kappa$-semi\-com\-plete} if every $A\subseteq P$ of cardinality $|A|<\kappa$,
with an upper bound in $P$, has a minimal upper bound in $P$. It is
\emph{semicomplete} if it is $\kappa$-semicomplete for all cardinals $\kappa$. 
\end{defn}

\noindent In the case were $P$ is a poset (and not just a quasi
order), minimal upper bounds are suprema. Then semicomplete is
precisely the order theoretic notion of a complete semilattice. 
Also note that complete semilattices always have a minimum element,
namely $\spr\emptyset$. Recall that a poset $P$ is \emph{pointed} if it
has a minimum element, which denote as $0_P$. 

\begin{example}
\label{x-10}
\emph{$\poset(\H)$ and $\spo(\H)$ are both complete semilattices}. 
If $A\subseteq\spo(\H)$ and $a\le p$ for all $a\in A$, 
then clearly $\spr A=\bigl(\bigcup_{a\in A}x_a,\bigcup_{a\in
  A}\X_a\bigr)$. 
Similarly for $\poset(\H)$, but taking the closure of $\bigcup_{a\in A}x_a$.
\end{example}

\begin{defn}
\label{d-8}
Recall that a subset  $A$ of a poset $P$ \emph{upwards
  order closed} if whenever $B\subseteq A$ is nonempty,
if $a=\spr B$ exists when taken in $P$, then $a\in A$. 
\end{defn}

\noindent The above notion is not to be confused with an \emph{upwards
  closed} subset, also called an \emph{upper set}. 

Recall that $e:P\to q$ is an \emph{order embedding} between to quasi
orders if it is both order preserving and reflecting, 
i.e.~$p\le q$ iff $e(p)\le e(q)$ for all $p,q\in P$. For a poset, this
means that $e$ is isomorphic to its range. 

\begin{lem}
\label{p-22}
Suppose that $P$ is a pointed poset and $Q$ is a complete semilattice. 
If $e:P\to Q$ is an order embedding 
with an upwards order closed range, then it has a
projection $\pi:Q\to P$ for a left inverse, given by 
\begin{equation}
  \label{eq:58}
  \pi(q)=\spr\{p\in P:e(p)\le q\}.
\end{equation}
\end{lem}
\begin{pf}
To check that the supremum always exists, take $q\in Q$.
Put $A=\{p\in P:e(p)\le q\}$. In the case $A=\emptyset$,
the supremum is $0_P$. Otherwise, since $q$ is an upper bound of $e[A]$,
$a=\spr e[A]=\spr_{p\in A}e(p)$ exists in $Q$ since it is a complete
semilattice; and then since $e[P]$ is
upwards order closed, there exists $p'\in P$ such that $e(p')=a$. 
Then $p'$ is an upper bound of $A$ since $e$ is order reflecting. 
And if $r\in P$ is an upper bound of $A$, then $e(r)$ is an upper
bound of $e[A]$ since $e$ is order preserving, and thus $e(p')=a\le e(r)$
implies $p'\le r$, proving that $p'$ is the least upper bound.

It is clear that $\pi$ is an order preserving left inverse of $e$. 
(And obviously $\pi$ is a surjection when it is a left inverse.)
\end{pf}

\begin{example}
\label{x-9}
Now let us see how this applies to say $\poset(\H)$. Let $W$ be a
transitive outer model of $V$ which has the same countable sequences
of ordinals as $V$.  $\poset(\H)$ is a complete semilattice in $V$,
and thus, in $W$, $\poset(\H)^V$ is also a complete semilattice
by the assumption on $W$, because  in example~\ref{x-10} we showed
that the suprema are given by countable unions which thus remain in
$V$. The fact that suprema remain in $V$, also implies that
$\poset(\H)^V$ is upwards order closed in $\poset(\H)^W$. 
Therefore, 
lemma~\ref{p-22} applies in $W$ to the identity mapping
$i:\poset(\H)^V\to\poset(\H)^W$, yielding a projection
$\pi:\poset(\H)^W\to\poset(\H)^V$ in $W$ that is the identity on
$\poset(\H)^V$. Exactly analogous facts hold for $\spo(\H)$. 
\end{example}

\fxnote{This should all be fixed up.}

\begin{notn}
\label{o-7}
For any map $e:P\to Q$ let $e^*$ be the corresponding mapping of
$P$-names to $Q$-names
(cf.~e.g.~\cite[Ch.~VII,~7.12]{MR597342}).\fxnote{Check citation.}
\end{notn}

\begin{prop}
\label{p-36}
Let $P\subseteq Q$. 
If $i:P\to Q$ is the inclusion map then $i^*$ is an inclusion map. 
\end{prop}

\fxnote{original in comment}

\begin{prop}
\label{p-34}
Let $M\prec H_\lambda$ with $P,Q\in M$.
Suppose that $e:P\to Q$ is a generic embedding in
$M$\footnote{\emph{Not} over $M$.} 
and that equation~\eqref{eq:58} defines $\pi:Q\to P$ as a left inverse
of $e$. If $q\in\gen(M,Q)$, then $\pi(q)\in\gen(M,P)$\tu; and if
$q\in\genc(M,Q)$, then $\pi(q)\in\genc(M,P)$. 
\end{prop}

\fxnote{Should be about Gen not gen. Then gen as a corollary.}

\begin{lem}
\label{l-53}
Let $M\prec H_\lambda$ with $P,Q\in M$
and $\lambda$ sufficiently large and regular.
Suppose that $e:P\to Q$ is a generic embedding in $M$ 
and equation~\eqref{eq:58} defines a left inverse $\pi:Q\to
P$ of $e$,
$\varphi(v_0,\dots,v_{n-1})$ is a formula which is absolute for
transitive models of $\zfc$ and $\dot x_0,\dots,\dot x_{n-1}\in M$ are
$P$-names. Then for all $q\in\gen(M,Q)$, $\pi(q)\in\gen(M,P)$ and
\begin{equation}
  \label{eq:115}
  \pi(q)\forces\varphi(\dot x_0,\dots,\dot x_{n-1})
  \Iff q\forces\varphi(e^*(\dot x_0),\dots,e^*(\dot x_{n-1})).
\end{equation}
\end{lem}
\begin{pf}
See e.g.~\cite[Ch.~{VII},~7.13]{MR597342}.
\end{pf}

\fxnote{original in comment}

Next we recall some basic forcing facts on maximal antichains and
generic embeddings. 

\begin{prop}
\label{p-39}
For $A\subseteq P\gdot \dot Q$, let $A\div P$ be the $P$-name 
$\{\dot q:(p,\dot q)\in A$, $p\in \dot G_P\}$ for a subset of $\dot Q$. 
Then the following are equivalent\tu:
\begin{enumerate}[label=\tu{(\alph*)}, ref=\alph*, widest=b]
\item\label{item:85}  $A\subseteq P\gdot \dot Q$ is a maximal antichain.
\item\label{item:87} $P\forces A\div P$ is a maximal antichain of  $\dot Q$. 
\end{enumerate}
\end{prop}

\begin{prop}
\label{p-58}
Suppose $R\forces \dot Q\embeds\dot P$. 
Then $\dot P\div\dot Q\cong (R\gdot\dot P)\div(R\gdot\dot Q)$. 
\end{prop}

\begin{remark}
\label{r-9}
Both sides of the equivalence in proposition~\ref{p-58} are
$R\gdot\dot Q$-names, and thus the equivalence should of course be
interpreted as $R\gdot\dot Q\forces \dot P\div\dot Q\cong(R\gdot\dot
P)\div(R\gdot\dot Q)$. We shall apply the equivalent statement
\begin{equation}
  \label{eq:125}
  R\gdot\dot P\cong R\gdot\dot Q\gdot(\dot P\div\dot Q).
\end{equation}
\end{remark}

\begin{notn}
\label{o-11} 
For a forcing notion $P$ and $p\in P$, we let $P_p$ denote the
principle filter $\{q\in P:q\ge p\}$. 
\end{notn}

\begin{prop}
\label{p-52}
Assume $Q\forces P\embeds[V]\dot R$. Then $P\embeds Q\gdot\dot R$. 
Moreover if $Q$ has a minimum element $0_Q$ and 
\tu{$Q\forces\ulc \dot e:P\to\dot R$ is a generic embedding over $V\urc$},
then $p\mapsto (0_Q,\dot e(p))$ defines a generic embedding\tu; 
hence, if $A\subseteq P$ is a maximal antichain then
\begin{equation}
  \label{eq:138}
  (Q\gdot\dot R)\div P
  \cong\coprod_{p\in A}\bigl(Q\gdot\dot R_{\dot e(p)}\bigr)\div P_{p}. 
\end{equation}
\end{prop}
\begin{pf}
Let $G\in\Gen(V,P)$ and $H\in\Gen(V[G],\dot R[G])$. 
In $V[G][H]$, $\Gen(V,P)\allowbreak\ne\emptyset$, because $P\embeds[V]\dot R[G]$.
Hence $P\embeds Q\gdot\dot R$ by definition. That the defined mapping
is a generic embedding is immediate from proposition~\ref{p-39}. 
\end{pf}

The following lemma well known, at least for the case $M=V$. 

\begin{lem}
\label{l-38}
Let $(P_\xi,\dot Q_\xi:\xi<\alpha)$ and $(P'_\xi,\dot Q'_\xi:\xi<\alpha)$ 
both be iterated forcing constructs 
with resulting forcing notions $P_\alpha$ and $P'_\alpha$, respectively\tu; 
and let $M$ be a transitive model of enough of $\zfc$. 
If $P_\xi\imbeds[M] P'_\xi$ for all $\xi<\alpha$, then $P_\alpha\imbeds[M] P'_\alpha$. 
\end{lem}

\begin{defn}
\label{d-6}
When we say that a forcing notion $P$ is \emph{densely included} in a
forcing notion $Q$, we mean that $P$ is a predense suborder of $Q$. 
\end{defn}

\begin{prop}
\label{p-68}
If $P$ is a suborder of $Q$ and $P\cong Q$ then $P$ is densely
included in $Q$. 
\end{prop}

\begin{lem}
\label{l-33}
Let $(P_\xi,\dot Q_\xi:\xi<\alpha)$ and $(P'_\xi,\dot Q'_\xi:\xi<\alpha)$ 
both be iterated forcing constructs 
with resulting forcing notions $P_\alpha$ and $P'_\alpha$,
respectively.
If $P_\xi$ is densely included in $P'_\xi$ for all $\xi<\alpha$, then
$P_\alpha\cong P'_\alpha$. 
\end{lem}

The next lemma is important for computing quotients. 

\begin{lem}
\label{l-55}
Let $O$ and $P\gdot\dot Q$ be forcing notions 
such that $O\embeds P\gdot\dot Q$.
Suppose that for all $G\in\Gen(V,O)$, 
there exists $f:(P\gdot\dot Q)\div G\to P$ such that 
\begin{enumerate}[label=\tu{(\alph*)}, ref=\alph*, widest=b]
\item\label{item:146} $f(p,\dot q)\ge p$, 
\item\label{item:144} $(p',\dot q)\in(P\gdot\dot Q)\div G$ for all
  $p'\ge f(p,\dot q)$,  
\item\label{item:145} 
$(f(p,\dot q),\dot q)$ and $(f(p',\dot q'),\dot q')$ 
  are $(P\gdot\dot Q)$-compatible 
  whenever $f(p,\dot q)$ and $f(p',\dot q')$ are compatible.
\end{enumerate}
Then $(P\gdot\dot Q)\div O\cong P$. 
\end{lem}
\begin{pf}
Let $G\in\Gen(V,O)$. 
Let $D=\{(f(p,\dot q),\dot q):(p,\dot q)\in(P\gdot\dot Q)\div G\}$. 
Then $D$ is a dense subset of $(P\gdot\dot Q)\div G$
by~\eqref{item:146} and~\eqref{item:144}. 
And by~\eqref{item:144} and~\eqref{item:145}, 
for all $(p,\dot q),(p',\dot q')\in(P\gdot\dot Q)\div G$, 
$(f(p,\dot q),\dot q)$ is $(P\gdot\dot Q)\div G$-compatible 
with $(f(p',\dot q'),\dot q')$ iff $f(p,\dot q)$ is $P$-compatible 
with $f(p',\dot q')$. Therefore, for all $(p,\dot q),(p',\dot q')\in D$,
$(p,\dot q)$ and $(p',\dot q')$ are $(P\gdot\dot Q)\div G$-compatible
iff $p$ and $p'$ are compatible. Hence, $(p,\dot q)\mapsto p$ is an
isomorphism between $D\div{\qsep}$ and $P\div{\qsep}$. Thus, as $D$ is
dense, $((P\gdot\dot Q)\div G)\div{\qsep}\cong D\div{\qsep}\cong P\div{\qsep}$.
Since $G$ is arbitrary, this proves that $O\forces (P\gdot\dot Q)\div
O\cong P$, as required.
\end{pf}

\subsection{Analysis of $\poset(\H)$ and $\spo(\H)$}
\label{sec:analysis-poseth-spoh}

To obtain generic embeddings of e.g.~$\poset(\H)^V$ into $\poset(\H)^W$,
we shall analyze the maximal antichains of $\poset(\H)$ and $\spo(\H)$. 
Clearly, if $A\subseteq\poset(\H)$ is a maximal antichain then
$\pi[A]=\{x_p:p\in A\}$ must be predense in $(\H,\subseteq)$; 
however, it need not form an antichain. 
For example, suppose that $\H$ is a $P$-ideal 
(and thus closed under addition of finite sets),
$y\in\shade_\H(\J)$ is closed and countable, 
and $\alpha>\max(y)$ is not in $\shade_\H(\J)$ 
but $\K$ is a $\subseteqfnt$-cofinal subset of $\H$
with $\alpha\in\shade_\H(\K)$. Then $(y,\{\J\})$ and
$(y\cup\{\alpha\},\{\K\})$ are incompatible conditions of $\poset(\H)$
even though $y\sqsubseteq y\cup\{\alpha\}$. Indeed, in analyzing the
sets $\pi[A]$, the difficulty is when $y$ is in the set and we
want to determine whether some $z\sqsubset y$ is also present in $\pi[A]$. 

For any $\H\subseteq\powcif \theta$, the following auxiliary family of
countable subsets of $\H$ allows us to analyze the maximality of
antichains in $\poset(\H)$ and $\spo(\H)$. 

\begin{defn}
\label{d-15}
Let $S$ be a set, and $\H\subseteq\power(S)$. 
Define a subcollection $\Psi(\H)\subseteq\powcnt{\H}$ of all $Z\in\powcnt{\H}$ for which
there exists $y\in\H$ such that every finite $s\subseteq y$ has a
finite $A_s\subseteq Z$ such that 
\begin{equation}
  \label{eq:108}
  \bigcup Z\setminus A_s \subseteq y\setminus s.
\end{equation}
\end{defn}

\noindent In particular, whenever $\H$ is an ideal,  $\bigcup Z\in\H$ for all $Z\in\Psi(\H)$.

\begin{lem}
\label{l-35}
Let $\H$ be a $\sigma$-directed subfamily of
$(\powcif S,\subseteqfnt)$.
Then $\Psi(\H)$ is a  $P$-ideal on $\H$.
\end{lem}
\begin{pf}
Assume $\H$ is $\sigma$-directed.
To verify that $\Psi(\H)$ is $\sigma$-directed,
let $Z_n$ ($n<\omega$) enumerate a subset of $\Psi(\H)$,
with $y_n\in\H$ witnessing $Z_n\in\Psi(\H)$ for each $n<\omega$.
Since $\H$ is $\sigma$-directed,
there exists $z\in\H$ with $y_n\subseteqfnt z$ for all $n<\omega$.
Choose finite subsets $t_n\subseteq S$ ($n<\omega$) so that
\begin{enumeq}
\item\label{item:104} $y_n\setminus t_n\subseteq z$, 
\item\label{item:105} $\bigcup_{n<\omega}t_n\supseteq z$,
\end{enumeq}
and put $s_n=y_n\cap\bigcup_{i=0}^n t_i$ for each $n$.
Let $A_{s_n}$ be the finite subset from equation~\eqref{eq:108} so
that $\bigcup Z_n\setminus A_{s_n}\subseteq y_n\setminus s_n\subseteq z$. 
Putting $Y=\bigcup_{n<\omega}Z_n\setminus A_{s_n}$, $Z_n\subseteqfnt
Y$ for  all $n<\omega$. 
And $Y\in\Psi(\H)$, because for any finite $u\subseteq z$,  
we can find $n$ so that $u\subseteq \bigcup_{i=0}^n t_i$,
and then $\bigcup\bigl(Y\setminus\bigcup_{i=0}^n A_{y_i\cap u}\bigr)\subseteq
z\setminus u$, where each $A_{y_i\cap u}$ is from
equation~\eqref{eq:108} with $y:=y_i$ and $Z:=Z_i$.

Moreover, $\Psi(\H)$ is obviously downwards closed, 
and it is an ideal because $\H$ is directed.
\end{pf}

\begin{lem}
\label{l-36}
Let $\H$ be a downwards closed $\sigma$-directed subfamily of
$\powcif S$. Then for $\K\subseteq\H$, the following are equivalent\tu:
\begin{enumerate}[label=\tu{(\alph*)}, ref=\alph*, widest=b]
\item\label{item:106} There exists a countable decomposition of $\K$
  into pieces orthogonal to $\Psi(\H)$.
\item\label{item:107} There exists a countable family $\X$ of  cofinal
  subsets of $(\H,\subseteqfnt)$ such that  
  $\K\cap\bigcap_{\J\in\X}\shade_\H(\J)\subseteq\{\emptyset\}$.
\end{enumerate}
\end{lem}
\begin{pf}
\eqref{item:106}\tar\eqref{item:107}: Let $\K=\bigcup_{n<\omega}\L_n$
with each $\L_n\perp\Psi(\H)$. Observe that every $y\in\H$ has a
finite $s_{yn}\subseteq y$ such that 
$\downcl{(y\setminus s_{yn})}\cap\L_n\subseteq\{\emptyset\}$:
Otherwise, letting $(\alpha_k:k<\nobreak\omega)$ enumerate $y$, there exists $z_k\in\L_n$ 
with $\emptyset\ne z_k\subseteq y\setminus\{\alpha_0,\dots,\alpha_k\}$ 
for all $k<\omega$, and then $\{z_k:k<\omega\}\in\powcif{\L_n}\cap\Psi(\H)$, 
contrary to the fact that $\L_n\perp\Psi(\H)$. 
Now for each $n$, let $\J_n=\{y\setminus s_{yn}:y\in\H\}$. 
Then every $\J_n$ is a cofinal subset of $\H$ as $\H$ is downwards closed,
and $\K\cap\bigcap_{n<\omega}\shade_\H(\J_n)
\subseteq\K\cap\bigcap_{n<\omega}\downcl{\J_n}\subseteq\{\emptyset\}$.

\eqref{item:107}\tar\eqref{item:106}: Let $\X=\{\J_0,\J_1,\dots\}$ 
be as in~\eqref{item:107}. Since $\H$ is $\sigma$-directed,
the noncofinal subsets of $\H$ form a $\sigma$-ideal; therefore, 
each $n$ and each $y\in\H$ has a finite $s_{yn}\subseteq y$ 
such that $y\setminus s_{yn}\in\shade_\H(\J_n)$. Then putting
$\K_n=\{y\setminus s_{yn}:y\in\H\}$ ($n<\omega$) we get
$\K\cap\bigcap_{n<\omega}\downcl\K_n\subseteq\{\emptyset\}$. 
For each~$n$, let $\L_n=\K\setminus\downcl\K_n$. 
Then $\K\setminus\{\emptyset\}=\bigcup_{n<\omega}\L_n$, 
and each $\L_n\perp\Psi(\H)$,
because $\K_n$ is cofinal, and hence if there were a $Z\in\powcif{\L_n}\cap\Psi(\H)$
witnessed by some $y\in\H$, then for any $z\in\K_n$ with
$y\subseteqfnt z$ we arrive at the contradiction that $\L_n\cap\downcl
z\ne\emptyset$. 
\end{pf}

\begin{notn}
\label{o-8}
For any $\S\subseteq\power(\theta)$ and $x\subseteq\theta$,
denote $\S_x=\{y\in\S:x\sqsubseteq y\}$ and $\S_{[x]}=\{y\setminus x:y\in\S_x\}$.
Analogously, we denote $\S_{(x)}=\{y\setminus x:y\in\S\cap\upcl x\}$
\tu($\S\cap\upcl x=\{y\in\S:x\subseteq y\}$\tu). 
\end{notn}

\begin{prop}
\label{p-41}
If $\H$ is a $P$-ideal then $\H_{[x]}$ and $\H_{(x)}$ are both
$P$-ideals for all~$x$.
\end{prop}

\fxnote{Might be useful to have a proposition relating cofinal subsets
  of $\H$ with $\H_{[x]}$.}

\begin{remark}
\label{r-3}
$\H_{[x]}$ and $\H_{(x)}$ may fail to be $\sigma$-directed if $\H$ is
not a $P$-ideal, even if $\H$ is $\sigma$-directed. 
\end{remark}

Next we see how to `freeze' a maximal antichain $A$, so that
the statement $\ulc A$ is a maximal antichain of $\poset(\H)\urc$ is
upwards absolute. 

We are assuming, until~\Section\ref{sec:product-p-ideals}, that $\H$ is a
$P$-ideal on $\theta$.

\begin{cor}
\label{c-16}
Let $A$ be a maximal antichain of $\poset(\H)$, 
and suppose $W\supseteq V$ is an outer model with $(\powcif\theta)^V=(\powcif\theta)^W$. 
If for every $x\in\H$, either
\begin{enumerate}[label=\tu{(\alph*)}, ref=\alph*, widest=b]
\item\label{item:108} there exists in $V$,
  a countable decomposition of $\pi[A]_{[x]}$ into
  pieces orthogonal to $\Psi(\H_{[x]})$, or
\item\label{item:110} there does not exist in $W$, 
a countable decomposition of $\pi[A]_{[x]}$ into
  pieces orthogonal to $\Psi(\H_{[x]})$,
\end{enumerate}
then $A$ is a maximal antichain of $\poset(\H)^W$.
\end{cor}
\begin{pf}
In $W$: We have already observed that $A$ is an antichain
of $\poset(\H)^W$, and hence it remains to establish its maximality. 
By assumption, $x_p\in V$ for all $p\in\poset(\H)^W$. 
Fix $p\in\poset(\H)^W$, and take any $\bar p\in\poset(\H)^V$ with
$x_{\bar p}=x_p$. 

\smallskip
\noindent\emph{Case} 1. $x_q\sqsubseteq x_p$ for some $q\in A$ compatible
with $\bar p$. 
\smallskip

\noindent Then from the definition of $\poset(\H)$,
$x_p\in\shade_\H(\J)$ for all $\J\in\X_q$, and conversely this proves
that $p$ is also compatible with $q$, as required.

\smallskip
\noindent\emph{Case} 2. There is no $q\in A$ compatible with $\bar p$ such that $x_q\sqsubseteq x_p$.
\smallskip

\noindent  Observe that, in $V$, there is no countable decomposition of
$\pi[A]_{[x_p]}$ into pieces orthogonal to $\Psi(\H_{[x_p]})$:
For if there was, then by lemma~\ref{l-36}, and also
proposition~\ref{p-41} using the fact that $\H$ is a $P$-ideal, 
there would be a countable
family $\X$ of cofinal subsets of $\H_{[x_p]}$ with 
$\pi[A]_{[x_p]}\cap\bigcap_{\J\in\X}\shade_{\H_{[x_p]}}(\J)\subseteq\{\emptyset\}$; 
then  $p'=\bigl(x_p,\X_{\bar p}\cup\bigl\{\{x_p\cup y:y\in\J\}:\J\in\X\bigr\}\bigr)$
is a condition of $\poset(\H)$ extending  $\bar p$. 
Since we are in Case~2, by maximality in $V$ there would
exist $q\in A$ compatible with $p'$ with $x_p\sqsubset x_q$.
However, this would entail that 
$x_q\setminus x_p\in\pi[A]_{[x_p]}\cap\bigcap_{\J\in\X}\shade_{\H_{[x_p]}}(\J)$, 
contradicting the fact that $x_q\setminus x_p\ne\emptyset$.

Therefore, by the hypothesis of the corollary,
in $W$ there can be no countable decomposition of
$\pi[A]_{[x_p]}$ into pieces orthogonal to $\Psi(\H_{[x_p]})$. 
Therefore, lemma~\ref{l-36} implies that
\begin{equation}
\pi[A]_{[x_p]}\cap\bigcap_{\J\in\X_p}\shade_{\H_{[x_p]}}(\J_{[x_p]})\ne\emptyset,
\label{eq:63}
\end{equation}
say $y$ is in the intersection. Then $y\in\pi[A]_{[x_p]}$ means that
there exists $q\in A$ such that $x_p\sqsubseteq x_q$ and
$y=x_q\setminus x_p$. 
And equation~\eqref{eq:63} implies that
$x_q\in\shade_\H(\J)$ for all $\J\in\X_p$, 
and thus $q$ is compatible with $p$ as required, because
$x_p\sqsubseteq x_q$. 
\end{pf}

Next we establish the analogous result for $\spo(\H)$. 

\begin{cor}
\label{l-37}
Let $A$ be a maximal antichain of $\spo(\H)$, 
and suppose $W\supseteq V$ is an outer model with $(\powcif\theta)^V=(\powcif\theta)^W$. 
If for every $x\in\H$, either
\begin{enumerate}[label=\tu{(\alph*)}, ref=\alph*, widest=b]
\item\label{item:66} there exists in $V$,
  a countable decomposition of $\pi[A]_{(x)}$ into
  pieces orthogonal to $\Psi(\H_{(x)})$, or
\item\label{item:67} there does not exist in $W$, 
a countable decomposition of $\pi[A]_{(x)}$ into
  pieces orthogonal to $\Psi(\H_{(x)})$,
\end{enumerate}
then $A$ is a maximal antichain of $\spo(\H)^W$.
\end{cor}
\begin{pf}
In $W$: We have already observed that $A$ is an antichain
of $\spo(\H)^W$, and hence it remains to establish its maximality. 
By assumption, $x_p\in V$ for all $p\in\spo(\H)^W$. 
Fix $p\in\spo(\H)^W$, and take any $\bar p\in\spo(\H)^V$ with
$x_{\bar p}=x_p$. 

\smallskip
\noindent\emph{Case} 1. $x_q\subseteq x_p$ for some $q\in A$ compatible
with $\bar p$. 
\smallskip

\noindent Then applying proposition~\ref{p-50} to $\bar p$ and $q$,
$x_p\in\shade_\H(\J)$ for all $\J\in\X_q$, and conversely this proves
that $p$ is also compatible with $q$, as required.

\smallskip
\noindent\emph{Case} 2. There is no $q\in A$ compatible with $\bar p$ 
such that $x_q\subseteq x_p$.
\smallskip

\noindent  Observe that, in $V$, there is no countable decomposition of
$\pi[A]_{(x_p)}$ into pieces orthogonal to $\Psi(\H_{(x_p)})$:
For if there was, then by lemma~\ref{l-36} there would be a countable
family $\X$ of cofinal subsets of $\H_{(x_p)}$ with 
$\pi[A]_{(x_p)}\cap\bigcap_{\J\in\X}\shade_{\H_{(x_p)}}(\J)\subseteq\{\emptyset\}$; 
then  $p'=\bigl(x_p,\X_{\bar p}\cup\bigl\{\{x_p\cup y:y\in\J\}:\J\in\X\bigr\}\bigr)$ 
is a condition of $\spo(\H)$ extending  $\bar p$. 
Since we are in Case~2, by maximality in $V$ there would
exist $q\in A$ compatible with $p'$ with $x_q\nsubseteq x_q$.
However, this would entail that 
$x_q\setminus x_p\in\pi[A]_{(x_p)}\cap\bigcap_{\J\in\X}\shade_{\H_{(x_p)}}(\J)$, 
contradicting the fact that $x_q\setminus x_p\ne\emptyset$.

Therefore, by the hypothesis of the corollary, 
in $W$ there can be no countable decomposition of
$\pi[A]_{(x_p)}$ into pieces orthogonal to $\Psi(\H_{(x_p)})$. 
Therefore, lemma~\ref{l-36} implies that
\begin{equation}
\pi[A]_{(x_p)}\cap\bigcap_{\J\in\X_p}\shade_{\H_{(x_p)}}(\J_{(x_p)})\ne\emptyset,
\label{eq:107}
\end{equation}
say $y$ is in the intersection. Then $y\in\pi[A]_{(x_p)}$ means that
there exists $q\in A$ such that $x_p\subseteq x_q$ and
$y=x_q\setminus x_p$. 
And equation~\eqref{eq:107} implies that
$x_q\in\shade_\H(\J)$ for all $\J\in\X_p$, 
and thus it follows from proposition~\ref{p-50} that 
$q$ is compatible with $p$, as required.
\end{pf}

\begin{cor}[$V\models\pstar$]
  \label{l-34}
Suppose that $W\supseteq V$ is an outer model with\linebreak
$(\powcif\theta)^V=(\powcif\theta)^W$.
Then $\poset(\H)^V\embeds[V]\poset(\H)^W$ 
\tu(and $\spo(\H)^V\embeds[V]\spo(\H)^W$\tu) via the inclusion
map. Furthermore, in the case $\theta=\oone$, we can weaken the
assumption to~$\ch+\pstar_\oone$. 
\end{cor}
\begin{pf}
Lemma~\ref{l-35} and corollary~\ref{c-16}. 
In the case $\theta=\oone$, under $\ch$, $\pi[A]$ has cardinality at most
$\aleph_1$ for every antichain $A$. 
\end{pf}

In particular, if $\H$ has no stationary orthogonal set, 
so that $\Q(\H)$ is completely proper and thus adds no new countable subsets of $\theta$,
then assuming $\pstar$, $\poset(\H)\forces\poset(\H)^V\embeds[V]\poset(\H)$.
Note that this is weaker than the statement
$\poset(\H)\forces\allowbreak
\poset(\H)^V\embeds\poset(\H)$,\fxwarning{typesetting}
that would in particular imply $\poset(\H)\times\poset(\H)$ is
proper by lemma~\ref{lem:proper:b}, because $\Q(\H)$ (in particular)
forces that $\H$ has no stationary orthogonal set.
This latter property, 
that is the square being proper, is the essence behind Shelah's $\nnr$ theory 
in~\cite[Ch.~XVIII,~\Section2]{MR1623206}. We have thus been led to
the following notion.

\begin{defn}
\label{d-24}
Let $\C$ be some class. 
Suppose that $P$ is a forcing notion definable from some parameter $a$.
We say that  $P(a)$ is \emph{$\C$-frozen} 
over a transitive model $M\ni P(a)$ of (enough of) $\zfc$, if 
for every outer model $N\supseteq M$ satisfying
\begin{enumerate}[label=(\roman*), ref=\roman*, widest=ii]
\item\label{item:126} $N\models\ulc \Gen(M,P(a)^M)\ne\emptyset\urc$,
\item\label{item:127} $\C^M=\C^N$, 
\end{enumerate}
we have $N\models\ulc P(a)^M\imbeds[M] P(a)^N\urc$.
\end{defn}

\noindent In other words, in every outer model $N$ extending some generic extension of
$M$ by~$P$ and preserving $\C$, $P(a)^M$ generically embeds via the
identity into $P(a)^N$ over~$M$. 

\begin{example}
\label{x-8}
$\pstar$ implies that $\poset(\H)$ is $\powcif\theta$-frozen over $V$,
for any $\sigma$-directed $\H\subseteq\powcif\theta$ 
with no stationary orthogonal set. 
In fact, corollary~\ref{l-34} says that this is true for every
outer model preserving $\powcif\theta$, and not just those also
satisfying~\eqref{item:126}. 
\end{example}

To obtain a model of $\pstarc_\oone$ with $\ch$ it is necessary (at least in our approach)
that for every forcing notion 
appearing in our iteration, of the form $\poset(\H)$, the property
that $\ulc A$ is a maximal antichain of $\poset(\H)\urc$ is upwards
absolute for forcing extensions derived in various ways from the
iteration (of course this will be made precise). 
Note that this entails preserving maximal antichains at every stage, 
because once the maximality of an antichain is lost it can never be restored. 

\fxnote{comment}

So far we have demonstrated that it is possible to freeze antichains
of $\poset(\H)$ by forcing uncountable sets locally in the appropriate
$\Psi(\H_{[x]})$. In fact, one can prove that $\Q(\H)$ itself forces
an uncountable set locally in each of the required~$\Psi(\H_{[x]})$, and
similarly for $\spo(\H)$. Thus we can strengthen example~\ref{x-8} by
eliminating $\pstar$, as follows,
although it should be noted that $\pstar$ cannot be eliminated from
corollary~\ref{l-34}. 

\begin{cor}
\label{c-12}
If $\H$ is a $P$-ideal on $\theta$ with stationary orthogonal set then
$\poset(\H)$ is $\powcif\theta$-frozen over $V$. Similarly, $\spo(\H)$
is $\powcif\theta$-frozen over $V$ for all $P$-ideals $\H$ on\/~$\theta$. 
\end{cor}

However, this approach cannot even handle two-stage iterations. 
By this we mean that it may not be possible to freeze all antichains
of say $\spo(\H)\gdot\spo(\dot \I)$. This is in spite of corollary~\ref{c-12}: 
For suppose $A\subseteq\spo(\H)\gdot\spo(\dot\I)$ is a maximal antichain.   
Let  $G\gdot H\in\Gen(V,\spo(\H)\gdot\spo(\dot\I))$. 
Then applying  corollary~\ref{c-12} in $V[G]$, $A\div G$
(cf.~proposition~\ref{p-39}) is frozen,
which means that $A\div G$ is a maximal antichain of $\spo(\dot\I[G])$
in every outer model of $V[G]$ preserving $\powcif\theta$. 
This does not however mean that the maximality of $A$ is preserved
because outer models of $V$ need not contain~$G$. What is needed, is an
$\spo(\H)$-name for an uncountable set locally in
$\Psi(\dot\I_{(x)})$, and we believe that this is generally impossible to
obtain.

What has been achieved in this section with corollaries~\ref{c-16}
and~\ref{l-37}, is that the problem of preserving the maximality
of antichains of $\poset(\H)$ and $\spo(\H)$ 
has been reduced to preserving the property that certain $P$-ideals have
no countable decompositions of their underlying set into pieces
orthogonal to them. 

\subsubsection{Products of $P$-ideals}
\label{sec:product-p-ideals}

\begin{notn}
\label{o-9}
Let $(\H_i:i\in I)$ be an indexed family where each $\H_i$ is a family
$\H_i\subseteq\power(S_i)$ of subsets of some fixed set $S_i$. 
Define
\begin{equation}
  \label{eq:82}
  \bigotimes_{i\in I}\H_i=\left\{\coprod_{i\in I}x(i):\vec x\in\prod_{i\in I}\H_i\right\},
\end{equation}
where $\coprod$ denotes the coproduct, i.e.~disjoint union.
Notice that $\bigotimes_{i\in I}\H_i\subseteq\power\bigl(\coprod_{i\in I}S_i\bigr)$.
\end{notn}

\begin{defn}
\label{d-14}
Let $(X_i:i\in I)$ be an indexed family of sets,
where each set has a zero $0_i\in X_i$. For $\vec x\in\prod_{i\in I}X_i$,
write $\supp(\vec x)=\{i\in I:x(i)\ne0_i\}$. 
The \emph{$\varSigma$-product} of $(X_i:i\in I)$ has the usual meaning:
\begin{equation}
  \label{eq:60}
  \sum\left(\prod_{i\in I}X_i\right)=\left\{\vec x\in\prod_{i\in I}X_i:
    \supp(\vec x)\text{ is countable}\right\}.
\end{equation}
\end{defn}

\begin{notn}
Suppose that $\H_i\subseteq\power(S_i)$, 
and moreover that $\emptyset\in\H_i$, for all $i\in I$. 
Taking $0_i=\emptyset$ for all $i\in I$, 
we extend the $\varSigma$-product notation as follows
\begin{equation}
  \label{eq:59}
  \sum\left(\bigotimes_{i\in I}\H_i\right)
  =\left\{\coprod_{i\in I}x(i):\vec x\in\sum\left(\prod_{i\in I}\H_i\right)\right\}. 
\end{equation}
\end{notn}

\begin{prop}
\label{p-30}
Let\/ $\H_i\subseteq\powcif{S_i}$ for every\/ $i\in I$. Then\/
$\sum\bigl(\bigotimes_{i\in I}\H_i\bigr)
\subseteq\bigl[\coprod_{i\in I}S_i\bigr]{}^{\aleph_0}$.
\end{prop}

\begin{prop}
\label{l-41}
Let $\H_0,\dots,\H_{n-1}$ be a finite sequence of $\sigma$-directed
subfamilies of $(\powcif{S_i},\subseteqfnt)$ for each $i=0,\dots,n-1$. 
Then $\bigotimes_{i=0}^{n-1}\H_i$ is a $\sigma$-directed subfamily of $(\powcif{\coprod_{i=0}^{n-1}S_i},\subseteqfnt)$.
\end{prop}

We can do better with $P$-ideals:

\begin{lem}
\label{l-1}
Let $(\I_i:i\in I)$ be an indexed family of $P$-ideals for some arbitrary~$I$. 
Then $\sum\bigl(\bigotimes_{i\in I}\I_i\bigr)$ is a $P$-ideal.
\end{lem}
\begin{pf}
Let $(y_n:n<\omega)$ be an enumeration of members of 
$\sum\bigl(\bigotimes_{i\in I}\I_i\bigr)$, 
say each $y_n=\coprod_{i\in I}x_n(i)$ 
for some $\vec x_n\in\sum\bigl(\prod_{i\in I}\I_i\bigr)$. 
Then $J=\bigcup_{n<\omega}\supp(\vec x_n)$ is countable, say $J=\{i_k:k<\omega\}$.
For each $i\in J$, since $\I_i$ is $\sigma$-directed,
there exists $z_i\in\I_i$ such that $x_n(i)\subseteqfnt z_i$ for all
$n<\omega$. Hence, as $\I_i$ is a $P$-ideal, 
$\bigcup_{n\in A}x_n(i)\cup z_i\in\I_i$ for every finite
$A\subseteq\omega$. Therefore, $\vec w\in\sum\bigl(\prod_{i\in
  I}\I_i\bigr)$, where $\supp(\vec w)\subseteq J$ is given by
\begin{equation}
  \label{eq:61}
  w(i_k)=\bigcup_{n=0}^{k-1} x_n(i_k)\cup z_{i_k}
\end{equation}
for each $k<\omega$. And clearly $y_n\subseteqfnt \coprod_{i\in
  I}w(i)$ for all $n<\omega$. 
\end{pf}

\begin{lem}
\label{l-60}
Suppose that\/ $\H_0,\dots,\H_{n-1}$ are $\sigma$-directed
with each $\H_i\subseteq\powcif{S_i}$, and each $S_i$ uncountable.
Let $J$ be the set of all $i=0,\dots,n-1$ for which 
$S_i$ has no countable decomposition into pieces orthogonal to $\H_i$.
Then\/ $\spo\bigl(\bigotimes_{i=0}^{n-1}\H_i\bigr)$ forces that there
exists $X\subseteq\coprod_{i=0}^{n-1}S_i$ locally in 
$\bigotimes_{i=0}^{n-1}\H_i$ such that 
$X\cap S_i$ is uncountable for all $i\in J$.
Similarly for $\varSigma$-products of $P$-ideals.
\end{lem}
\begin{pf}
Essentially the same as for lemma~\ref{l-7}.
\end{pf}

\begin{prop}
\label{p-43}
Let $\H_0,\dots,\H_{n-1}$ be a finite sequence with each
$\H_i\subseteq\power(S_i)$. Then
$(\H_0\otimes\cdots\otimes\H_{n-1},\subseteqfnt)$ is order isomorphic
to $\H_0\times\cdots\times\H_{n-1}$ with the product order obtained
from $(\H_i,\subseteqfnt)$. 
\end{prop}

Recall the notion from~\cite{MR0002515}, 
where a map $f:D\to E$ between two directed
posets is called \emph{convergent} if every $e\in E$ has a $d\in D$
such that $f(a)\ge e$ for all $a\ge d$. Notice that $f$ is convergent
iff it maps cofinal subsets of $D$ to cofinal subsets of~$E$.
We say that $D$ is \emph{cofinally finer} than~$E$, written $E\cofle D$, 
if there exists a convergent map from $D$ into $E$. It was established by
Tukey (in~\cite{MR0002515}) that $E\cofle D$ is equivalent to the existence of a map
$g:E\to D$ that maps unbounded subsets of $E$ to unbounded subsets of $D$. 
Then $\cofle$ is a quasi ordering of the class directed posets, which
we refer to as the \emph{Tukey order}. 
For two directed quasi orders  $A$ and $B$, we use the same definition
of convergent maps. Then the existence of a convergent map from $A$
into~$B$ is equivalent to the existence of a convergent map from the poset
$A\div{\sim_{\mathrm{asym}}}$ into  the poset
$B\div{\sim_{\mathrm{asym}}}$, i.e. the antisymmetric quotient. 
Thus the Tukey ordering $\cofle$ also makes sense between directed
quasi orders. The notation $D\cofeq E$ indicates that $D$ is
\emph{cofinally equivalent} to $E$, i.e.~$D\cofle E$ and $E\cofle D$. 
Then $\cong$ is an equivalent relation, and the equivalence classes
are called \emph{cofinal types}. 

A basic result on this is as follows.

\begin{lem}[Tukey]
\label{l-64}
For any finite sequence $D_0,\dots,D_{n-1}$ of directed sets,
$D_0\times\cdots\times D_{n-1}$ is their least upper bound in the
Tukey order. 
\end{lem}

\begin{example}
\label{x-3}
$1$, $\omega$, $\oone$, $\omega\times\oone$ and $\powfin\oone$ are
five distinct cofinal types, where the first four orders are given
by the $\in$ relation and $\powfin\oone$ is ordered by $\subseteq$. 
It is proved in~\cite{MR792822} that: $\pfa$ implies that these five
are the only cofinal types of cardinality at most $\aleph_1$, 
while on the other hand, $\ch$ implies that there are~$2^{\aleph_1}$ 
many cofinal types of cardinality $\aleph_1$. 
\end{example}

\begin{example}[$\ch$]
\label{x-7}
\emph{If $\H\subseteq\powcif\oone$ and $(\H,\subseteqfnt)$ is
  $\sigma$-directed, then $(\H,\subseteqfnt)$ is of cofinal type 
  either $1$ or $\oone$. It has cofinal type $1$ iff
  $(\H,\subseteqfnt)$ has a maximal element.}
\end{example}

\begin{prop}
\label{p-45}
If $D$ is a directed set and $\kappa\cofle D$ for some infinite
cardinal $\kappa$ \tu(ordered by $\in$\tu),
then no bounded
subset of\/ $\kappa$ can be mapped onto a cofinal subset of~$D$. 
\end{prop}
\begin{pf}
If $f:\kappa\to D$ maps a bounded subset of $\kappa$ onto a cofinal
subset of $D$, then for any convergent $g:D\to\kappa$, $g\circ f$ maps
a bounded subset of $\kappa$ onto a cofinal subset of~$\kappa$, which is
impossible if $\kappa$ is an infinite cardinal.
\end{pf}

\begin{lem}
\label{p-44}
For any finite sequence $\H_0,\dots,\H_{n-1}$ where each $\H_i$ is a
directed subfamily of $(\powcif{S_i},\subseteqfnt)$, 
$\bigotimes_{i=0}^{n-1}\H_i$ is the $\cofle$-least upper bound of
the sequence, under the almost inclusion order.
\end{lem}
\begin{pf}
By proposition~\ref{p-43} and lemma~\ref{l-64}. 
\end{pf}

\begin{cor}
\label{p-42}
If $\H$ and $\I$ are $\subseteqfnt$-directed subfamilies of $\powcif
S$ and $\powcif T$, 
respectively, and $\I\cofle\H$, then $\H\cofeq\H\otimes\I$.
\end{cor}
\begin{pf}
Since $\H\otimes\I$ is the least upper bound of $\H$ and $\I$ by
lemma~\ref{p-44}. 
\end{pf}

We need something more specific.

\begin{lem}
\label{l-63}
Suppose that $\H$ and $\I$ are directed subfamilies of $(\powcif
S,\subseteqfnt)$ and $(\powcif T,\subseteqfnt)$, 
respectively, and both $\H$ and $\I$ have cofinal type $\kappa$
for some infinite cardinal $\kappa$.
Then every cofinal $\K\subseteq\H\otimes\I$ has a cofinal subset
$\L\subseteq\K$ such that for every cofinal subset
$\J\subseteq\pi[\L]=\{x\in\H:x\amalg y\in\L\tu{ for some }y\in\I\}$,
$(\J\otimes\I)\cap\L$ is cofinal in $\H\otimes\I$.
\end{lem}
\begin{pf}
Since $\H\otimes\I\cong\kappa$ by lemma~\ref{p-44},
there is a convergent map $g:\kappa\to\H\otimes\I$. 
For each $\alpha<\kappa$, since $\K$ is cofinal we can find
$x_\alpha\amalg y_\alpha\in\K$ such that
\begin{equation}
  \label{eq:93}
  g(\alpha)\subseteqfnt x_\alpha\amalg y_\alpha.
\end{equation}
We claim that $\L=\{x_\alpha\amalg y_\alpha:\alpha<\kappa\}$ satisfies
the conclusion: $\L$ is cofinal because $g$ is convergent. 
Suppose $\J\subseteq\pi[\L]$ is cofinal. Then $\J$ is cofinal in $\H$
as $\pi[\L]$ is by proposition~\ref{p-43}.
Thus $\J=\{x_\alpha:\alpha\in A\}$ for some cofinal $A\subseteq\kappa$,
because $g[B]$ is noncofinal for all bounded $B\subseteq\kappa$ by
proposition~\ref{p-45} as $\kappa\altle\H$. And
$(\J\otimes\I)\cap\L\supseteq\{x_\alpha\amalg y_\alpha:\alpha\in A\}$,
which is cofinal by~\eqref{eq:93} since $g$ is convergent.  
\end{pf}

\begin{lem}
\label{l-62}
Let $\H$ and $\I$ be $\sigma$-directed subfamilies of $(\powcif
S,\subseteqfnt)$ and\linebreak $(\powcif T,\allowbreak\subseteqfnt\nobreak)$, respectively.
If both $\H$ and $\I$ have cofinal type $\kappa$ for some infinite
cardinal~$\kappa$, then $\spo(\H)$ generically embeds into $\spo(\H\otimes\I)$. 
\end{lem}
\begin{pf}
Define $e:\spo(\H)\to\spo(\H\otimes\ideal)$ by $e(p)=(x_p,\Y_p)$ where 
\begin{equation}
  \label{eq:72}
  \Y_p=\{\J\otimes\I:\J\in\X_p\}.
\end{equation}
Given a maximal antichain $A\subseteq\spo(\H)$, we need to show that
$e[A]$ is a maximal antichain of $\spo(\H\otimes\ideal)$.
Take $q\in\spo(\H\otimes\ideal)$. Write $x_q=y\amalg z$ ($y\in\H$, $z\in\ideal$). 
For each $\K\in\X_q$, apply lemma~\ref{l-63} to the
cofinal set $\{w\in\K:x_q\subseteq w\}$ to obtain a cofinal subset
$\L_\K$ as in the conclusion of that lemma. 
Then clearly
\begin{equation}
q'=(y,\{\pi[\L_\K]:\K\in\X_q\})\label{eq:97}
\end{equation}
is a condition of $\spo(\H)$. 
Hence there must be $p\in A$ compatible with $q'$.
For all $\J\in\X_p$, $x_p\cup y\in\shade_\H(\J)$, 
and therefore $x_p\cup x_q=x_p\cup(y\amalg z)\in\shade_{\H\otimes\I}(\J\otimes\ideal)$; 
and for all $\K\in\X_q$, $\J_\K=\{x\in\pi[\L_\K]:x_p\cup y\subseteq x\}$ is cofinal, 
and therefore $(\J_\K\otimes\I)\cap\L_\K$ is cofinal, which implies
that $x_p\cup x_q\in\shade_{\H\otimes\I}(\K)$,
because $x_p\cup(y\amalg z)\subseteq w$ for all $w\in\L_\K$ with $\pi(w)\in\J_\K$.
This proves that $e(p)$ is compatible with $q$, as required.
\end{pf}

\begin{remark}
\label{r-11}
We do not believe that there is any analogue of lemma~\ref{l-62} for
$\poset$. This embedability of $\spo(\H)$ is the primary reason we are
interested in the forcing notion $\spo(\H)$ when we are only trying to
force clubs with the forcing notion $\poset(\I)$. For example, it figures in the analysis
of properties of the forcing notion $\poset(\I)$ in corollary~\ref{c-15}. 
There is also a secondary use of the forcing notion $\spo(\H)$ in
section~\ref{sec:proof-cons-ppst} where it is used to force stationary
sets. 
\end{remark}

\begin{cor}[$\ch$]
\label{c-10}
Let $\H$ and $\ideal$ be $\sigma$-directed subfamilies of $(\powcif\oone,\subseteqfnt)$. 
If $\I$ has  no countable decomposition of $\oone$ into orthogonal
pieces, then $\spo(\H)$ forces that $\I$ has no countable
decomposition of $\oone$ into orthogonal pieces.
\end{cor}
\begin{pf}
Let $\H$ and $\I$ be as in the hypothesis. 
By lemma~\ref{l-60}, $\spo(\H\otimes\I)$ forces that there is an
uncountable $X\subseteq\oone$ locally in $\I$ (meaning uncountable in
the forcing extension, i.e.~$\spo(\H\otimes\I)$ does not collapse $\aleph_1$). 
In particular, $\spo(\H\otimes\I)$ forces that there is no countable decomposition of
$\oone$ into pieces orthogonal to $\I$. 
By $\ch$ and example~\ref{x-7}, we know that both
$\H$ and $\I$ are of cofinal type either $1$ or $\oone$,
and by the hypothesis we know further that $\I\cong\oone$.
In the case $\H\cong1$, $\spo(\H)$ is the trivial forcing notion and
thus the conclusion of the corollary is trivial. Assume then that
 $\H\cong\I\cong\oone$ in the Tukey order. 
Then by lemma~\ref{l-62}, $\spo(\H)\embeds\spo(\H\otimes\I)$,
and thus $\spo(\H)$ cannot introduce a countable decomposition of
$\oone$ into pieces orthogonal to~$\I$.
\end{pf}

\begin{cor}[$\ch$]
\label{c-6}
Let $\H$ and $\ideal$ be $\sigma$-directed subfamilies of
$(\powcif\oone,\subseteqfnt)$, with $\I$ moreover a $P$-ideal.
Then $\spo(\H)\forces\spo(\I)^V\imbeds[V]\spo(\I)$\tu; hence\tu,
$\spo(\I)\embeds\spo(\H)\gdot\spo(\I)$.
\end{cor}
\begin{pf}
Let $G\in\Gen(V,\spo(\H))$, and set $W=V[G]$. 
In $V$, let $A$ be a maximal antichain of $\spo(\I)$.
Then for all $x\in\I$, if $\pi[A]_{(x)}$ has no countable
decomposition, in $V$, into pieces orthogonal to $\Psi(\I_{(x)})$,
then by corollary~\ref{c-10} with $\I:=\Psi(\I_{(x)})$, which is
$\sigma$-directed by proposition~\ref{p-41} and lemma~\ref{l-35}, 
$\pi[A]_{(x)}$ has no countable decomposition, in $W$, into pieces
orthogonal to $\Psi(\I_{(x)})$. Therefore, $A$ is a maximal antichain
of $\spo(\I)^W$ by corollary~\ref{l-37}. 

$\spo(\I)\embeds \spo(\H)\gdot\spo(\I)$ is immediate from proposition~\ref{p-52}.
\end{pf}

\begin{cor}[$\ch$]
\label{c-14}
Let $\H$ and $\ideal$ be $\sigma$-directed subfamilies of
$(\powcif\oone,\subseteqfnt)$, with $\I$ moreover a $P$-ideal.
Then $\spo(\H)\forces\poset(\I)^V\imbeds[V]\poset(\I)$\tu; 
hence, $\poset(\I)\embeds\spo(\H)\gdot \poset(\I)$. 
\end{cor}
\begin{pf}
This is the same as the proof of corollary~\ref{c-6} but using
corollary~\ref{c-16}.
\end{pf}

\begin{cor}[$\ch$]
\label{c-15}
Let $\H$ and $\ideal$ be $\sigma$-directed subfamilies of
$(\powcif\oone,\subseteqfnt)$, with $\I$ moreover a $P$-ideal.
Suppose $\H$ has no countable decomposition of~$\oone$ into orthogonal
pieces. 
If $\spo(\H)$ forces that $\ideal$ has no stationary orthogonal subset
of~$\oone$, then $\poset(\I)$ forces that $\H$ has no countable
decomposition of $\oone$ into orthogonal pieces. 
\end{cor}
\begin{pf}
By corollary~\ref{c-14}, $\poset(\I)\embeds\spo(\H)\gdot\poset(\I)$.
Now, if $\spo(\H)$ forces that $\ideal$ has no stationary orthogonal
set, then $\spo(\H)\gdot\poset(\I)$ is proper by lemmas~\ref{l-24}
and~\ref{lem:proper:b}, and hence does not collapse $\aleph_1$. 
Therefore, by our assumption on $\H$, $\spo(\H)\gdot\poset(\I)$ forces
that there exists an uncountable set locally in $\downcl\H$. 
Thus $\poset(\I)$ cannot force a countable decomposition of~$\oone$
into pieces orthogonal to $\H$. 
\end{pf}

\begin{cor}[$\ch$]
\label{c-17}
Let $\I$ be a $P$-ideal on $\oone$.
Suppose that $\spo(\H)$ forces that there is no stationary subset
of\/ $\oone$ orthogonal to $\I$, for every $\sigma$-directed subfamily $\H$ of $(\powcif\oone,\subseteqfnt)$
having no countable decomposition of\/ $\oone$ into orthogonal pieces. 
Then for every $P$-ideal $\J$ on $\oone$, 
$\poset(\I)\forces\spo(\J)^V\imbeds[V]\spo(\J)$\tu; hence,
$\spo(\J)\embeds\poset(\I)\gdot\spo(\J)$. 
\end{cor}
\begin{pf}
Let $G\in\Gen(V,\poset(\I))$, and set $W=V[G]$. 
In $V$, let $A$ be a maximal antichain of $\spo(\J)$.
Then for all $x\in\J$, if $\pi[A]_{(x)}$ has no countable
decomposition, in $V$, into pieces orthogonal to $\Psi(\J_{(x)})$,
then by corollary~\ref{c-15} with $\H:=\Psi(\J_{(x)})$, which is
$\sigma$-directed by proposition~\ref{p-41} and lemma~\ref{l-35}, 
$\pi[A]_{(x)}$ has no countable decomposition, in $W$, into pieces
orthogonal to $\Psi(\J_{(x)})$. Therefore, $A$ is a maximal antichain
of $\spo(\J)^W$ by corollary~\ref{l-37}. 
\end{pf}

Similarly:

\begin{cor}[$\ch$]
\label{c-18}
Let $\I$ be a $P$-ideal on $\oone$.
Suppose that $\spo(\H)$ forces that there is no countable
decomposition of $\oone$ into pieces orthogonal to $\I$, for every
$\sigma$-directed subfamily $\H$ of $(\powcif\oone,\subseteqfnt)$
having no countable decomposition of $\oone$ into orthogonal pieces. 
Then for every $P$-ideal $\J$ on $\oone$, 
$\poset(\I)\forces\poset(\J)^V\imbeds[V]\poset(\J)$\tu; hence,
$\poset(\J)\embeds\poset(\I)\gdot\poset(\J)$. 
\end{cor}

In corollaries~\ref{c-6},~\ref{c-14},~\ref{c-17} and~\ref{c-18}, all four
permutations of $\spo(\I)$ ($\poset(\I)$) in the extension by
$\spo(\H)$ [$\poset(\H)$] have been considered. The next step is to
consider iterations. It is immediate from two applications of
corollary~\ref{c-6}, that
$\spo(\H)\gdot\spo(\dot\I)\forces\spo(\J)^V\imbeds[V]\spo(\J)$. But
this begs the question of whether 
\begin{equation}
\spo(\J)\forces[\spo(\H)\gdot\spo(\dot\I)]^V\imbeds[V][\spo(\H)\gdot\spo(\dot\I)]\tu?
\label{eq:127}
\end{equation}
The key to answering this is to establish that, under certain
conditions, the forcing notions $\poset(\H)$ and $\spo(\H)$ commute
among themselves. 

\begin{notn}
\label{o-10}
For any forcing notion $P$ and any $P$-name $\dot A$, 
we denote
\begin{equation}
  \label{eq:132}
  \dot A[p]=\{x:q\forces x\in\dot A\text{ for some }q\ge p\},
\end{equation}
for each $p\in P$.
\end{notn}

\begin{prop}
\label{p-59}
If $P\forces\dot A\subseteq V$ 
then $p\forces\dot A\subseteq\dot A[p]$, for all $p\in P$. 
\end{prop}

\begin{remark}
\label{r-6}
Note that when we say ``$p$ decides $\dot A$'' this means the same
thing as ``$p\forces\dot A=\dot A[p]$''. 
\end{remark}

\begin{prop}
\label{p-63}
Let $\H\subseteq\powcif\theta$. Suppose $P$ is a forcing notion that
adds no new countable subsets of $\theta$.
If $(p,\dot q)$ is a condition of $P\gdot\poset(\H)$ and $p$ decides
$x_{\dot q}$, then $(x_{\dot q}[p],\{\dot\J[p]:\dot\J\in\X_{\dot
  q}\})\in\poset(\H)$. Similarly, for $\spo(\H)$.  
\end{prop}
\begin{pf}
By the assumption on $P$, 
 $P\forces\dot\J\subseteq V$ for all $\dot\J\in\X_{\dot q}$.
The result thus follows from proposition~\ref{p-59}. 
\end{pf}

\noindent Note that we are implicitly assuming an enumeration of
$\X_{\dot q}$ by $\aleph_0$ when referring to $\dot\J\in\X_{\dot q}$.  

\begin{lem}
\label{l-40}
Let $\H\subseteq\powcif\theta$. 
For every forcing notion $P$ that adds no new countable subsets of~$\theta$, 
if $(p,\dot q)\in P\gdot\poset(\H)$ and $D\subseteq\{r\in P:r\ge p\}$ is predense above $p$ 
then\fxwarning{Note in the introduction that BA's are upside down for us}
\begin{equation}
\overline{(P\gdot\poset(\H))\div{\qsep}}
\models(p,\dot q)\ge\ifm_{d\in D}
\bigl(d,(x_{\dot q},\{\dot\J[d]:\dot\J\in\X_{\dot q}\})\bigr), 
\label{eq:133}
\end{equation}
Similarly, for $\spo(\H)$. 
\end{lem}

\begin{remark}
\label{r-5}
Equation~\eqref{eq:133} is equivalent to: 
for every $(p',\dot q')\ge (p,\dot q)$ there exists $d\in D$ such that $(p',\dot q')$ is
compatible with $\bigl(d,(x_{\dot q},\{\dot\J[d]:\dot\J\in\X_{\dot q}\})\bigr)$. 
\end{remark}

\begin{pf}[of lemma~\tu{\ref{l-40}}]
We establish equation~\eqref{eq:133} using remark~\ref{r-5}.
Given $(p',\dot q')\ge(p,\dot q)$, since $D$ is predense above $p$, there
exists $d\in D$ compatible with $p'$, say with common extension $d'$. 
By the assumption on $P$, proposition~\ref{p-59} applies, 
and then by proposition~\ref{p-60}, 
\begin{equation}
d\forces\overline{\poset(\H)\div{\qsep}}
\models \dot q\ge(x_{\dot q},\{\dot\J[d]:\dot\J\in\nobreak\X_{\dot q}\}),\label{eq:134}
\end{equation}
which implies that there is an $\dot s$ such that $d\forces\dot s$ is a
common extension of $\dot q'$ and 
$(x_{\dot q},\{\dot\J[d]:\dot\J\in\nobreak\X_{\dot q}\})$. 
Therefore, $(d',\dot s)\ge (p',\dot q')$ and 
$(d',\dot s)\ge\bigl(d,(x_{\dot q},\{\dot\J[d]:\dot\J\in\nobreak\X_{\dot q}\})\bigr)$, 
concluding the proof.
\end{pf}

\begin{lem}
\label{l-72}
Let $\H\subseteq\powcif\theta$. If $P$ is a forcing notion that adds
no new countable subsets of\/ $\theta$, and
$P\forces\poset(\H)^V\imbeds[V]\poset(\H)$,
then $P\gdot\poset(\H)\div\poset(\H)\cong P$. 
\end{lem}
\begin{pf}
First we deal with the pathological case where some $q\in\poset(\H)$
forces that $C_{\dot G_{\poset(\H)}}$ is countable. 
Let $A$ be a antichain maximal with respect to every $q\in A$ having this property. 
Then let $B\subseteq\poset(\H)$ satisfy $A\cup B$ is a maximal
antichain. For all  $a\in A$, clearly $\poset(\H)_a$ is the trivial
forcing notion and moreover this is upwards absolute, and thus
$P\gdot\poset(\H)_a\div\poset(\H)_a\cong P\gdot 1\div 1\cong P$. Now,
by equation~\eqref{eq:138}, it suffices to prove that
$P\gdot\poset(\H)_b\div\poset(\H)_b\cong P$ for all $b\in
B$. Henceforth, we assume without loss of generality that
$\poset(\H)\forces C_{\dot G_{\poset(\H)}}$ is uncountable. 

By proposition~\ref{p-52}, the map $e:\poset(\H)\to P\gdot\poset(\H)$ given by
\begin{equation}
e(q)=(0_{\poset(\H)},q)\label{eq:137}
\end{equation}
defines a generic embedding of $\poset(\H)$ into $P\gdot\poset(\H)$.
Let $G\in\Gen(V,\poset(\H))$. 

In $V[G]$:  The representation of the quotient given by $e$ is
\begin{multline}
  \label{eq:136}
  (P\gdot \poset(\H))\div G=\{(p,\dot q)\in Q\gdot\dot R:\\
  (p,\dot q)\text{ is $Q\gdot\dot R$-compatible with every member of }e[G]\},
\end{multline}
with the order inherited from $P\gdot\poset(\H)$. 
Thus, as $e[G]=\{0_Q\}\times G$, 
$(p,\dot q)\in (P\gdot \poset(\H))\div G$ 
iff every $r\in G$ has a $p'\ge p$ forcing that $r$ is compatible with
$\dot q$. 

\begin{lemclaim}
\label{a-3}
For all $(p,\dot q)\in (P\gdot \poset(\H))\div G$, 
there exists $p'\ge p$ such that $p'$ decides $x_{\dot q}$ 
and $(p'',\dot q)\in (P\gdot\poset(\H))\div G$ for all $p''\ge p'$. 
\end{lemclaim}
\begin{pf}
Let $(p,\dot q)\in (P\gdot\poset(\H))\div G$ be given. 
Then letting $D$ be the set of all $d\ge\nobreak p$ deciding $x_{\dot q}$,
$D$ is dense above $p$ by our assumption on $P$. 
Note that $(x_{\dot q}[d],\{\dot\J[d]:\dot\J\in\X_{\dot q}\})\in\poset(\H)$ 
for all $d\in D$ by proposition~\ref{p-63}. 

Now assume towards a contradiction that claim~\ref{a-3} fails.
Let $E$ be the set of all $d\in D$ 
for which there is an $r_d\in G$ such that
\begin{equation}
  \label{eq:135}
  r_d\perp (x_{\dot q}[d],\{\dot\J[d]:\dot\J\in\X_{\dot q}). 
\end{equation}

\begin{lemsubclaim}
\label{a-5}
$E$ is dense above $p$. 
\end{lemsubclaim}
\begin{pf}
Take $p_0\ge p$. Pick $d\ge p_0$ in $D$. 
By our assumption that the claim fails, 
there exists $p_1\ge d$ 
such that $(p_1,\dot q)\notin(P\gdot\poset(\H))\div G$. 
Hence there is an $r\in G$ and $p_2\ge p_1$ forcing that $r$ is
incompatible with $\dot q$. Since $p_2\forces x_{\dot q}=x_{\dot q}[d]$,
by proposition~\ref{p-64}, either $x_{\dot q}[d]$ is not comparable
under end-extension with $x_r$, in which case it is clear that
$r_d:=r$ witnesses that $d\in E$, or else $p_2$ forces that there
exists $\J\in\X_{\dot q}\cup\X_{r}$ 
with $x_{\dot q}[d]\cup x_r\notin\shade_\H(\J)$. 
If it is the case that $x_{\dot q}[d]\notin\shade_\H(\J)$ for some $\J\in\X_r$, 
then $r_d:=r$ witnesses that $d\in E$. Otherwise, in the remaining case there exists
$p_3\ge p_2$ and $\dot\J\in\X_{\dot q}$ such that 
$p_3\forces x_r\notin\shade_\H(\dot\J)$. Hence, there exists $y\in\H$
and $p_4\ge p_3$ forcing that there is no $z\supseteqfnt y$ in
$\dot\J$ with $x_r\subseteq z$. Therefore, there is no
$z\in\dot\J[p_4]$ with $y\subseteqfnt z$ and $x_r\subseteq z$,
i.e.~$x_r\notin\shade_\H(\dot\J[p_4])$. Now $d:=p_4$ and $r_d:=r$
witness that $p_4\in E$. 
\end{pf}

\begin{lemsubclaim}
\label{subc:1}
There exists $d\in E$ such that 
$(x_{\dot q}[d],\{\dot\J[d]:\dot\J\in\nobreak\X_{\dot q}\})\in G$.
\end{lemsubclaim}
\begin{pf}
Suppose not. Then there exists $r\in G$ such that
$r\forces (x_{\dot q}[d],\{\dot\J[d]:\dot\J\in\nobreak\X_{\dot
  q}\})\notin\dot G_{\poset(\H)}$ for all $d\in E$. This means that 
\begin{equation}
  \label{eq:139}
  r\perp (x_{\dot q}[d],\{\dot\J[d]:\dot\J\in\nobreak\X_{\dot q}\})
  \espc\text{for all $d\in E$}.
\end{equation}
However,  $(p,\dot q)\in(P\gdot\poset(\H))\div G$ implies that there
exists $(p',\dot q')\ge (p,\dot q)$ in $ P\gdot\nobreak\poset(\H)$ 
such that $p'\forces\dot q'\ge r$. And by subclaim~\ref{a-5},
and lemma~\ref{l-40} with $D:=E$,
 there exists $d\in E$ such that $(p',\dot q')$ is compatible with
$\bigl(d,(x_{\dot q},\{\dot\J[d]:\dot\J\in\nobreak\X_{\dot q}\})\bigr)
=\bigl(d,(x_{\dot q}[d],\{\dot\J[d]:\dot\J\in\nobreak\X_{\dot q}\})\bigr)$. 
This clearly implies that $r$ is compatible with 
$(x_{\dot q}[d],\{\dot\J[d]:\dot\J\in\X_{\dot q}\})$, contradicting~\eqref{eq:139}.
\end{pf}

Let $d$ be as in subclaim~\ref{subc:1}.
Then as in equation~\eqref{eq:135}, there exists $r\in G$ such that
$r\perp (x_{\dot q}[d],\{\dot\J[d]:\dot\J\in\X_{\dot q}\})$. 
This obviously contradicts $(x_{\dot q},\{\dot\J[d]:\dot\J\in\nobreak\X_{\dot
  q})\in G$. 
\end{pf}

Claim~\ref{a-3} allows us to define $f:(P\gdot\poset(\H))\div G\to P$
so that for all $(p,\dot q)\in(P\gdot\poset(\H))\div G$,  
$f(p,\dot q)\in D$, $f(p,\dot q)\ge p$ and $(p',\dot q)\in(P\gdot\poset(\H))\div G$ 
for all $p'\ge f(p,\dot q)$. 
Thus $f$ satisfies clauses~\eqref{item:146} and~\eqref{item:144} of
lemma~\ref{l-55} with $O:=\poset(\H)$ and $\dot Q:=\poset(\H)^{V[\dot G_P]}$. 
Observe that for all $(p,\dot q)\in (P\gdot\poset(\H))\div G$,
\begin{equation}
  \label{eq:141}
  f(p,\dot q)\forces\dot q\text{ is compatible with }r\espc\text{for
    all $r\in G$}.
\end{equation}

It remains to verify clause~\eqref{item:145} that $(f(p,\dot q),\dot
q)$ and $(f(p',\dot q'),\dot q')$ are $P\gdot\dot Q$-compatible 
whenever $f(p,\dot q)$ and $f(p',\dot q')$ are compatible. 
Then lemma~\ref{l-55} will yield $(P\gdot\poset(\H))\div\poset(\H)\cong\poset(\H)$.

Suppose then that $f(p,\dot q)$ and $f(p',\dot q')$ are compatible,
say $p''$ is a common extension. 
By our  assumption that $C_G$ is uncountable, there exists $r\in G$
such that $x_r\nsubseteq x_{\dot q}[f(p,\dot q)]$ 
and $x_r\nsubseteq x_{\dot q'}[f(p',\dot q')]$. Therefore,
by~\eqref{eq:141} and proposition~\ref{p-64}, $x_{\dot q}[f(p,\dot
q)]$ and $x_{\dot q'}[f(p',\dot q')]$ are both initial segments of
$x_r$ and are thus comparable under end-extension. Again
by~\eqref{eq:141} and proposition~\ref{p-64}, 
$f(p,\dot q)\forces x_r\in\shade_\H(\dot\J)$ for all $\dot\J\in\X_{\dot q}$ 
and $f(p',\dot q')\forces x_r\in\shade_\H(\dot\J)$ for all $\dot\J\in\X_{\dot q'}$. 
Thus $p''\forces\ulc x_{\dot q}\cup x_{\dot q'}\subseteq
x_r\in\shade_\H(\J)$ for all $\J\in\X_{\dot q}\cup\X_{\dot q'}\urc$,
proving that $p''\forces \dot q$ and $\dot q'$ are compatible, by
proposition~\ref{p-64}. 
\end{pf}

\begin{lem}
\label{l-73}
Let $\H\subseteq\powcif\theta$. If $P$ is a forcing notion that adds
no new countable subsets of\/~$\theta$, and
$P\forces\spo(\H)^V\imbeds[V]\spo(\H)$, then
$P\gdot\spo(\H)\div\spo(\H)\cong P$. 
\end{lem}
\begin{pf}
Essentially the same as lemma~\ref{l-72} but using
proposition~\ref{p-50}. 
\end{pf}

\begin{cor}[$\ch$]
\label{c-9}
Let $\H$ and $\I$ be $\sigma$-directed subfamilies of
$(\powcif\oone,\subseteqfnt)$, with $\I$ moreover a $P$-ideal.
Then
\begin{enumerate}[label=\tu{(\alph*)}, ref=\alph*, widest=b]
\item \label{item:147} $\spo(\H)\gdot\spo(\I)\div\spo(\I)\cong \spo(\H)$,
\item\label{item:148}  $\spo(\H)\gdot\poset(\I)\div\poset(\I)\cong\spo(\H)$.
\end{enumerate}
\end{cor}
\begin{pf}
For conclusion~\eqref{item:147}, we can apply lemma~\ref{l-73} with
$P:=\spo(\H)$ and $\H:=\I$, because $\spo(\H)\forces\spo(\I)^V\imbeds[V]\spo(\I)$
by corollary~\ref{c-6}. 

For~\eqref{item:148}, we can apply lemma~\ref{l-72} with
$P:=\spo(\H)$, by corollary~\ref{c-14}. 
\end{pf}

\begin{cor}[$\ch$]
\label{c-22}
Let $\H$ and $\I$ be $P$-ideals on $\oone$.
Suppose that $\spo(\J)$ forces that there is no stationary set
orthogonal to $\H$ for every $\sigma$-directed subfamily $\J$ of
$(\powcif\oone,\subseteqfnt)$ having no countable decomposition of
$\oone$ into orthogonal pieces. 
Then
\begin{enumerate}[label=\tu{(\alph*)}, ref=\alph*, widest=b]
\item\label{item:149}  $\poset(\H)\gdot\spo(\I)\div\spo(\I)\cong \poset(\H)$,
\item\label{item:150} $\poset(\H)\gdot\poset(\I)\div\poset(\I)\cong\poset(\H)$. 
\end{enumerate}
\end{cor}
\begin{pf}
By the hypothesis, corollary~\ref{c-17} applies so
that
$\poset(\H)\forces\spo(\I)^V\imbeds[V]\spo(\I)$, and
thus
conclusion~\eqref{item:149} holds by lemma~\ref{l-73}.

Similarly, corollary~\ref{c-18} applies so that
$\poset(\H)\forces\poset(\I)^V\imbeds[V]\poset(\I)$,
and thus
conclusion~\eqref{item:150} holds by lemma~\ref{l-72}.
\end{pf}

We have now achieved commutativity.

\begin{cor}[$\ch$]
\label{c-19}
Let $\H$ and $\I$ be $P$-ideals on $\oone$.
Then
\begin{equation}
  \label{eq:128}
  \spo(\H)\gdot\spo(\I)\cong\spo(\H)\times\spo(\I)\cong\spo(\I)\times\spo(\H)
  \cong\spo(\I)\gdot\spo(\H).
\end{equation}
\end{cor}
\begin{pf}
Corollary~\ref{c-9}\eqref{item:147} is equivalent to
$\spo(\H)\gdot\spo(\I)\cong\spo(\I)\times\spo(\H)$.
The remaining equivalences are by the commutativity of products and
another application of corollary~\ref{c-9}\eqref{item:147}. 
\end{pf}

\begin{cor}[$\ch$]
\label{c-23}
Let $\H$ and $\I$ be $P$-ideals on $\oone$.
Suppose that $\spo(\J)$ forces that there is no stationary set
orthogonal to $\H$ for every $\sigma$-directed subfamily $\J$ of
$(\powcif\oone,\subseteqfnt)$ having no countable decomposition of
$\oone$ into orthogonal pieces. 
Then
\begin{equation}
  \label{eq:43}
  \poset(\H)\gdot\spo(\I)\cong\poset(\H)\times\spo(\I)\cong\spo(\I)\times\poset(\H)
  \cong\spo(\I)\gdot\poset(\H).
\end{equation}
\end{cor}
\begin{pf}
We have $\poset(\H)\gdot\spo(\I)\cong\spo(\I)\times\poset(\H)$ 
by corollary~\ref{c-22}\eqref{item:149}, and
$\spo(\I)\gdot\poset(\H)\cong\poset(\H)\times\spo(\I)$ 
by corollary~\ref{c-9}\eqref{item:148}. 
\end{pf}

\begin{cor}[$\ch$]
\label{c-24}
Let $\H$ and $\I$ be $P$-ideals on $\oone$.
Suppose that $\spo(\J)$ forces that there is no stationary set
orthogonal to $\H$ and that there  is no stationary set orthogonal to
$\I$ for every $\sigma$-directed subfamily $\J$ of
$(\powcif\oone,\subseteqfnt)$ having no countable decomposition of
$\oone$ into orthogonal pieces. Then
\begin{equation}
  \label{eq:124}
  \poset(\H)\gdot\poset(\I)\cong\poset(\H)\times\poset(\I)\cong\poset(\I)\times\poset(\H)
  \cong\poset(\I)\gdot\poset(\H). 
\end{equation}
\end{cor}
\begin{pf}
By two applications of corollary~\ref{c-22}\eqref{item:150}. 
\end{pf}

\begin{remark}
\label{r-12}
This is already very significant. 
For example, by corollary~\ref{c-19},
$\spo(\H)\times\nobreak\spo(\I)\cong\spo(\H)\gdot\spo(\I)$ which is proper. 
This can easily be extended arbitrary finite products, 
whence $\spo(\H_0)\times\cdots\times\spo(\H_{n-1})$ is proper. 
This strongly suggests that Shelah's NNR theory
from~\cite[Ch.~{XVIII}, \Section2]{MR1623206} applies to our
classes of forcing notions (see section~\ref{sec:trind-properness} for
more discussion). 
This would be the first instance we are aware of
where the theory applies to forcing notions of cardinality $\aleph_2$
or greater. All of the examples 
in~\cite[Ch.~{XVIII}, \Section1,2]{MR1623206} are forcing
notions of cardinality $\aleph_1$. 
\end{remark}

Note for example that, at least when dealing with $P$-ideals, corollary~\ref{c-6}
strengthens to:

\begin{cor}[$\ch$]
 \label{c-20}
Let $\H$ and $\I$ be $P$-ideals on $\oone$.
Then $\spo(\H)\forces\spo(\I)^V\cong\spo(\I)$\tu; hence,
$\spo(\H)\forces\ulc\spo(\I)^V$\tu{ is densely included in }$\spo(\I)\urc$.
\end{cor}
\begin{pf}
By corollary~\ref{c-19} and proposition~\ref{p-68}. 
\end{pf}

\fxnote{comment}

We are going to extend e.g.~corollary~\ref{c-10} to countable support
iterations. For example, we shall prove: \fxnote{Need to get
  hypothesis right. Later ...}

\begin{thm}[$\ch$]
\label{u-5}
Suppose that $(P_\xi,\spo(\dot\H_\xi):\xi<\delta)$ is a countable
support iteration, where each $\dot\H_\xi$ is a $P_\xi$-name for
a\/~$P$-ideal on $\oone$, 
and $\I$ is a\/~$P$-ideal on $\oone$ with no countable decomposition of
$\oone$ into orthogonal pieces. Let 
Then the limit $P_\delta$ of the iteration forces that $\I$ has no
countable decomposition of $\oone$ into pieces orthogonal to $\I$. 
\end{thm}

\noindent We do not however obtain a preservation theorem for
countable support iterations not decomposing $\oone$ into countably
many pieces orthogonal to $\ideal$, and we doubt that this property is
preserved under the iteration of any general class of proper forcing
notions (as opposed to the specific class $\spo$). 

\subsection{Coding iterations}
\label{sec:coding-iterations}

While the forcing notions $\poset$ and $\spo$ 
are viewed as classes with one parameter,
we need to generalize definability to iterations, to also allow iterations of
$\poset$ and $\spo$ to be interpreted in the relevant model. This is
necessary for our analysis of embedability, and will be necessary
for our handling of the $\nnr$ iteration as well. 

\begin{defn}
\label{d-26}
Let $\theta$ be an ordinal of uncountable cofinality.
We describe a coding of those iterations
consisting of combinations of the forcing notions $\poset(\H)$ and
$\spo(\H)$, with $\H\subseteq\powcif\theta$.
We define a class  $\code\theta$ of sequences, or codes,
and forcing notions $P(\vec a)$ for each $\vec a\in\code\theta$, by
recursion on $\xi=\length(\vec a)$.
Let $\code\theta\restriction0$ be the singleton containing the null sequence $\<\>$ 
and let $P(\emptyset)$ be the trivial forcing notion.
Having defined $\code\theta\restriction\xi$, let $\code\theta\restriction\xi+1$ 
be the collection of all sequences of the form $\vec
a\bigexta(\dot\H,\O)$  where $\vec a\in\code\theta\restriction\xi$,
 $\dot\H$ is a $P(\vec a)$-name for a $\sigma$-directed subfamily of  $(\powcif\theta,\subseteqfnt)$ 
and $\O$ is either $\poset$ or $\spo$; then let
\begin{equation}
  \label{eq:85}
  P(\vec a\bigexta(\dot\H,\O))=
  \begin{cases}
    P(\vec a)\gdot\poset(\dot\H),&\text{if $\O=\poset$,}\\
    P(\vec a)\gdot\spo(\dot\H),&\text{if $\O=\spo$.}
    \end{cases}
\end{equation}
For limit $\delta$, let $\code\theta\restriction\delta=\liminv\xi\delta\code\theta\restriction\xi$ be the inverse limit, 
i.e.~all sequences $\vec a$ of length $\delta$  with $\vec a\restriction\xi\in\code\theta$ for all $\xi<\delta$;
then for each $\vec a\in\code\theta\restriction\delta$, we let $P(\vec a)$ be the
corresponding countable support iteration. 
Thus $P(\vec a)$ is the limit of $(P_\xi:\xi<\delta)$ of the iterated forcing $(P_\xi,\dot Q_\xi:\xi<\delta)$,
where each $\dot Q_\xi$ is the second iterand in
equation~\eqref{eq:85} plugging in $\vec a:=\vec a\restriction\xi$ and $(\dot\H,\O):=\vec a(\xi)$, 
inverse limits are taken at limits of countable cofinality 
and direct limits are taken at limits of uncountable cofinality. 
Denote the class $\code\theta=\bigcup_{\xi\in\ord}\code\theta\restriction\xi$. 
For each $\vec a\in\code\theta$ and each $\xi<\length(\vec a)$, 
we let $\dot\H(\vec a(\xi))=\dot\H$ where $\vec a(\xi)=(\dot\H,\O)$.

Let $\pcode\theta\subseteq\code\theta$ be the set of all codes $\vec
a$ such that for all $\xi<\length(\vec a)$,
if $\vec a(\xi)$ is of the form $(\dot\H,\poset)$ then 
$P(a\restriction\xi)\forces\ulc$there is no stationary subset of $\theta$ 
orthogonal to $\dot\H\urc$.

For $\vec a\in\code\theta$, let $D(\vec a)$ be the set of codes
generated by the operation $\vec a\bigexta(\dot\H,\spo)$, where
$P(\vec a)$ forces that $\dot\H$ is $\sigma$-directed 
(and thus $\length(\vec b)<\length(\vec a)+\omega$ for all $\vec b\in\nobreak D(\vec a)$). Then define
$\qcode\theta\subseteq\code\theta$ as the set of all codes $\vec a$ such that
for all $\xi<\length(\vec a)$, if $\vec a(\xi)$ is of the form
$(\dot\H,\poset)$, then for all $\vec b\in\nobreak D(\vec a\restriction\xi)$,
 $P(\vec b)$ forces that there is no stationary
subset of $\theta$ orthogonal to $\dot\H$. 

We also define $\gcode\theta\subseteq\qcode\theta$ as the set of all
codes $\vec a$ such that for all $\xi<\length(\vec a)$, if $\vec
a(\xi)$ is of the form $(\dot\H,\poset)$, 
then for every $\vec c\in
D(\vec a\restriction\xi)$ and every $\vec b\in C(\vec c)$ with 
 $\vec a\restriction\xi\subseteq\vec b$ (cf.~definition~\ref{d-29}),
 $P(\vec b)$ forces that there is no stationary subset of $\theta$
 orthogonal to $\dot\H$. 

Define $\codep\theta\subseteq\code\theta$ as the set of all codes $\vec a$ 
such that $P(\vec a\restriction\xi)\forces\ulc\dot\H(\vec a(\xi))$  is a $P$-ideal$\urc$ for all $\xi<\length(\vec a)$. 
Let $\pcodep\theta=\pcode\theta\cap\codep\theta$, $\qcodep\theta=\qcode\theta\cap\codep\theta$ 
and   $\gcodep\theta=\gcode\theta\cap\codep\theta$. 
\end{defn}

\begin{prop}
\label{p-25}
If $\vec a,\vec b\in\code\theta$ then $\vec a\bigexta\vec
b\in\code\theta$  and $P(\vec a\bigexta\vec b)=P(\vec a)\gdot P(\vec
b)$. 
More generally, if $(\vec a_\xi:\xi<\mu)$ is a sequence of elements of
$\code\theta$, then the concatenation 
$\vec b=\vec a_0\bigexta\vec a_1\bigexta\cdots\bigexta\vec
a_\xi\bigexta\cdots$ is in $\code\theta$ and $P(\vec b)$ is the limit of
the countable support iteration determined by $(P(\vec
a_\xi):\xi<\mu)$. 
\end{prop}

\begin{lem}
\label{l-74}
For all $\vec a\in\pcode\theta$, $P(\vec a)$ is proper. 
\end{lem}
\begin{pf}
By lemmas~\ref{lem:proper:b} and~\ref{l-24}. 
\end{pf}

\begin{prop}
\label{p-65}
$\gcode\theta\subseteq\qcode\theta\subseteq\pcode\theta$.
\end{prop}
\begin{pf}
For all $\xi<\length(\vec a)$, $\vec a\restriction\xi\in D(\vec
a\restriction\xi)$ and $\vec a\restriction\xi\in C(\vec
a\restriction\xi)$. 
\end{pf}

\begin{prop}
\label{p-66}
For all $\vec a\in\code\theta$, for all $\xi<\length(\vec a)$,
$P(\vec a\restriction\xi)\forces\vec a\restriction[\xi,\length(\vec a))\in\code\theta$,
i.e.~we are taking $\code\theta$ as a class with parameter $\theta$ that is being
interpreted in the forcing extension by $P(\vec a\restriction\xi)$. 
Similarly, for all $\vec a\in\pcode\theta$ \tu($\qcode\theta$\tu) 
\tu[$\codep\theta$\tu], for all $\xi<\length(\vec a)$,
$P(\vec a\restriction\xi)\forces\vec a\restriction[\xi,\length(\vec a))\in\pcode\theta$
\tu($\qcode\theta$\tu) \tu[$\codep\theta$\tu].
\end{prop}
\begin{pf}
These are immediate from the associativity of iterated forcing. 
\end{pf}

\begin{remark}
\label{r-10}
Proposition~\ref{p-66} may fail for $\gcode\theta$, because in some
forcing extension by $P(\vec a\restriction\xi)$ there may be new
elements of $C(\dot{\vec c})$ that do not correspond to elements of~$C(\vec a)^V$, 
because for example elements of $C(\dot{\vec c})$
include uncountable concatenations. 
\end{remark}

We also have a converse. 

\begin{prop}
\label{p-67}
For all $\vec a\in\code\theta$, and every $P(\vec a)$-name $\dot{\vec
  c}$, if $P(\vec a)\forces\dot{\vec c}\in\code\theta$ 
then $\vec a\bigexta\dot{\vec c}\in\code\theta$ \tu(assuming a suitable
representation of $\dot{\vec c}$\tu). Similarly, for $\pcode\theta$, $\qcode\theta$ and $\codep\theta$. 
\end{prop}

Now we can generalize corollary~\ref{l-37} using our coding of
iterations in the definition of frozen
(definition~\ref{d-24}).\fxwarning{Check that this is what wanted to
  say.}

\begin{lem}[$\ch$]
\label{l-45}
Let $\I$ be a $P$-ideal on $\oone$ and let $\vec a\in\qcodep\oone$.
If $P(\vec a)$ adds no new reals,  then all of the following are true\tu:
\begin{enumerate}[label=\tu{(\alph*)}, ref=\alph*, widest=b]
\item\label{item:18} If there is no countable decomposition of $\oone$ into pieces
  orthogonal to $\I$, then $P(\vec a)$ forces that $\I$ 
  has no countable decomposition of $\oone$ into orthogonal pieces.
\item\label{item:19}  $P(\vec a)\forces \spo(\I)^V\imbeds[V]\spo(\I)$\tu; 
  hence\tu, $\spo(\I)\embeds P(\vec a)\gdot\spo(\I)$.
\item\label{item:21}   $\spo(\I)\forces P(\vec a)^V\imbeds[V] P(\vec a)$\tu;
  hence\tu, $P(\vec a)\embeds\spo(\I)\gdot P(\vec a)$.
\item\label{item:22}  $P(\vec a)\gdot\spo(\I)\div\spo(\I)\cong P(\vec a)$.
\item\label{item:23} $P(\vec a)\gdot\spo(\I)\cong P(\vec a)\times\spo(\I)
  \cong\spo(\I)\times P(\vec a)$.
\item\label{item:55}  $P(\vec a)\gdot\spo(\I)\cong\spo(\I)\gdot P(\vec a)$.
\item\label{item:155} $\spo(\I)\forces\ulc P(\vec a)^V$\tu{ is densely
  included in }$P(\vec a)\urc$.
\item\label{item:154}  $(\I,\spo)\bigexta\vec a\in\pcodep\oone$.
\item\label{item:73}   Let $G\in\Gen(V,\spo(\I))$ and 
$G\gdot H\in\Gen(V,\spo(\I)\gdot P(\vec a))$.
Then for every $P(\vec a)$-name $\dot\H$ for a $\sigma$-directed subfamily
of $\powcif\oone$, if \tu{$V[H]\models\ulc\dot\H[H]$ has no countable
decomposition into orthogonal sets$\urc$}, 
then so does \tu{$V[G\gdot H]\models\ulc \dot\H[H]$ has no 
such countable decomposition of~$\oone\urc$}.
\end{enumerate}
\end{lem}
\begin{pf}
All clauses~\eqref{item:18}--\eqref{item:73} 
are proved simultaneously by induction on $\length(\vec a)$. 

\smallskip
\noindent{\emph{Base case}:} $\length(\vec a)=0$.
\smallskip

\noindent $P(\vec a)$ is the trivial forcing notion. 
Thus clauses~\eqref{item:18}--\eqref{item:154} are trivial, while~\eqref{item:73}
reduces to corollary~\ref{c-10}.

\smallskip
\noindent{\emph{Successor case}:} $\length(\vec a)=\xi+1$. 
\smallskip

\noindent Then $\vec a$ is either 
of the form $\vec b\bigexta(\dot\J,\spo)$ or $\vec b\bigexta(\dot\J,\poset)$, 
i.e.~$P(\vec a)$ is either of the form $P(\vec b)\gdot\spo(\dot\J)$ 
or $P(\vec b)\gdot\poset(\dot\J)$ (possibly $\vec b=\<\>$). 

For clause~\eqref{item:18}, 
let $G\in\Gen(V,P(\vec b))$ and let $G\gdot H\in\Gen(V,P(\vec a))$. 
By the induction hypothesis, $\I$ has no countable decomposition,
in $V[G]$, of $\oone$ into orthogonal pieces. 
In the first case $\vec a=\vec b\bigexta(\dot\J,\spo)$, 
applying corollary~\ref{c-10} in $V[G]$, this remains true in
$V[G\gdot H]$. In the other case $\vec a=\vec b\bigexta(\dot\J,\poset)$.
Then, in $V[G]$, $\dot\J[G]$ is a $P$-ideal since $\vec a\in\codep\oone$; 
and $\spo(\I)$ forces there is no stationary set
orthogonal to $\dot\J[G]$, because, in $V$, $\vec a\in\qcode\oone$ and
thus $\vec a\restriction\xi\bigexta(\I,\spo)\in D(\vec
a\restriction\xi)$ implies that $P(\vec b)\gdot\spo(\I)$ forces there
is no stationary subset of $\oone$ orthogonal to $\dot\J$. 
Therefore, corollary~\ref{c-15} applies in $V[G]$, establishing that
in $V[G\gdot H]$ there is no countable decomposition of $\oone$ into
pieces orthogonal to $\I$. 

Clause~\eqref{item:19} follows from clause~\eqref{item:18},
just as in the proof of corollary~\ref{c-6}. 

For clause~\eqref{item:21}, given a maximal antichain $A\subseteq P(\vec a)$, 
we need to show that $\spo(\I)\forces\ulc A$ is a maximal antichain of $P(\vec a)\urc$. 
First suppose $\vec a$ is of the form $\vec b\bigexta(\dot\J,\spo)$. 
Fix $I\in\Gen(V,\spo(\I))$.  
Then take $J\in\Gen(V[I],P(\vec b)^{V[I]})$, so
that $I\gdot\nobreak J\in\Gen(V,\spo(\I)\gdot P(\vec b))$. 
By proposition~\ref{p-39}, 
we have that $P(\vec b)$ forces $A\div P(\vec b)$ is a maximal
antichain of $\spo(\dot\J)$. Therefore, by the induction hypothesis that
clause~\eqref{item:21} holds for $\vec b$, $J\in\Gen(V,P(\vec b))$ and hence
putting $B=(A\div P(\vec b))[J]$, $V[J]\models\ulc B$ is a maximal
antichain of $\spo(\dot\J[J])\urc$. 
We  apply corollary~\ref{l-37} with $V:=V[J]$, 
$\H:=\dot\J[J]$, $A:=B$ and $W:=V[I\gdot J]$.
For any $x\in\dot\J[J]$, applying the induction hypothesis that 
clause~\eqref{item:73} holds for~$\vec b$, with
$\dot\J:=\Psi(\dot\J_{(x)})$,\fxwarning{Should note somewhere that
  $P$-ideal on set of size $\aleph_1$.} we see that if, in $V[J]$, there is no
countable decomposition into sets orthogonal to 
$\Psi(\dot\J_{(x)})[J]=\Psi(\dot\J[J]_{(x)})$, 
then, in $V[I\gdot J]$, there also no such decomposition.
Therefore, corollary~\ref{l-37} 
yields $V[I\gdot J]\models\ulc (A\div
P(\vec b))[J]$ is a maximal antichain of $\spo(\dot\J[J])\urc$. 
Since $J$ is arbitrary, this proves that $V[I]\models\ulc A$ is a maximal
antichain of $\spo(\vec a)\urc$ by proposition~\ref{p-39}, as desired. The other case where
$\vec a=\vec b\bigexta (\dot\J,\poset)$, is exactly the same but
corollary~\ref{c-16} is used instead.

For~\eqref{item:22}, 
the hypothesis of lemma~\ref{l-73}, with $P:=P(\vec a)$, is satisfied
because $P(\vec a)$ adds no new countable subsets of $\oone$ and by
clause~\eqref{item:19}. Then clause~\eqref{item:22} is the conclusion
of the lemma.

Clause~\eqref{item:23} is a restatement of clause~\eqref{item:22}
together with the fact that products commute.

Clause~\eqref{item:55} is proved algebraically.
First consider $\vec a=\vec b\bigexta(\dot\J,\spo)$. 
Since $P(\vec b)$ adds no new reals, $P(\vec b)\forces \ch$, 
and thus
\begin{equation}
P(\vec b)\forces \spo(\dot\J)\gdot\spo(\I)\cong\spo(\I)\gdot\spo(\dot\J)\label{eq:129}
\end{equation}
by applying corollary~\ref{c-19} in this forcing extension. 
Now by associativity of iterated forcing for the first equivalence, 
by equation~\eqref{eq:129} for the second equivalence, and by the
induction hypothesis that~\eqref{item:55} holds for $\vec b$ for the
third equivalence,
\begin{equation}
  \label{eq:130}
  \begin{split}
    P(\vec a)\gdot\spo(\I)&\cong P(\vec b)\gdot[\spo(\dot\J)\gdot\spo(\I)]\\
    &\cong [P(\vec b)\gdot\spo(\I)]\gdot\spo(\dot\J)\\
    &\cong \spo(\I)\gdot [P(\vec b)\gdot\spo(\dot\J)]\\
    &=\spo(\I)\gdot P(\vec a),
  \end{split}
\end{equation}
as required. 

Now we consider 
the other case $\vec a=\vec b\bigexta(\dot\J,\poset)$. If $\dot\H$ is
a $P(\vec b)$-name for a $\sigma$-directed family with no countable
decomposition of $\oone$ into orthogonal pieces, 
then $\vec b\bigexta(\dot\H,\spo)\in D(\vec b)$ and 
hence $P(\vec b)\forces\spo(\dot\H)\forces\ulc$there is no
stationary set orthogonal to $\dot\J\urc$ because
$\vec a\in\qcode\oone$. Therefore,
the hypothesis of corollary~\ref{c-23} holds in the extension by
$P(\vec b)$, and hence by the corollary,
\begin{equation}
  \label{eq:131}
  P(\vec b)\forces \poset(\dot\J)\gdot\spo(\I)\cong\spo(\I)\gdot\poset(\dot\J).
\end{equation}
Now we can obtain the result in exactly the same manner as
equation~\eqref{eq:130}.

Clause~\eqref{item:155} is an immediate consequence
of~\eqref{item:23} and~\eqref{item:55} and proposition~\ref{p-68}.

For~\eqref{item:154}, put $\vec c=(\I,\spo)\bigexta\vec a$. 
First of all note that $\vec c\in\code\oone$ by clause~\eqref{item:21}. 
Take $\xi<\length(\vec c)$. We can assume $\xi>0$ since $\vec c(0)=(\I,\spo)$
is not of the form $(\dot\H,\poset)$, say $\xi=1+\eta$.
Then $\eta<\length(\vec a)$. 
We have to deal with the situation where $\vec a(\eta)$ is of the form
$(\dot\H,\poset)$, in which case we must show
that $\spo(\I)\gdot P(\vec a\restriction\eta)$ forces there is no
stationary set orthogonal to $\dot\H$. Applying the induction
hypothesis that clause~\eqref{item:55} holds for 
$\vec a\restriction\eta$, $\spo(\I)\gdot P(\vec a\restriction\eta)\cong
P(\vec a\restriction\eta)\gdot\spo(\I)$. 
Now $\vec d=(\vec a\restriction\eta)\bigexta
(\I,\spo)\in D(\vec a\restriction\eta)$, 
and thus $P(\vec d)$ forces there is no stationary set orthogonal to $\dot\H$,  because $\vec a\in\qcode\oone$. 
Since $P(\vec c\restriction\xi)=\spo(\I)\gdot P(\vec
a\restriction\eta)\cong P(\vec d)$,  this concludes the proof that $\vec c\in\pcode\oone$.

\fxnote{attempting to prove in $\qcode\oone$.}

For clause~\eqref{item:73}, let $\dot\H$ be a $P(\vec a)$-name for a
$\sigma$-directed family. Let $G\in\Gen(V,\allowbreak\spo(\I))$ and
$H\in\Gen(V[G],P(\vec a)^{V[G]})$ 
(thus $G\gdot H\in\Gen(V,\spo(\I)\gdot P(\vec a))$).
Then $H\in\Gen(V,P(\vec a))$ by clause~\eqref{item:21}.
We assume that, in $V[H]$, there is no countable decomposition
into sets orthogonal to $\dot\H[H]$. 
But then by~\eqref{item:55}, we know that 
$V[G\gdot H]$ is an $\spo(\I)$-generic extension of $V[H]$, 
and therefore there is no countable decomposition, in $V[G\gdot H]$,
into sets orthogonal to $\dot\H[H]$ by corollary~\ref{c-10}. 

\fxnote{Might be better approach to proving 2nd lemma: By the induction hypothesis, we have $\spo(\I)^{V}\imbeds[V]\spo(\I)^{V[G]}$. 
And applying corollary~\ref{c-6} in $V[G]$, 
we obtain $\spo(\I)^{V[G]}\imbeds[{V[G]}]\spo(\I)^{V[G\gdot H]}$.}

\smallskip
\noindent\emph{Limit case}: $\length(\vec a)$ equals some limit ordinal $\delta$. 
\smallskip

First we establish clause~\eqref{item:21}. 
Let $G\in\Gen(V,\spo(\H))$. In $V[G]$:
$P(\vec a)$ is the limit of $(P(\vec a\restriction\xi):\xi<\delta)$
(proposition~\ref{p-25}).  And by the induction hypothesis,   
$P(\vec a\restriction\xi)^V\imbeds[V]P(\vec a\restriction\xi)$ for all
$\xi<\delta$. 
Therefore, by lemma~\ref{l-38}, the limit, let us call it $Q$, 
of $(P(\vec a\restriction\xi)^V:\xi<\delta)$ is generically included
in $P(\vec a)$ over $V$. Since we are dealing with countable support
iterations, and since $\spo(\H)$ adds no new reals, the limit of
$(P(\vec a\restriction\xi)^V:\xi<\delta)$ is the same whether taken
here in $V[G]$ or in the ground model $V$. Hence $P(\vec
a)^V=Q\imbeds[V] P(\vec a)^{V[G]}$. 

For clause~\eqref{item:154}, we first of all have
$(\I,\spo)\bigexta\vec a\in\codep\oone$ by clause~\eqref{item:21}.
It then follows immediately from the induction hypothesis that
$(\I,\spo)\bigexta P(\vec a\restriction\xi)\in\pcode\oone$ for all
$\xi<\delta$, that $(\I,\spo)\bigexta P(\vec a)\in\pcode\oone$.

Next we deal with clause~\eqref{item:18}. If there is no countable
decomposition of $\oone$ into pieces orthogonal to $\I$, then 
$\spo(\I)$ forces an uncountable set locally in $\I$. 
Moreover, $\spo(\I)\gdot P(\vec a)$ is proper, and in particular does
not collapse $\aleph_1$, by lemma~\ref{l-74}, 
because $(\I,\spo)\bigexta\vec a\in\pcode\oone$ by
clause~\eqref{item:154}. 
Therefore, as $P(\vec a)\embeds\spo(\I)\gdot P(\vec a)$ by~\eqref{item:21}, 
$P(\vec a)$ cannot force a countable decomposition of $\oone$
into pieces orthogonal to $\I$.

Clause~\eqref{item:19} follows from clause~\eqref{item:18} as before.

Clause~\eqref{item:22} follows from clause~\eqref{item:19} exactly as
in the successor case; similarly for clause~\eqref{item:23}. 

For clause~\eqref{item:55}, let $G\in\Gen(V,\spo(\I))$. 
In $V[G]$: By the induction hypothesis that clause~\eqref{item:155} 
holds for $\vec a\restriction\xi$ for all $\xi<\length(\vec a)$, 
we have that $P(\vec a\restriction\xi)^V$ is densely included in
$P(\vec a\restriction\xi)$ for all $\xi<\delta$. 
Therefore, by lemma~\ref{l-33}, the limit of $(P(\vec
a\restriction\xi):\xi<\delta)$, 
call it $Q$, is isomorphic as a
forcing notion to $P(\vec a)$. Since $\spo(\I)$ adds no reals, and the
iterations are of countable support, $P(\vec a)^V=Q\cong P(\vec a)$.
Now, back in $V$, we have established that $\spo(\I)\times P(\vec
a)\cong\spo(\I)\gdot P(\vec a)$, 
and thus the result is now a consequence of~\eqref{item:23}.

Clause~\eqref{item:155} follows as for the successor case. 

Clause~\eqref{item:73} follows from~\eqref{item:21},~\eqref{item:55} 
and corollary~\ref{c-10} identically as for the successor case. 
\end{pf}

\begin{lem}[$\ch$]
\label{l-46}
Let $\I$ be a $P$-ideal on $\oone$ and let $\vec a\in\qcodep\oone$.
If $P(\vec a)$ adds no new reals,  then all of the following are true\tu:
\begin{enumerate}[label=\tu{(\alph*)}, ref=\alph*, widest=b]
\item\label{item:159}  $P(\vec a)\forces \poset(\I)^V\imbeds[V]\poset(\I)$\tu; 
  hence\tu, $\poset(\I)\embeds P(\vec a)\gdot\poset(\I)$.
\item\label{item:160} If $\vec a\bigexta(\I,\poset)\in\qcode\oone$ 
  then  $\poset(\I)\forces P(\vec a)^V\imbeds[V] P(\vec a)$\tu;
  hence\tu, $P(\vec a)\embeds\poset(\I)\gdot P(\vec a)$.
\item\label{item:151}  $P(\vec a)\gdot\poset(\I)\div\poset(\I)\cong P(\vec a)$.
\item\label{item:152} $P(\vec a)\gdot\poset(\I)\cong P(\vec a)\times\poset(\I)
  \cong\poset(\I)\times P(\vec a)$.
\item\label{item:153} If $\vec a\bigexta(\I,\poset)\in\qcode\oone$ 
then  $P(\vec a)\gdot\poset(\I)\cong\poset(\I)\gdot P(\vec a)$.
\item\label{item:156} If $\vec a\bigexta(\I,\poset)\in\qcode\oone$ 
then  $\poset(\I)\forces\ulc P(\vec a)^V$\tu{ is densely  included in }$P(\vec a)\urc$.
\item\label{item:158} If $\vec a\bigexta(\I,\poset)\in\qcode\oone$ 
then $(\I,\spo)\bigexta\vec a\in\pcodep\oone$.
\item\label{item:157}   Let $G\in\Gen(V,\poset(\I))$ and 
$G\gdot H\in\Gen(V,\poset(\I)\gdot P(\vec a))$.
Then for every $P(\vec a)$-name $\dot\H$ for a $\sigma$-directed subfamily
of $\powcif\oone$, if
\begin{enumerate}[label=\tu{(\arabic*)}, ref=\arabic*, widest=2]
\item\label{item:161} \tu{$V[H]\models\ulc\dot\H[H]$ has no countable
orthogonal decomposition$\urc$},
\item\label{item:162} \tu{$V[H]\models\ulc\spo(\dot\H[H])$ forces
    that there is no stationary set orthogonal to $\ideal\urc$}, 
\end{enumerate}
then so does \tu{$V[G\gdot H]\models\ulc \dot\H[H]$ has no 
countable orthogonal decomposition of~$\oone\urc$}.
\end{enumerate}
\end{lem}
\begin{pf}
The proof is by induction on $\length(\vec a)$. 
The base case $\length(\vec a)=0$ is completely straightforward, and
the limit case is the same as for the proof of lemma~\ref{l-45}. Hence
we only deal with the successor case $\length(\vec a)=\xi+1$. 

Clause~\eqref{item:159} follows from lemma~\ref{l-45}\eqref{item:18},
just as in the proof of corollary~\ref{c-14}. 

For clause~\eqref{item:160}, given a maximal antichain $A\subseteq P(\vec a)$, 
we need to show that $\poset(\I)\forces\ulc A$ is a maximal antichain of $P(\vec a)\urc$. 
First suppose $\vec a$ is of the form $\vec b\bigexta(\dot\J,\spo)$. 
Fix $I\in\Gen(V,\poset(\I))$.  
Then take $J\in\Gen(V[I],P(\vec b)^{V[I]})$, so that $I\gdot\nobreak J\in\Gen(V,\poset(\I)\gdot P(\vec b))$. 
By proposition~\ref{p-39}, 
we have that $P(\vec b)$ forces $A\div P(\vec b)$ is a maximal
antichain of $\spo(\dot\J)$. 
Therefore, by the induction hypothesis that
clause~\eqref{item:160} holds for $\vec b$, $J\in\Gen(V,P(\vec b))$ and hence
putting $B=(A\div P(\vec b))[J]$, $V[J]\models\ulc B$ is a maximal antichain of $\spo(\dot\J[J])\urc$. 
 We  apply corollary~\ref{l-37} with $V:=V[J]$, 
$\H:=\dot\J[J]$, $A:=B$ and $W:=V[I\gdot J]$.
For any $x\in\dot\J[J]$, suppose that, in $V[J]$, there is no
countable decomposition into sets orthogonal to $\Psi(\dot\J[J]_{(x)})$. 
Since $\vec b\bigexta(\Psi(\dot\J_{(x)}),\spo)\in D(\vec b)$, 
and since the hypothesis on $\I$ clearly entails that 
$\vec b\bigexta(\I,\poset)\in\qcode\oone$, we have that $P(\vec
b)\gdot\spo(\Psi(\dot\J_{(x)}))$ forces there is no stationary set
orthogonal to $\I$, i.e.~$V[J]\models\ulc\spo(\Psi(\dot\J[J]_{(x)}))$ 
forces there is no stationary set orthogonal to $\I\urc$. 
Therefore, the induction hypothesis~\eqref{item:157} applies to $\vec b$ with
$\dot\H:=\nobreak\Psi(\dot\J_{(x)})$, and thus, in $V[I\gdot J]$, 
there is also no countable decomposition into sets orthogonal to $\Psi(\dot\J[J]_{(x)})$.
Therefore, corollary~\ref{l-37} 
yields $V[I\gdot J]\models\ulc (A\div
P(\vec b))[J]$ is a maximal antichain of $\spo(\dot\J[J])\urc$. 
Since $J$ is arbitrary, this proves that $V[I]\models\ulc A$ is a maximal
antichain of $\spo(\vec a)\urc$ by proposition~\ref{p-39}, as desired. The other case where
$\vec a=\vec b\bigexta (\dot\J,\poset)$, is exactly the same but corollary~\ref{c-16} is used instead.

Clause~\eqref{item:151} is a consequence of lemma~\ref{l-72} with
$P:=P(\vec a)$, by the hypothesis that $P(\vec a)$ does not add reals
and clause~\eqref{item:159}. 

Clause~\eqref{item:152} is a restatement of clause~\eqref{item:151}.

Clause~\eqref{item:153} is proved algebraically.
First consider $\vec a=\vec b\bigexta(\dot\J,\spo)$. 
Since $P(\vec b)$ adds no new reals, $P(\vec b)\forces \ch$;
and for every $P(\vec b)$-name $\dot\H$ for a $\sigma$-directed
family, it follows from the fact that $\vec
a\bigexta(\I,\poset)\in\qcode\oone$ that $P(\vec
b)\forces\spo(\dot\H)\forces\ulc$there is no stationary set orthogonal
to $\I\urc$; 
and thus
\begin{equation}
\label{eq:73}
P(\vec b)\forces \spo(\dot\J)\gdot\poset(\I)\cong\poset(\I)\gdot\spo(\dot\J)
\end{equation}
by applying corollary~\ref{c-23} the forcing extension by $P(\vec b)$.
Using by equation~\eqref{eq:73} for the second equivalence, and 
the  induction hypothesis~\eqref{item:153} for $\vec b$ for the third equivalence,
\begin{equation}
  \label{eq:78}
  \begin{split}
    P(\vec a)\gdot\poset(\I)&\cong P(\vec b)\gdot[\spo(\dot\J)\gdot\poset(\I)]\\
    &\cong [P(\vec b)\gdot\poset(\I)]\gdot\spo(\dot\J)\\
    &\cong \poset(\I)\gdot [P(\vec b)\gdot\spo(\dot\J)]\\
    &=\poset(\I)\gdot P(\vec a),
  \end{split}
\end{equation}
as required. 

Now we consider 
the other case $\vec a=\vec b\bigexta(\dot\J,\poset)$. If $\dot\H$ is
a $P(\vec b)$-name for a $\sigma$-directed family,
then $\vec b\bigexta(\dot\H,\spo)\in D(\vec b)$ and 
hence $P(\vec b)\forces\spo(\dot\H)\forces\ulc$there is no
stationary set orthogonal to $\dot\J\urc$ because
$\vec a\in\qcode\oone$; and furthermore, we saw above that $P(\vec
b)\forces\spo(\dot\H)\forces\ulc$there is no stationary set orthogonal to $\I\urc$. 
Therefore, the hypothesis of corollary~\ref{c-24} holds 
in the extension by $P(\vec b)$, and hence by the corollary,
\begin{equation}
  \label{eq:140}
  P(\vec b)\forces \poset(\dot\J)\gdot\poset(\I)\cong\poset(\I)\gdot\poset(\dot\J).
\end{equation}
Now we can obtain the result in exactly the same manner as
equation~\eqref{eq:78}.

Clause~\eqref{item:156} is an immediate consequence of proposition~\ref{p-68}.

For~\eqref{item:158} is immediate from~\eqref{item:153}.

For clause~\eqref{item:157}, let $\dot\H$ be a $P(\vec a)$-name for a
$\sigma$-directed family. Let $G\in\Gen(V,\allowbreak\poset(\I))$ and
$H\in\Gen(V[G],P(\vec a)^{V[G]})$. 
Then $H\in\Gen(V,P(\vec a))$ by clause~\eqref{item:160};
and by~\eqref{item:153}, we know that 
$V[G\gdot H]$ is a $\poset(\I)$-generic extension of $V[H]$. 
We assume that~\eqref{item:161} and~\eqref{item:162} hold, 
and therefore, in $V[H]$, the hypotheses of corollary~\ref{c-15} hold
with $\H:=\dot\H[H]$, and hence,
in $V[G\gdot H]$, there is no countable decomposition
into sets orthogonal to $\dot\H[H]$ by the corollary.  
\end{pf}

The following theorem is the absoluteness result we have been working
towards. 

\begin{thm}[$\ch$]
\label{u-10}
Let $\vec a,\vec b\in\codep\oone$.  Suppose\/ $\vec a\bigexta\vec b\in\qcode\oone$. 
If $P(\vec a\bigexta \vec b)$ adds no new reals,
then $P(\vec a)\forces P(\vec b)^V\imbeds[V] P(\vec
b)$, and hence $P(\vec b)\embeds P(\vec a\bigexta\vec b)$. 
Moreover, 
\begin{equation}
\begin{split}
P(\vec a\bigexta\vec b)
=P(\vec a)\gdot P(\vec b)
&\cong P(\vec a)\times P(\vec b)\\
&\cong P(\vec b)\times P(\vec a)
\cong P(\vec b)\gdot P(\vec a)
= P(\vec b\bigexta\vec a).
\end{split}
\end{equation} 
\end{thm}
\begin{pf}
This is proved by a straightforward induction from lemmas~\ref{l-45} and~\ref{l-46}. 
\end{pf}

\fxnote{nice notation}

\fxnote{comment ???}

Let us next describe how theorem~\ref{u-10} is applied, after
introducing notation for concatenating sequences of sequences.

\begin{defn}
\label{d-4}
For $X\subseteq\ord$ and any sequence $\vec{\vec x}=(\vec
x_\gamma:\gamma\in X)$ of \emph{sequences} (i.e.~functions whose domains
are ordinals), 
let $\rho(\vec{\vec x})$ be the concatenation under the ordinal ordering,
i.e.~$\rho(\vec{\vec x})$ 
is a sequence of length $\sum_{\gamma\in X}\length(\vec x_\gamma)$ and 
$\rho(\vec{\vec x})\restriction[\zeta_\gamma,\zeta_{\gamma+1})
=\vec x_{\gamma}$ for all $\gamma\in X$, 
where $\zeta_\gamma=\sum_{\xi\in X,\xi<\gamma}\length(\vec x_\xi)$. 
\end{defn}

\begin{prop}
\label{p-1}
Suppose that $\vec{\vec a}=(\vec a_\gamma:\gamma\in X)$ where each
$\vec a_\gamma\in\code\theta$. Then
every $p\in P(\rho(\vec{\vec a}))$ is of the form $\rho(\vec p)$ where
$\vec p=(p_\gamma:\gamma\in X)$ and each $p_\gamma\in P(\vec a_\gamma)$.
\end{prop}

\begin{defn}
\label{d-1}
For all $\vec a,\vec b\in\code\theta$, let $e(\vec a,\vec b):P(\vec
b)\to P(\vec a\bigexta\vec b)$ be given by
\begin{equation}
  \label{eq:3}
  e(\vec a,\vec b)(p)={0_{P(\vec a)}}\bigext p
\end{equation}
for all $p\in P(\vec b)$.

More generally, suppose that $\vec{\vec a}=(\vec a_\gamma:\gamma<\delta)$  
is a sequence with each $\vec a_\gamma\in\code\theta$. 
For $X\subseteq\delta$, let $f(\vec{\vec a},X):P(\rho(\vec{\vec
  a}\restriction X))\to P(\rho(\vec{\vec a}))$ be given by
$f(\vec{\vec a},X)(p)=\rho(\vec q)$ where $\vec
q=(q_\gamma:\gamma<\delta)$ is given by
\begin{equation}
  \label{eq:8}
  q_\gamma=\begin{cases}
    p_\gamma,&\text{if $\gamma\in X$},\\
    0_{P(\vec a_\gamma)},&\text{if $\gamma\notin X$},
    \end{cases}
\end{equation}
and $p=\rho(\vec p)$ as in proposition~\tu{\ref{p-1}}.
\end{defn}

The following are corollaries of theorem~\ref{u-10}.

\begin{cor}[$\ch$]
\label{c-8}
Let $\vec a,\vec b\in\codep\oone$. 
Suppose $\vec a\bigexta\vec b\in\qcode\oone$ 
and $P(\vec a\bigexta\vec b)$ adds no new reals. 
Then $e(\vec a,\vec b)$ is a generic embedding. 
\end{cor}

\begin{cor}[$\ch$]
\label{c-11}
Let $\vec a_\gamma\in\codep\oone$ \tu($\gamma<\delta$\tu).
Suppose that $\rho(\vec{\vec a})\in\qcode\oone$ 
and $P(\rho(\vec{\vec a}))$ adds no new reals. 
Then $f(\vec{\vec a},X)$ is a generic embedding for all
$X\subseteq\delta$. 
\end{cor}

\begin{remark}
\label{r-2}
The argument in example~\ref{x-9} applies so that the antisymmetric
quotient of $P(\rho(\vec{\vec a}))$ is a complete semilattice and the
map from $P(\rho(\vec{\vec a}\restriction X))\div{\qasym}$ into
$P(\rho(\vec{\vec a}))\div{\qasym}$ induced by $f(\vec{\vec a},X)$ has
an upward order closed range. Hence lemma~\ref{p-22} justifies the
following definition.
\end{remark}

\begin{defn}
\label{d-22}
Let  $\pi(\vec a,\vec b):P(\vec a\bigexta \vec b)\to P(\vec b)$ 
be the projection defined in equation~\eqref{eq:58} 
from $e(\vec a,\vec b)$; 
and let $\nu(\vec{\vec a},X):P(\rho(\vec{\vec a}))\to
P(\rho(\vec{\vec a}\restriction X))$ be the projection defined in
equation~\eqref{eq:58} from $f(\vec{\vec a},X)$. 
\end{defn}

\begin{prop}
\label{p-6}
$\pi(\vec a,\vec b)$ is a left inverse of $e(\vec a,\vec b)$. 
\end{prop}

\begin{prop}
\label{p-7}
$\nu(\vec{\vec a},X)$ is a left inverse of $f(\vec{\vec a},X)$. 
\end{prop}

The important properties of $\pi$, and more generally of  $\nu$ are:

\begin{lem}[$\ch$]
\label{p-57}
Let $\vec a,\vec b\in\codep\oone$.  Suppose $\vec a\bigexta\vec
b\in\qcode\oone$ and $P(\vec a\bigexta \vec b)$ adds no new reals. 
If $q\in\genc(M,P(\vec a\bigexta\vec b))$ 
then $\pi(\vec a,\vec b)(q)\in\genc(M,P(\vec b))$. 
\end{lem}
\begin{pf}
Proposition~\ref{p-34}, corollary~\ref{c-8} and proposition~\ref{p-6}.
\end{pf}

\begin{lem}[$\ch$]
\label{c-2}
Suppose that $\vec{\vec a}$ is a sequence of members of $\codep\oone$,
with $\rho(\vec{\vec a})\in\qcode\oone$, and $P(\rho(\vec{\vec a}))$
adds no new reals.
If $q\in\genc(M,P(\rho(\vec{\vec a})))$, 
then for all $X\subseteq\delta$, $\rho(\vec
q)\in\genc(M,P(\rho(\vec{\vec a}\restriction X)))$
where
\begin{equation}
  \label{eq:6}
  q_\gamma=\pi\bigl(\rho(\vec{\vec a}\restriction\gamma),\vec a_\gamma\bigr)
  \bigl(q\restriction\length(\rho(\vec {\vec a}\restriction
  \gamma)\bigext\vec a_\gamma)\bigr)
\end{equation}
for all $\gamma\in X$. 
\end{lem}
\begin{pf}
By proposition~\ref{p-34}, corollary~\ref{c-11} and
proposition~\ref{p-7},
$\nu(\vec{\vec a},X)(q)\in\genc(M,\allowbreak P(\rho(\vec{\vec a}\restriction X)))$. 
Equation~\eqref{eq:6} is established by verifying that
\begin{equation}
  \label{eq:10}
  \pi\bigl(\rho(\vec{\vec a}\restriction\gamma),\vec a_\gamma\bigr)
  \bigl(q\restriction\length(\rho(\vec {\vec a}\restriction
  \gamma)\bigext\vec a_\gamma)\bigr)\qsep p_\gamma
  \espc\text{for all $\gamma\in X$},
\end{equation}
where $\nu(\vec{\vec a},X)(q)=\rho(\vec p)$ and $\vec
p=(p_\gamma:\gamma\in X)$. 
\end{pf}

\fxwarning{Some interesting remarks, inside comment}

\subsection{$\mathrm{trind}$-properness}
\label{sec:trind-properness}

In~\cite[Ch.~XVIII, Definition~2.1]{MR1623206}, the notation
$\trind_\alpha(t)$ is used to denote all labelings
$\vec\beta=(\beta_x:x\in t)$ of some finite tree
$t$ with ordinals at most $\alpha$, 
i.e.~$\beta_x\le\alpha$ for all $x\in t$, so that
\begin{equation}
  \label{eq:39}
  x <_t y\impls\beta_x\le\beta_y.
\end{equation}
An operation is defined on iterations $\vec P=(P_\xi,\dot Q_\xi:\xi<\alpha)$ of
length $\alpha$ by members of $\vec\beta\in \trind_\alpha(t)$, 
where $\vec P_{\vec\beta}$  is the collection of all sequences $(p_x:x\in t)$
such that
\begin{samepage}
\begin{enumerate}[label=(\roman*), ref=\roman*, widest=ii]
\item $p_x\in P_{\beta_x}$ for all $x\in t$,
\item $x <_t y$ implies $p_y\restriction\beta_x=p_x$.
\end{enumerate}
\end{samepage}
Thus for example, if $t$ 
is a finite tree of height $1$ then for every
$\vec\beta\in\trind_\alpha(t)$,
 $\vec P_{\vec\beta}$ is
a finite product of the form $P_{\beta_0}\times \cdots \times P_{\beta_{n-1}}$
where $\beta_i\le\alpha$ for all $i=0,\dots,n-1$. 

Then (in~\cite[Ch.~XVIII, Definition~2.2]{MR1623206}) the notion of a
\emph{$\nnr_2$-iteration} is defined, which in particular entails that
the iteration is completely proper. Then the new theorem for
iterations not adding new reals is~\cite[Ch.~XVIII, Main
Lemma~2.8]{MR1623206} stating that if $(P_\xi:\xi\le\delta)$ an
iteration where $P_\xi$ is $\nnr_2$ for all $\xi$ less than the limit
$\delta$, then $P_\delta$ is $\nnr_2$. Without reviewing the details
of the definition of $\nnr_2$, we refer to this theory as the
\emph{$\trind$-properness} $\nnr$ theory. 

Unexpectedly, in overcoming the difficulties in constructing a properness parameter
suitable for forcing $\pstarc$, we came very close to satisfying the
hypotheses for the $\trind$-properness $\nnr$ theory. Indeed, using
the methods we have already presented,
our theorem~\ref{u-10} can be extended to say: $P(\vec a)_{\vec\beta}$
is proper for all $\vec a\in\gcodep\oone$ and all
$\vec\beta\in\trind_\alpha(t)$ for every finite tree $t$. Thus our
iteration, which will be of the form $P(\vec a)$ for some $\vec
a\in\gcodep\oone$, is \emph{$\trind$-proper}, i.e.~it remains proper after
operating on it with members of $\trind_\alpha(t)$. 

We think it is most likely that Shelah's above mentioned theorem can
be strengthened to something like: if $(P_\xi,\dot Q_\xi:\xi<\delta)$
is a countable support iteration such that 
$P_\xi\forces\ulc\dot Q_\xi$ is $\mathbb D\text{-complete}\urc$ for all $\xi<\delta$, 
and $(P_\xi,\dot Q_\xi:\xi<\delta)$ is $\trind$-proper, then $P_\delta$ adds no new
reals (probably this would require a slightly more general operation
than $\trind$). This seems to agree with his description of the essence
of the theory in~\cite[page~868]{MR1623206}; however, at present we do
not have a good enough understanding of his proof to make a
conjecture. 

Such a theorem would
result in a better (or at least shorter) proof of
theorem~\ref{u-3} than the one here using properness parameters. 
However, as it stands, the definition of $\vec P=(P_\xi,\dot Q_\xi:\xi<\delta)$ being
$\nnr_2$ requires the properness
of $\vec P'_{\vec\beta}$ for $\beta\in\trind_\alpha(t)$ where $\vec
P'$ is some arbitrary completely proper \emph{extension} of some
initial segment of~$\vec P$. Our iteration will not
satisfy this requirement of $\nnr_2$. 

\section{Model of $\ch$}
\label{sec:proof-cons-ppst}

\fxwarning{interesting comment}

\fxnote{comment: old attempt}

We begin with an arbitrary ground model $V$ (of enough of $\zfc$) satisfying $\gch$.
Set $\kappa=\aleph_2$ and $\lambda=\aleph_3$ as in equation~\eqref{eq:88}.

\fxnote{comment}

As usual, $\ns(\powcif A)$ denotes the nonstationary ideal on $\powcif
A$ and $\ns^*(\powcif A)$ is the dual filter, and thus is generated by the family
of closed cofinal subsets of $\powcif A$. 

\begin{defn}
\label{d-17}
Whenever $V\models\ulc\E\in\ns^*(\powcif{H_\kappa})\urc$,
let
$\ns^*(\E;V)=\bigl\{\F\subseteq\E:\F\cap\S\ne\emptyset$
for all $\S\in\bigl(\ns^*(\powcif{H_\kappa})\bigr)^V\bigr\}$. 
Let $\ns^*(V)$ denote $\ns^*((\powcif{H_\kappa})^V,V)$, when
$H_\kappa$ is understood. 
\end{defn}

\begin{prop}
\label{p-35}
Suppose $P$ is proper and $\dot\T$ is a $P$-name where
$P\forces\dot\T\subseteq\E[\dot G_P]$ \tu(cf.~notation~\tu{\ref{o-3}}\tu). The following are equivalent:
\begin{enumerate}[label=\tu{(\alph*)}, ref=\alph*, widest=b]
\item $P\forces\dot\T \in\ns^*(\E;V)$. 
\item For all $M\prec H_\lambda$ with $P,\E\in M$ and $M\cap H_\kappa\in\E$,
every $p\in P\cap M$ has 
an $(M,P)$-generic extension $q$ such that $q\forces M\cap H_\kappa\in\dot\T$. 
\end{enumerate}
\end{prop}

\fxnote{comment}

\subsection{The properness parameters}
\label{sec:properness-parameter}

The cardinal sequence
\begin{equation}
  \label{eq:5}
  \mu_\alpha=\aleph^+_{2+\alpha}\espc\text{($\alpha<\oone$)}
\end{equation}
is suitable for a $\lambda$-properness parameter, and we let $\vec\A$
be any fixed skeleton, e.g.~$\A_\alpha=\{M\in\powcif{H_{\mu_\alpha}}:
M\prec\nobreak H_{\mu_\alpha}\}$ ($\alpha<\oone$). 

We begin by motivating the definitions to follow. 
Let $\H$ be a $\sigma$-directed subfamily of
$(\powcif\oone,\subseteqfnt)$ with no stationary subset of $\oone$
orthogonal to it. In lemma~\ref{l-14} we saw that for a given map
$\Omega$ on $\injlim\A$, a sufficient condition for $\poset(\H)$ to be
$(\vec\A,\idealmodp\Omega)$-proper is that there exist
$y_M\in\Omega(M)\cap\downcl\H$ for all $M\in\injlim\A$ of positive
rank satisfying $\varphi_*(M,\H,y)$.
Conversely, suppose that $X\subseteq\genmod(q)\cap M$ for some
$q\in\gen(M,\poset(\H))$. Since $C_{\dot G_{\poset(\H)}}$ names a club
(cf.~corollary~\ref{c-1}), $q\forces\trsup_\oone(X)\subseteq C_{\dot
  G_{\poset(\H)}}$. 
Therefore, every initial segment of
$\trsup_\oone(X)$ must be in $\downcl\H$. 

Consider the next simplest case: an iteration of the form
$\spo(\I)\gdot\poset(\dot\H)$ where $\dot\H$ names an $\H$ as
above. In order that $\spo(\I)\gdot\poset(\dot\H)$ is
$\idealmodp\Omega$-proper, we must in particular have for every $M$ of positive
rank, that every finite sequence 
$(p_0,\dot s_0),\dots,\allowbreak(p_{n-1},\dot s_{n-1})
\in\spo(\I)\gdot\poset(\dot\H)\cap M$
has $(q_i,\dot t_i)\in\gen(M,\spo(\I)\gdot\poset(\dot\H),p_i)$ ($i=0,\dots,n-1$)
and an $X\in\idealmodp\Omega(M)$ such that
$X\subseteq\bigcap_{i=0}^{n-1}\genmod(q_i,\dot t_i)$. 
Let us focus on the case $n=2$. 
We shall need $y_0,y_1\in\Omega(M)$ such that 
$q_i\forces y_i\in\downcl{\dot\H}$ and $q_i\forces\varphi_*(M,\H,y_i)$ for $i=0,1$. 
Then to ensure that $\idealmodp\Omega(M)\ne\emptyset$ we would apply lemma~\ref{l-5}. 

In particular, to satisfy property~\eqref{item:38} in the definition
of instantiation, this means that we must be able to find
cofinally many $K\in\injlim\A\cap M$ with
\begin{equation}
\sup(\oone\cap K)\in y_0\cap y_1.\label{eq:12}
\end{equation}
This can be achieved as follows. Let $\vec a$ be the code for
$\spo(\I)$, let $\vec b$ be the code for
$\spo(\I)\gdot\poset(\dot\H)$. Assume that $\vec b\in\gcode\oone$.  
First of all we choose $r\in\genc(M,P(\vec a\bigexta\vec a),
\allowbreak{p_0}\bigext{p_1})$.
Then we let $q_0=\pi(\emptyset,\vec a)(r\restriction\length(\vec
a))\in P(\vec a)$
and $q_1=\pi(\vec a,\vec a\bigexta\vec a)(r)\in P(\vec a)$. It follows
from an application of lemma~\ref{c-2} that
${q_0}\bigexta q_1\in\genc(M,P(\vec a\bigexta\vec a))$.  
Then by extending $q_0$ and $q_1$ we may assume that there exist
$y_0$ and $y_1$ as above. Now for some fixed $b\in M\cap H_\lambda$
and $\xi<\rank(M)$, 
suppose that we want to find $K\in\A_{\xi}\cap M$ with $b\in K$
satisfying~\eqref{eq:12}. Since $\A_\xi\in M$ is stationary,
$S=\trsup(\A_\xi)\subseteq\oone$ is a stationary set in $M$. By the
assumption that $P(\vec a)$ forces that $\dot\H$ has no stationary
orthogonal set and by lemma~\ref{l-70}, $P(\vec a)\gdot\spo(\dot\H)$
forces that $S\cap X_{\ddot G_{\spo(\dot\H)}}$ is stationary. Since in
particular, $\vec b\in\qcode\oone$, 
we know that $P(\vec a)\gdot\spo(\dot\H)$ 
forces there is no stationary set orthogonal to $\dot\H$. 
It then follows from lemma~\ref{l-53} 
that $P(\vec a)\gdot\spo(\dot\H)$ forces that there is no stationary set
orthogonal to $e^*(\vec a,\vec a)(\dot\H)$. 
Hence applying lemma~\ref{l-70} again, 
$P(\vec a)\gdot\spo(\dot\H)\gdot\spo(e^*(\vec a,\vec a)(\dot\H))$
forces that 
\begin{equation}
\label{eq:4}
S\cap X_{\ddot G_{\spo(\dot\H)}}\cap X_{\ddot G_{\spo(e^*(\vec a,\vec a)(\dot\H))}}
\espc\text{is stationary}. 
\end{equation}
The whole point of invoking the embedding $e(\vec a,\vec a)$ is that
we want the name $\dot\H$ to be interpreted according to $q_1$ (in
particular, $e^*(\vec a,\vec a)(\dot\H)$ is independent of $q_0$,
unlike $\dot\H\cap M$ which is determined by $q_0$). 
It is now straightforward to produce $K\in\A_\xi\cap M$ with $b\in K$
satisfying~\eqref{eq:12}. 

The argument just outlined is a simplified version of our main lemma,
lemma~\ref{l-51}. The general case, were $\vec a$ codes an
initial segment of our iteration, is where we need to use $\gcode\oone$. 

\begin{defn}
\label{d-29}
Suppose $\vec a\in\code\theta$.
Let $\codeg{\vec a}$ be the set of all codes generated by
restriction and concatenations of arbitrary length, i.e.~$\vec b\in\codeg{\vec a}$ implies
that $b\restriction\xi\in\codeg{\vec a}$ for all $\xi<\length(\vec b)$, 
and $\vec b_\gamma\in\codeg{\vec a}$ ($\gamma<\delta$)
implies that $\vec c\in\codeg{\vec a}$, where
\begin{equation}
  \label{eq:67}
  \vec c=\vec b_0\bigexta\cdots\bigext \vec b_\gamma
  \bigexta\cdots\espc\text{($\gamma<\delta$)}.
\end{equation}
For $\kappa$ an infinite cardinal, 
we let $\codegfin{\vec a}\kappa$ be the subfamily of $\codeg{\vec a}$
generated by restrictions, and concatenations of length less than $\kappa$.
\end{defn}

\begin{prop}
\label{p-48}
Let $\kappa$ be an infinite regular cardinal. Then
$\codegfin{\vec a}{\kappa}$ consists of all codes of the form 
$(\vec a\restriction\xi_0)\bigexta\cdots\bigext(\vec
a\restriction\xi_{\gamma})\bigexta\cdots$ \tu($\gamma<\delta$\tu)  
where each\/ $\xi_\gamma<\length(\vec a)$ and\/ $\delta<\kappa$.
\end{prop}

\begin{notn}
\label{o-1}
For each $\vec a\in\code\theta$, and each ordinal $\gamma$, 
we let $\vec a^\gamma$ denote the concatenation $\vec a\bigexta\vec
a\bigexta\cdots$ iterated $\gamma$ times,
i.e.~$\length(\vec
a^\gamma)=\length(\vec a)\cdot\gamma$ 
and  $\vec a^\gamma(\length(\vec a)\cdot \zeta+\rho)=\vec a(\rho)$ 
for all $\zeta<\gamma$ and $\rho<\length(\vec a)$ 
(and considering $P(\vec a)$-names to also be $P(\vec a\bigexta \vec b)$-names). 
\end{notn}

\begin{defn}
\label{d-23}
Define a (class) function
$\psi=\psi_\theta:\code\theta\times\ord\to\code\theta$
by recursion on $\length(\vec a)$ by $\psi(\<\>,\gamma)=\<\>$, and 
\begin{equation}
  \label{eq:65}
  \begin{split}
  \psi(\vec a,\gamma)
  &:=\bigcup_{\xi<\length(\vec a)}\psi(\vec a\restriction\xi)\\
  &\hphantom{:}=(\vec a\restriction 1)^\gamma\bigexta(\vec a\restriction
  2)^\gamma\bigexta
  \cdots\bigext(\vec
  a\restriction\xi)^\gamma\bigexta\cdots\espc\text{($\xi<\length(\vec  a)$)}.
  \end{split}
\end{equation}
\end{defn}

\begin{prop}
\label{p-49}
$\psi_\theta(\vec a,\gamma)
\in\codegfin{\vec a}{\max\{\length(\vec a),|\gamma|\}^+}$ 
for all\/~$\vec a\in\code\theta$. 
\end{prop}

Henceforth, $\theta=\oone$. 

\begin{defn}
\label{d-27}
We let 
\begin{multline*}
\label{eq:53}
\Phi\bigl(\vec a,
(r^M,q^M_{\xi\gamma},y^M_{\xi\gamma}:M\in\injlim\A\text{, }\vec a\in M\text{, }
\xi<\length(\vec a)\text{ and $\gamma<\oone$ are in }M)\bigr)
\end{multline*}
be a formula expressing the following state of affairs:
 $\vec a\in\gcode\oone$;
and for all $M\in\injlim\A$ with $\vec a\in M$,
\begin{enumerate}[leftmargin=*, label=(\roman*), ref=\roman*, widest=iii]
\item\label{item:3} $r^M\in\genc(M,P(\psi(\vec a,\oone)))$,
\save
\end{enumerate}
and for all $\xi<\length(\vec a)$ and all $\gamma<\oone$ with $\xi,\gamma\in M$,
\begin{enumerate}[leftmargin=*, label=(\roman*), ref=\roman*, widest=iii]
\restore
\item\label{item:125} $q^M_{\xi\gamma}\in P(\vec a\restriction\xi)$,
\item\label{item:130} 
$q^M_{\xi\gamma}\ge \pi\bigl(\vec a\restriction\xi,
\bigcup_{\zeta<\xi}\psi(\vec a\restriction\zeta)
\bigext(\vec a\restriction\xi)^\gamma\bigr)
\bigl(r^M\restriction\length\bigl(
\bigcup_{\zeta<\xi}\psi(\vec a\restriction\zeta)
\bigext(\vec a\restriction\xi)^\gamma\bigr)\bigr)$,
\item\label{item:122} $q^M_{\xi\gamma}\forces
  y^M_{\xi\gamma}\in\dot\H(\vec a(\xi))$,
\item\label{item:123} $q^M_{\xi\gamma}\forces x\subseteqfnt y^M_{\xi\gamma}$ 
  for all $x\in\dot\H(\vec a(\xi))\cap M$.
\end{enumerate}
We abbreviate the above expression as 
$\Phi\bigl(\vec a,(\vec r,\vec q,\vec y)\bigr)$. 

Assuming $\Phi\bigl(\vec a,(\vec r,\vec q,\vec y)\bigr)$, 
for each $\vec c\in\codegfin{\vec a}{\aleph_0}$,
say $\vec c=(\vec a\restriction \xi_0)\bigexta\cdots\bigext(\vec a\restriction\xi_{k-1})$
(cf.~proposition~\ref{p-48}), and each $\alpha<\oone$, 
define a $P(\vec c)$-name 
\begin{multline}
  \dot\B^{\vec c}_\alpha=\bigl\{M\in\A_\alpha:
  q^M_{\xi_0 \gamma_0}{}\bigext\cdots\bigext q^M_{\xi_{k-1}\gamma_{k-1}}
  \in\dot G_{P(\vec c)}\\
  \text{ for some }\gamma_0<\cdots<\gamma_{k-1}\text{ in $\oone\cap M$}\bigr\}.
\end{multline}
\fxnote{Should specify exactly how $\B^{b}$ is to be used.}
We also define 
$\Omega(\vec y)(M)\in\powcnt{\powcif\theta}$  by
\begin{equation}
  \label{eq:119}
  \Omega(\vec y)(M)=
  \bigl\{y^M_{\xi\gamma}:\xi<\length(\vec a),
  \gamma<\oone\text{, }\xi,\gamma\in M\bigr\},
\end{equation}
and we put
\begin{equation}
\Z(\vec a)=\bigl\{(P(\vec c),\dot\B^{\vec c}):
\vec c\in C(\vec a,\aleph_0)\bigr\}.\label{eq:120}
\end{equation}
\end{defn}

\begin{lem}
\label{p-46}
For all $\vec c=(\vec a\restriction\xi_0)\bigexta\cdots\bigext(\vec a\restriction\xi_{k-1})$ in $C(\vec a,\aleph_0)$,
\begin{equation}
q^M_{\xi_0\gamma_0}{}\bigext\cdots\bigext
q^M_{\xi_{k-1}\gamma_{k-1}}\in\genc(M,P(\vec c))\label{eq:121}
\end{equation}
for all $\gamma_0<\cdots<\gamma_{k-1}$ in $\oone\cap M$. 
\end{lem}
\begin{pf}
This is a straightforward application of lemma~\ref{c-2}. 
\end{pf}

The reason that the codes are repeated $\oone$ times (rather than just
$\omega$ times) is so that we have the following.

\begin{lem}
\label{p-47}
For all $\vec c=(\vec a\restriction\xi_0)\bigexta\cdots\bigext(\vec a\restriction\xi_{k-1})$ 
in $C(\vec a,\aleph_0)$, every $p\in P(\vec c)\cap M$ has
$\gamma_0<\cdots<\gamma_{k-1}$ in $\oone\cap M$ such that
\begin{equation}
  \label{eq:122}
  q^M_{\xi_0\gamma_0}{}\bigext\cdots\bigext
  q^M_{\xi_{k-1}\gamma_{k-1}}
  \ge p.
\end{equation}
\end{lem}
\begin{pf}
Standard density argument since we have countable supports with an
iteration of uncountable cofinality. 
\end{pf}

\begin{lem}
\label{l-31}
$P(\vec c)\forces\dot\B^{\vec c}_\alpha\in\ns^*(\A_\alpha, V)$ for all
$\alpha<\oone$.
\end{lem}
\begin{pf}
We apply proposition~\ref{p-35}. Find $N\prec H_{\mu_{\alpha+1}}$ with
$\A_\alpha\in N$ and $M=\utilde N\in\A_\alpha$.
Take $p\in P(\vec c)\cap M$. Then
$q^M_{\xi_0\gamma_0}\bigexta\cdots\bigext
q^M_{\xi_{k-1}\gamma_{k-1}}\ge p$ for some
$\gamma_0<\cdots<\gamma_{k-1}$ in $\oone\cap M$ by lemma~\ref{p-47}. 
Then $q^M_{\xi_0\gamma_0}\bigexta\cdots\bigext
q^M_{\xi_{k-1}\gamma_{k-1}}\forces M\in\dot\B^{\vec c}_\alpha$ as wanted. 
\end{pf}

\fxnote{comment}

\begin{notn}
\label{o-4}
For an iterated forcing notion of the form
 $R=P_0\gdot \dot Q_0\gdot\dot Q_1\gdot\cdots\gdot\dot Q_n$, 
a $R$-name $\dot A$ and $r=(p,\dot q(0),\dots,\dot
q(n))\in\genc(M,R)$, we let $\dot A[p,\dot q(0),\dots,\dot q(n)]$
denote the interpretation of $\dot A$ by $\dot G_{R}[M,r]$ 
\tu(cf.~\tu{\Section\ref{sec:terminology})}.
\end{notn}

\begin{lem}
\label{l-51}
$\Phi\bigl(\vec a,\Dot{\vec\H},(\vec r,\vec q,\vec y)\bigr)$ 
implies that for all $M\in\injlim\A$ with $\rank(M)>0$, 
for all $\vec c\in C(\vec a,\aleph_0)\cap M$, say as in~\eqref{eq:67},
for all $\gamma_0<\cdots<\gamma_{k-1}$ in $\oone\cap M$,
for all $b\in M\cap H_\lambda$, for all $\alpha<\rank(M)$
there exists $K\in\A_\alpha\cap M$ such that
\begin{enumerate}[label=\tu{(\alph*)}, ref=\alph*, widest=b]
\item\label{item:112} $b\in K$,
\item\label{item:113} $\sup(\theta\cap K)\in\bigcap_{i=0}^{k-1}
  y^M_{\xi_i\gamma_i}$,
\item\label{item:114} $K\in\dot\B^{\vec b}_\alpha\bigl[q^M_{\xi_0\gamma_0},\dots,
  q^M_{\xi_{k-1}\gamma_{k-1}}\bigr]$.
\end{enumerate}
\end{lem}
\begin{pf}
Working in 
$M[\dot G_{P(\vec c)}[q_{\xi_0\gamma_0}^M,\dots,q_{\xi_{k-1}\gamma_{k-1}}^M]]$: 
Lemma~\ref{l-31} in particular implies that
$\C=\bigl\{K\in\dot\B^{\vec
  c}_\alpha[q_{\xi_0\gamma_0}^M,\dots,q_{\xi_{k-1}\gamma_{k-1}}^M]:b\in K\bigr\}$
is a cofinal subset of $\powcif{H_{\mu_\alpha}}$.
Let $S=\trsup_\theta(\C)$, which is thus stationary. 
We define $\dot S_{n}$ and $\vec d_n\in D(\<\>)$ 
by recursion on $n=0,\dots,k$ so that $\dot
S_0=S$, $\vec d_0=\<\>$ and
\begin{enumeq}
\item\label{item:100}
 $\dot S_{n+1}$ is 
 a $P(\vec c)\gdot\spo\bigl(e^*(\vec d_0,\vec a\restriction\xi_0)(\dot\H(\vec a(\xi_0)))\bigr)
 \gdot
 \cdots\gdot\spo\bigl(e^*(\vec d_n,\vec a\restriction
 \xi_n)(\dot\H(\vec a(\xi_n)))\bigr)$-name 
 for a stationary subset of $\dot S_{n}$ locally in 
  $e^*(\vec d_n,\vec a\restriction \xi_n)(\dot\H(\vec a(\xi_n)))$,
\item $\vec d_{n+1}=\vec d_n\bigexta\bigl(e^*(\vec d_n,\vec a\restriction\xi_n)
(\dot\H(\vec a(\xi_n))),\spo\bigr)$.
\end{enumeq}
This possible by lemma~\ref{l-70}, by the hypothesis that $\vec a\in\gcode\oone$, and  
thus forcing notions as in~\eqref{item:100} do not add stationary
subsets of $\oone$ orthogonal to any $\dot\H(\vec a(\xi))$,
and therefore do not add stationary subsets orthogonal to any 
$e^*(\vec d_n,\allowbreak\vec a\restriction\xi_n)(\dot\H(\vec a(\xi_n)))$
by lemma~\ref{l-53}.

We can find (an infinite) $x\in\powcif\theta$ and
$\vec p\in\spo\bigl(e^*(\vec d_0,\vec a\restriction\xi_0)(\dot\H(\vec a(\xi_0)))\bigr)\gdot
\cdots\gdot\spo\bigl(e^*(\vec d_{k-1},\vec
a\restriction\xi_{k-1})(\dot\H(\vec a(\xi_{k-1})))\bigr)$ so that $\vec p\forces x\subseteq\dot S_k$. 
Now by equation~\eqref{item:100}, 
$x\in e^*(\vec d_n,\vec a\restriction\xi_0)(\dot\H(\vec a(\xi_n)))
[e(\vec d_n,\vec a\restriction\xi_0)(q_{\xi_n\gamma_n}^M)]
=\dot\H(\vec a(\xi_n))[q_{\xi_n\gamma_n}^M]$
for all $n=0,\dots,k-1$. Thus, as $x\in M$ by complete properness,
$x\subseteqfnt y^M_{\xi_n\gamma_n}$ by equation~\eqref{item:123}, for all $n$.
Hence there exists $\delta\in x\cap y^M_{\xi_0\gamma_0}\cap\cdots
\cap y^M_{\xi_{k-1}\gamma_{k-1}}$.   
And then by elementarity, there exists $K\in\C\cap M$ 
with $\sup(\theta\cap K)=\delta$. 
\end{pf}

\begin{cor}
\label{c-13}
$\Phi\bigl(\vec a,\Dot{\vec\H},(\vec r,\vec q,\vec
y)\bigr)$    implies that $\idealmodp{\Omega(\vec y)}(\A;\Z(\vec a))$ 
is a properness parameter. 
\end{cor}
\begin{pf}
We apply lemma~\ref{l-61}. Let $M\in\injlim\A$ with $\rank(M)>0$ be given.
Each $(P,\dot\B)\in\Z\cap M$ is of the form $(P(\vec c),\dot\B^{\vec c})$ 
for some $\vec c\in C(\vec a)\cap M$, 
say $\vec c=(\vec a\restriction\nobreak\xi_0)
\bigext\cdots\bigext(\vec a\restriction \xi_{n^{\vec c}-1})$.
Using lemmas~\ref{p-46} and~\ref{p-47}, 
we can find pairwise disjoint sequences $\vec\gamma^{\vec c}_p
\in\omega_1^{n^{\vec c}}\cap M$ 
($p\in P(\vec c)\cap M$). We can also arrange that (the ranges of)
$\vec\gamma^{\vec c}_p$ and $\vec\gamma^{\vec c\,'}_{p'}$ are disjoint
whenever $\vec c\ne\vec c\,'$. Define 
\begin{equation}
\bar q^M_{P(\vec c),\dot\B^{\vec c}}(p)
={q^M_{\xi_0\gamma_p^{\vec c}(0)}}{}\bigext\cdots\bigext
{q^M_{\xi_{n^{\vec c}-1}\gamma_p^{\vec c}(n^{\vec c}-1)}}
\label{eq:9}
\end{equation}
for each $\vec c\in C(\vec a)\cap M$ and $p\in P(\vec c)\cap M$. 

To apply lemma~\ref{l-61},
let $A\subseteq\Omega(\vec y)$ be finite, 
say $A=\{y^M_{\xi_0\gamma_0},\dots,y^M_{\xi_{k-1}\gamma_{k-1}}\}$, 
$\vec c_0,\dots,\vec c_{m-1}$ be codes for members of $\Z\cap M$, 
let $O_i\subseteq P_i\cap M$ be finite for each $i=0,\dots,m-1$, 
let $b\in M\cap H_\lambda$ and $\xi<\rank(M)$. 
By extending both $A$ and the subset of $\Z\cap M$, 
we may assume without loss of generality that
$\{(\xi_0,\gamma_0),\dots,\allowbreak(\xi_{k-1},\gamma_{k-1})\}
=\bigcup_{i=0}^{m-1}\bigcup_{p\in O_i}
\{(\xi_0,\gamma_p^{\vec c_i}(0)),\dots,
(\xi_{n^{\vec c}-1},\gamma_p^{\vec c_i}(n^{\vec c_i}-1))\}$. 
Then an application of lemma~\ref{l-51} yields 
 $K\in\A_\xi\cap M$ with $b\in K$,
$\sup(\theta\cap K)\in\bigcap_{i=0}^{k-1}y^M_{\xi_i\gamma_i}=\bigcap A$
and  $K\in\dot\B^{{\vec c_0}\bigext\cdots\bigext(\vec c_{m-1})}_\xi
[q^M_{\xi_0\gamma_0},\dots,q^M_{\xi_{k-1}\gamma_{k-1}}]$. 
It follows that 
$\bar q^M_{P(\vec c_i),\dot\B^{\vec c_i}}(p)\forces K\in\dot\B_\xi^{\vec c_i}$ 
for all $i=0,\dots,m-1$ and all $p\in O_i$. 
We have therefore found $K$
witnessing~\eqref{item:68},~\eqref{item:76},~\eqref{item:83} 
and~\eqref{item:132} of lemma~\ref{l-61}.
Moreover, conditions~\eqref{item:75} and~\eqref{item:1}
automatically follow from the definitions of $\Omega(\vec y)$ and
$\dot\B^{\vec c_i}_\xi$. 
\end{pf}

\begin{cor}
\label{l-69}
$\Phi\bigl(\vec a,(\vec r,\vec q,\vec y)\bigr)$ implies
that $P(\vec c)$ is $\idealmodp{\Omega(\vec y)}(\vec\A;\Z(\vec a))$-proper for
all $\vec c\in C(\vec a)$.
\end{cor}
\begin{pf}
By corollary~\ref{c-13}, corollary~\ref{p-29},
equation~\eqref{eq:120},
the definition of $\dot\B_\alpha^{\vec c}$ and
proposition~\ref{p-46}. 
\end{pf}

\begin{remark}
\label{r-14}
What we actually need (see the proof of lemma~\ref{l-58}), is that $P(\vec c)$ is
\emph{long} $\idealmodp{\Omega(\vec y)}(\vec\A;\Z(\vec a))$-proper. 
This can be proved using the ideas already presented.
\end{remark}

The following says that $\Phi$ is ``preserved'' at successors.

\begin{lem}
\label{l-32}
Assume $\Phi\bigl(\vec a,(\vec r,\vec q,\vec y)\bigr)$.
If $\vec a\bigexta(\dot\H,\poset)\in\gcode\oone$,
then there exists $(\vec r_*,\vec q_*,\allowbreak\vec y_*)$ such that
$\Phi\bigl(\vec a,(\vec r_*,\vec q_*,\vec y_*)\bigr)$ holds.
\end{lem}
\begin{pf}[Sketch of proof]
Set $\vec b=\vec a\bigexta(\dot\H,\poset)$. 
By corollary~\ref{l-69} and theorem~\ref{u-1}, we can find
$r^M_p\in\genc(M,P(\psi(\vec a,\oone)),p)$ for all $M\in\injlim\A$ with
$\vec a\in M$, and all $p\in P(\psi(\vec a,\oone))\cap M$. 
Then for each $p$, we can find $q^M_p\ge r^M_p$ and $y^M_p\in\powcif\oone$ 
such that $q^M_p\forces y^M_p\in\dot\H$ and $x\subseteqfnt y^M_p$ for
all $x\in\dot\H\cap M$. For each $\alpha<\oone$, 
define a $P(\psi(\vec a,\oone))$-name
\begin{equation}
  \label{eq:69}
  \dot\C_\alpha=\bigl\{M\in\A_\alpha:q^M_p\in\dot G_{P(\psi(\vec a,\oone))}
\text{ for some $p\in P(\psi(\vec a,\oone))$}\bigr\}.
\end{equation}

Let $G\in\Gen(V,P(\psi(\vec a,\oone)))$.  
It is easy to see that $\dot\C_\alpha[G]$ is stationary for all
$\alpha<\oone$. Then defining $\Omega(M)\in\powcnt{\powcif\oone}$ 
by $\Omega(M)=\Omega(\vec y)(M)\cup\{y^M_p:p\in P(\psi(\vec
a,\oone))\cap M\}$, $\poset(\dot\H[G])$ is
$(\vec\C[G],\idealmodp\Omega)$-proper by lemma~\ref{l-14}. This proves
that $P(\vec b)$ is $(\vec\A,\idealmodp\Omega)$-proper.

Now this allows us to use the parameterized properness theory to find 
$r^M_*\in\genc(M,\psi(\vec b,\oone))$ for all $M$. It is then clear
how to find $\vec q_*$ and $\vec y_*$ so that $\Phi\bigl(\vec
a,(\vec r_*,\vec q_*,\vec y_*)\bigr)$ holds. 
\end{pf}

\fxnote{maybe need?}

The following says that $\Phi$ is ``preserved'' at limits. 

\begin{lem}
\label{l-58}
Let $\vec a\in\code\theta$. 
If for all $\xi<\length(\vec a)$, 
there exists $(\vec r_\xi,\vec q_\xi,\vec y_\xi)$ 
satisfying $\Phi\bigl(\vec a\restriction\xi,(\vec
r_\xi,\vec q_\xi,\vec y_\xi)\bigr)$, then there exists $(\vec r,\vec
q,\vec y)$ satisfying $\Phi\bigl(\vec a,(\vec r,\vec
q,\vec y)\bigr)$. 
\end{lem}
\begin{pf}
This is a straightforward application of  lemma~\ref{u-7}. 
\end{pf}

\fxnote{comment}

\fxnote{comment}

\fxnote{to be removed}

\fxnote{to be removed}

\fxnote{tbr}

\fxnote{comment}

\fxnote{old, remove}

\fxnote{older version, to be removed}

\fxnote{old stuff. To be removed}

\fxnote{comment}

\begin{pf}[of theorem~\tu{\ref{u-3}}]
We are going to recursively define an iterated forcing construct
$(P_\xi,\dot Q_\xi:\xi<\omega_2)$ of length $\omega_2$ with countable
supports, and let $P_{\omega_2}$ denote the limit of the iteration.
At the same time, we are going to
choose $\vec a_\xi\in\gcode\oone$ such that
\begin{enumerate}[label=(\roman*), ref=\roman*]
\item\label{item:128} $\length(\vec a_\xi)<\omega_2$,
\item\label{item:129} $P_\xi=P(\vec a_\xi)$,
\item\label{item:137} $a_\xi\subseteq a_\eta$ for all $\xi\le\eta$; 
\save
\end{enumerate}
we will also find $(\vec r_\xi,\vec q_\xi,\vec y_\xi)$ as in
definition~\ref{d-27}, so that
\begin{enumerate}[label=(\roman*), ref=\roman*]
\restore
\item\label{item:141} $\Phi\bigl(\vec a_\xi,(\vec r_\xi,\vec
  q_\xi,\vec y_\xi)\bigr)$ holds.  
\save
\end{enumerate}

Observe that from this information we can already deduce that
\begin{enumerate}[label=(\roman*), ref=\roman*]
\restore
\item\label{item:139} $P_\xi$ has the $\aleph_2$-cc for all $\xi\le\omega_2$, 
\item\label{item:140} $P_\xi$ has a dense suborder of cardinality at
  most $\aleph_2$ for all $\xi\le\omega_2$, 
\item\label{item:138} $P_\xi\forces 2^{\aleph_1}=\aleph_2$ for all $\xi$,
\item\label{item:142} $P_\xi$ is completely proper for all $\xi\le\omega_2$. 
\save
\end{enumerate}
This is so because~\eqref{item:128} and~\eqref{item:129} imply that
$P_\xi$ is an iteration of length at most $\omega_2$,
where each iterand satisfies the properness
isomorphism condition by lemmas~\ref{l-13} and~\ref{l-9}; 
hence, we can conclude condition~\eqref{item:139}. Conditions~\eqref{item:140}
and~\eqref{item:138} are established simultaneously by induction as
usual: If $P_\xi\forces 2^{\aleph_1}=\aleph_2$, then $P_\xi\forces|\dot
Q_\xi|=|P(\vec a_{\xi+1})\div P(\vec a_\xi)|\le\aleph_2$, and
therefore by the $\aleph_2$-cc, 
$P_{\xi+1}$ satisfies~\eqref{item:140} and~\eqref{item:138}.
Condition~\eqref{item:142} is of course by the parameterized
properness theory: By~\eqref{item:141} and corollaries~\ref{c-13}
and~\ref{l-69}, $\idealmodp{\Omega(\vec y_\xi)}(\vec\A;\Z(\vec
a_\xi))$ is a properness parameter for which $P_\xi$ is proper. 
Since
$P_\xi=P(\vec a_\xi)$ is an iteration with \deecmp iterands by
lemmas~\ref{l-8} and~\ref{l-39},
 $P_\xi$ adds no new reals by the $\nnr$ theorem (theorem~\ref{u-1}). 

Using conditions~\eqref{item:139}--\eqref{item:138}, 
by standard bookkeeping, and regarding
$P_\xi$-names as also being $P_\eta$-names for $\xi\le\eta$, 
we can arrange an enumeration $(\dot\H_\xi:\xi<\nobreak\omega_2)$ of
$P_\xi$-names in advance such that,
for every $\xi<\omega_2$ and every $P_\xi$-name $\dot\H$ for a
$\sigma$-directed subfamily of $\powcif\oone$, 
\begin{enumerate}[label=(\roman*), ref=\roman*]
\restore
\item\label{item:143} there exists $\aleph_2$ many  $\eta\ge\xi$ such that
$P_\eta\forces\dot \H_\eta=\dot\H$. 
\save
\end{enumerate}

Now we describe the construction. First we deal with the successor
stage $\xi+1$ of the construction. We separate into two cases:
\begin{description}
\item[Case 1] $\vec a_\xi\bigexta(\dot\H_\xi,\poset)\in\gcode\oone$.
\item[Case 2] $\vec a_\xi\bigexta(\dot\H_\xi,\poset)\notin\gcode\oone$.
\end{description}

In Case 1, we put $\vec a_{\xi+1}=\vec a_\xi\bigexta(\dot\H_\xi,\poset)$. 
Therefore,
\begin{equation}
  \label{eq:44}
  P_{\xi+1}\text{ forces that there exists a club locally in }\dot\H_\xi
\end{equation}
by lemma~\ref{l-16}.\fxnote{does lemma mention club?}  
And there exists $(\vec r_{\xi+1},\vec q_{\xi+1},\vec y_{\xi+1})$
satisfying $\Phi\bigl(\vec a_{\xi+1},(\vec r_{\xi+1},\vec
q_{\xi+1},\allowbreak\vec y_{\xi+1})\bigr)$ by lemma~\ref{l-32}. 

In Case 2, there exists $\vec c\in D(\vec a)$ and $\vec b\in C(\vec
c)$ with $\length(\vec b)<\omega_2$, 
with a condition $p\in P(\vec b)$ such that
\begin{equation}
  \label{eq:86}
  p\forces\text{there exists a stationary set orthogonal to }\dot\H_\xi.
\end{equation}
We set $\vec a_{\xi+1}=\vec b$. By corollary~\ref{l-69}, $P(\vec b)$
is $\idealmodp{\Omega(\vec y)}$-proper and thus we can 
take $(\vec r_{\xi+1},\vec q_{\xi+1},\vec y_{\xi+1})
=(\vec r_\xi,\vec q_\xi,\vec y_\xi)$. 

At limit stages $\delta$, we let $\vec
a_\delta=\bigcup_{\xi<\delta}\vec a_\xi$.
Then there exists $(\vec r_\delta,\vec q_\delta,\vec y_\delta)$
satisfying~\eqref{item:141} by lemma~\ref{l-58}.

Having completed the construction, let $G\in\Gen(V,P_{\omega_2})$.
Then $\aleph_1$ is not collapsed, i.e.~$\aleph_1^{V[G]}=\aleph_1$,
and $V[G]\models\ch$ by condition~\eqref{item:142}. 
Since $V\models\ch$, by the $\aleph_2$-cc and by
condition~\eqref{item:143}, 
every $\sigma$-directed family $\H$ of $(\powcif\oone,\subseteqfnt)$ is equal to
$\dot\H_\xi[G]$ for cofinally many $\xi<\omega_2$. Then assuming
standard bookkeeping, we can ensure that there exists $\xi<\omega_2$
such that $\H=\dot\H_\xi[G]$ and either equation~\eqref{eq:44} holds,
or else there exists $p\in G$ as in equation~\eqref{eq:86}. Therefore,
$V[G]\models\ulc\pstarc_\oone\urc$. 
\end{pf}

\bibliographystyle{elsart-harv}
\bibliography{database}

\begin{thebibliography}{13}
\expandafter\ifx\csname natexlab\endcsname\relax\def\natexlab#1{#1}\fi
\expandafter\ifx\csname url\endcsname\relax
  \def\url#1{\texttt{#1}}\fi
\expandafter\ifx\csname urlprefix\endcsname\relax\def\urlprefix{URL }\fi

\bibitem[{Abraham(2006)}]{hbst}
Abraham, U., 2006. Proper forcing. In: Foreman, Kanamori, Magidor (Eds.),
  Handbook of Set Theory. To appear.

\bibitem[{Abraham and Todor{\v{c}}evi{\'c}(1997)}]{MR1441232}
Abraham, U., Todor{\v{c}}evi{\'c}, S., 1997. Partition properties of
  {$\omega\sb 1$} compatible with {CH}. Fund. Math. 152~(2), 165--181.

\bibitem[{Hirschorn(2007{\natexlab{a}})}]{Hir-comb}
Hirschorn, J., 2007{\natexlab{a}}. Combinatorial and hybrid principles for
  $\sigma$-directed families of countable sets modulo finite,
  arXiv:0706.3729v1.

\bibitem[{Hirschorn(2007{\natexlab{b}})}]{H2}
Hirschorn, J., 2007{\natexlab{b}}. Random trees under $\ch$. Israel J. Math.
  157, 123--153.

\bibitem[{Ihoda and Shelah(1988)}]{MR973109}
Ihoda, J.~I., Shelah, S., 1988. Souslin forcing. J. Symbolic Logic 53~(4),
  1188--1207, {S}helah~[JdSh:292].

\bibitem[{Kunen(1980)}]{MR597342}
Kunen, K., 1980. Set theory. Vol. 102 of Studies in Logic and the Foundations
  of Mathematics. North-Holland Publishing Co., Amsterdam, {A}n introduction to
  independence proofs.

\bibitem[{Shelah(1984)}]{rst:S}
Shelah, S., 1984. Can you take {S}olovay's inaccessible away? Israel J. Math.
  48~(1), 1--47, {S}helah~[Sh:176].

\bibitem[{Shelah(1998)}]{MR1623206}
Shelah, S., 1998. Proper and improper forcing, 2nd Edition. Perspectives in
  Mathematical Logic. Springer-Verlag, Berlin.

\bibitem[{Shelah(2000{\natexlab{a}})}]{math.LO/0003115}
Shelah, S., 2000{\natexlab{a}}. {NNR Revisited}. arXiv:math.LO/0003115,
  {S}helah~[Sh:656].

\bibitem[{Shelah(2000{\natexlab{b}})}]{MR1804704}
Shelah, S., 2000{\natexlab{b}}. On what {I} do not understand (and have
  something to say). {I}. Fund. Math. 166~(1-2), 1--82, {S}aharon Shelah's
  anniversary issue, {S}helah~[Sh:666].

\bibitem[{Todor{\v{c}}evi{\'c}(1985)}]{MR792822}
Todor{\v{c}}evi{\'c}, S., 1985. Directed sets and cofinal types. Trans. Amer.
  Math. Soc. 290~(2), 711--723.

\bibitem[{Tukey(1940)}]{MR0002515}
Tukey, J.~W., 1940. Convergence and {U}niformity in {T}opology. Annals of
  Mathematics Studies, no. 2. Princeton University Press, Princeton, N. J.

\bibitem[{Woodin(1999)}]{MR1713438}
Woodin, W.~H., 1999. The axiom of determinacy, forcing axioms, and the
  nonstationary ideal. Vol.~1 of de Gruyter Series in Logic and its
  Applications. Walter de Gruyter \& Co., Berlin.

\end{thebibliography}

\end{document}